\documentclass[hidelinks,onefignum,onetabnum]{siamart250211}

\usepackage{graphicx}
\usepackage{epstopdf}
\ifpdf
  \DeclareGraphicsExtensions{.eps,.pdf,.png,.jpg}
\else
  \DeclareGraphicsExtensions{.eps}
\fi

\usepackage{amsmath,amssymb,amsfonts}
\usepackage{enumitem}
\usepackage{tikz}
\usepackage{bm}
\usepackage{algpseudocode}

\newcommand{\img}{\mathbf{i}}
\newcommand{\rme}{\mathrm{e}}
\newcommand{\rmd}{\mathrm{d}}

\newcommand{\sss}[1]{{\scriptscriptstyle(#1)}}

\newsiamremark{remark}{Remark}
\newsiamremark{hypothesis}{Hypothesis}
\crefname{hypothesis}{Hypothesis}{Hypotheses}
\newsiamthm{claim}{Claim}
\newsiamremark{fact}{Fact}
\crefname{fact}{Fact}{Facts}

\headers{Multi-Coefficient Inversion}{S. Lu, and B. Xu}

\title{A Recursive Algorithm for Multi-Coefficient Inversion in Nonlinear Helmholtz Equations\thanks{Submitted to the editors DATE.
\funding{The first author was supported by National Key Research and Development Programs of China (No.
 2023YFA1009103), Science and Technology Commission of Shanghai Municipality (23JC1400501). The second author was supported by the NSFC (12171301 and 12571459).}}}

\author{Shuai Lu\thanks{School of Mathematical Sciences, Fudan University, No. 220 Handan Road, Shanghai, China 
  (\email{slu@fudan.edu.cn}).}
\and Boxi Xu\thanks{Corresponding author. School of Mathematics, Shanghai University of Finance and Economics, No. 777 Guoding Road, Shanghai, China 
  (\email{xu.boxi@mail.sufe.edu.cn}).}
}

\usepackage{amsopn}

\ifpdf
\hypersetup{
  pdftitle={Multi-Coefficient Inversion},
  pdfauthor={S. Lu, and B. Xu}
}
\fi

\begin{document}

\maketitle

\begin{abstract}
We present a recursive algorithm for multi-coefficient inversion in nonlinear Helmholtz equations with polynomial-type nonlinearities, utilizing the linearized Dirichlet-to-Neumann map as measurement data.
To achieve effective recursive decoupling and simultaneous recovery of multiple coefficients, we develop a novel Fourier-based approach that combines the principle of inclusion-exclusion with carefully constructed complex exponential solutions.
This methodology not only ensures unique identifiability but also yields progressively increasing stability with enhanced wavenumbers.
Comprehensive numerical experiments demonstrate the algorithm's computational efficiency and excellent reconstruction accuracy.
\end{abstract}



\section{Introduction}

In this paper, we investigate the inverse problem of simultaneously recovering multiple coefficients in a nonlinear Helmholtz equation with polynomial-type nonlinearity:
\begin{align}\label{eqn:sec1}
\Delta u + k^{2} u = {\textstyle\sum\limits_{\ell=1}^{m}} c_{\ell}(x) \, u^{\ell} \quad \text{in\ } \Omega \subset \mathbb{R}^{n},
\end{align}
where $k > 0$ represents the wavenumber, the integer $m \geqslant 2$ characterizes the highest order of nonlinearity, $\{ c_{\ell}(x) : x \in \Omega,\ \ell = 1,2,\dots,m \}$ denote the unknown spatially-varying coefficients to be recovered, and  $\Omega \subset \mathbb{R}^{n}$ ($n \geqslant 2$) is an open bounded domain with smooth $C^{\infty}$ boundary $\partial\Omega$.
Nonlinear wave equations of this type describe wave propagation in a medium with intensity-dependent responses, where time-harmonic fields are governed by the nonlinear Schr\"{o}dinger or Helmholtz equations \cite{FT2005, YL2017, WZ2018}.
Such equations and their corresponding inverse problems play a crucial role in characterizing nonlinear optical phenomena induced by high-intensity radiation, with applications spanning nonlinear optics, photonics, and material science.
For further discussion on physical motivations and applications, we refer to \cite{SJ2004, BFT2009, XB2010, BPS2011}.

In the nonlinear Helmholtz equation \eqref{eqn:sec1}, the multiple coefficients $\{ c_{\ell}(x) \}_{\ell=1}^{m}$ associated with the nonlinear terms represent unknown spatially-varying functions that characterize the nonlinear responses within the domain $\Omega$.
The main inverse problem that we consider involves reconstructing these unknown coefficients $\{ c_{\ell}(x) \}_{\ell=1}^{m}$ from boundary measurements encoded in the linearized Dirichlet-to-Neumann (DtN) map, whose precise definition will be given subsequently.
To the best of our knowledge, this problem remains largely unexplored in the existing literature, with neither theoretical stability estimates nor effective numerical reconstruction algorithms currently available for this nonlinear inverse problem.
Under suitable regularity assumptions, we establish both uniqueness and increasing stability estimates for the multi-coefficient identification problem, demonstrating that the stability improves significantly with larger wavenumber $k$.
These theoretical results not only enable the development of an explicit recursive reconstruction algorithm but also suggest the feasibility of obtaining high-resolution reconstructions in practical applications.

The reconstruction of interior unknown coefficients from boundary measurements has remained a long-standing challenge in inverse problems research.
In what follows, we briefly review related developments pertaining to a single-coefficient identification problem for both linear and nonlinear Helmholtz equations.

The inverse problem of single-coefficient inversion for the linear Helmholtz equation (corresponding to $m = 1$ in \eqref{eqn:sec1}),
\begin{align}\label{eqn:single_c1}
\Delta u + k^{2} u = c_{1}(x) \, u \quad \text{in\ } \Omega,
\end{align}
has been extensively studied, yielding crucial insights into both theoretical stability and numerical reconstruction.
In the static case ($k = 0$), the problem reduces to electrical impedance tomography, where the seminal logarithmic stability results were first obtained in \cite{A1988} and later proven optimal in \cite{M2001}.
{Several hybrid methods have been proposed for conductivity identification, employing an augmented Lagrangian formulation to synergistically integrate output least-squares and equation error approaches.
These methods have proven effective in reconstructing both smooth \cite{IK1990, IKK1991} and non-smooth \cite{CZ1999, KZ2001} conductivity functions.
In the context of an inhomogeneous medium, the simultaneous inversion of the diffusion and absorption coefficients constitutes the basis of diffuse optical tomography.
A more comprehensive discussion can be found in \cite{AL1998, H2009, CHZ2021}.} 
For propagating waves ($k > 0$), the stability properties improve substantially, with \cite{I2011} first establishing increasing stability estimates that demonstrate enhanced resolution at higher wavenumbers.
This wavenumber-dependent stability framework has been successfully generalized to various linear inverse problems \cite{SI2007, NUW2013, I2015, CIL2016, LY2017, KUW2021}, providing a significant advantage over the elliptic case ($k = 0$) by enabling more robust reconstruction algorithms through a linearized DtN map.
A particularly numerical development is the Fourier-based approach in \cite{ILX2020} for equation \eqref{eqn:single_c1} with $k > 1$, which achieves stable recovery of coefficient $c_{1}(x)$ via its Fourier modes within a wavenumber-dependent frequency band.
The success of the linearization methodology has spurred applications to diverse inverse potential problems \cite{ILX2022, ZLX2022}, while recent work \cite{CFO2022} has further characterized the DtN map's local convexity near zero potential, establishing uniqueness guarantees.

Regarding the inverse problem of single-coefficient inversion for the nonlinear Helmholtz equation ($m \geqslant 2$) with power-type nonlinearity in \eqref{eqn:sec1},
\begin{align}\label{eqn:single_cm}
\Delta u + k^{2} u = c_{m}(x) \, u^{m} \quad \text{in\ } \Omega,
\end{align}
the analysis naturally divides into two main approaches: first-order linearization with respect to the unknown coefficient $c_{m}(x)$ and higher-order linearization with respect to the Dirichlet boundary data.
The first approach was originally developed for nonlinear parabolic inverse problems \cite{I1993} and subsequently extended to various nonlinear models \cite{IS1994, SU1997, S2010, SZ2012, LL2019}.
Recent work in \cite{LSX2022, ZLX2024} has combined this approach with the principle of inclusion-exclusion, enabling stable reconstruction of the single coefficient $c_{m}(x)$ in equation \eqref{eqn:single_cm} while establishing rigorous increasing stability estimates.
The higher-order linearization method, initially introduced for parameter identification in wave equations on Lorentzian manifolds in \cite{KLU2018}, was later adapted to inverse problems for equation \eqref{eqn:sec1} with $k = 0$ and positive integer $m$ in \cite{FO2020, FLL2023}.
This method constructs solutions corresponding to Dirichlet boundary data $u |_{\partial \Omega} = \sum_{j=1}^{m} \varepsilon_{j} f_{j}$ with small parameters $\varepsilon_{j} > 0$.
By computing the $m$-th order mixed derivative with respect to $\bm{\varepsilon} = (\varepsilon_{1},\varepsilon_{2},\dots,\varepsilon_{m})$ and evaluating at $\bm{\varepsilon} = 0$, one obtains a linearized equation for coefficient recovery.
This higher-order linearization technique has proven effective for coefficient reconstruction of $c_{m}(x)$ in \cite{LSX2022} and has been extended to various nonlinear inverse problems in \cite{CFKKU2021, LLLS2021, LLLS2021b, LLPT2022b, LLST2022}, including recent numerical implementations for $1$-D nonlinear wave equations using the Dirichlet-to-Neumann map in \cite{LLPT2023}.

Although much progress has been made in the investigation of single-coefficient identification problems for nonlinear equations, the polynomial-type nonlinearity in equation \eqref{eqn:sec1} presents substantial mathematical challenges, primarily due to the coupled interactions among coefficients $\{ c_{\ell}(x) \}_{\ell=1}^{m}$ through the polynomial structure.
Our main theoretical contribution establishes that the first-order linearization method can achieve stable recovery of these coefficients for general polynomial nonlinearities.
Through innovative combinatorial techniques grounded in the principle of inclusion-exclusion, we prove wavenumber-dependent stability estimates for degree-$m$ polynomial nonlinearities under first-order linearization, thereby addressing a crucial theoretical gap in the literature.
The cornerstone of our analysis is a combinatorial framework based on the principle of inclusion-exclusion that effectively decouples polynomial interactions and extends Alessandrini-type identities to arbitrary coefficients $\{ c_{\ell}(x) \}_{\ell=1}^{m}$.
These results not only maintain consistency with known stability properties from both first-order and higher-order linearization methods \cite{ILX2020, LSX2022, ZLX2024} for single-coefficient identification problems, but also uncover fundamental mathematical connections between these distinct approaches to nonlinear inverse problems.


The remainder of this paper is organized as follows.
Section \ref{sec:setup} presents the mathematical framework for the inverse problem, starting with the formulation of the first-order linearized Dirichlet-to-Neumann map.
In Section \ref{sec:scheme}, we derive a novel Alessandrini-type identity incorporating the principle of inclusion-exclusion, which forms the basis for recovering the unknown coefficients $\{ c_{\ell}(x) \}_{\ell=1}^{m}$ through systematic combinations of boundary measurements.
Section \ref{sec:theorem} contains our main theoretical contributions: \textbf{Theorem \ref{thm:Uniqueness}} establishes uniqueness for the multi-coefficient inversion; \textbf{Theorem \ref{thm:Lipschitz}} proves Lipschitz stability for Fourier modes $\mathcal{F}[c_{\ell}](\xi)$ within the frequency band $|\xi| \leqslant (\ell+1)k$; and \textbf{Theorem \ref{thm:main}} provides wavenumber-dependent increasing stability estimates.
Finally, Section \ref{sec:algorithm} presents a recursive reconstruction algorithm (\textbf{Algorithm \ref{alg:main}}) that implements these stability results computationally, accompanied by numerical examples demonstrating both the theoretical stability estimates and the algorithm's practical effectiveness.

\section{Problem setup}\label{sec:setup}

Let $m \geqslant 2$ be an integer throughout this paper.
Given a polynomial of degree $m$ in $u$, which is written in the form
\begin{align}\label{eqn:poly_Pm}
P_{m}(x,u) := {\textstyle\sum\limits_{\ell = 1}^{m}} c_{\ell}(x) \, u^{\ell},
\end{align}
where $\{ c_{\ell}(x) : \ell = 1,2,\dots,m \}$ denote the coefficient functions of the polynomial $P_{m}(x,u)$.
We consider the following nonlinear Helmholtz equation with Dirichlet boundary condition, that is,
\begin{align}\label{eqn:probl_I}
\text{($I$)}~
& \left\{~
\begin{aligned}
\Delta u + k^{2} u &= P_{m}(x,u)
& & \text{in\ } \Omega, \\
u &= f
& & \text{on\ } \partial\Omega. \\
\end{aligned}
\right.
\end{align}
Here, $k > 1$ is the wavenumber, $\Omega \subset \mathbb{R}^{n}$ is an open bounded domain, $n \geqslant 2$ is the number of dimensions, and $\partial\Omega$ is the boundary of domain $\Omega$.
We assumed that $\partial\Omega$ is smooth enough, i.e., $C^{\infty}$.

We further assume that each coefficient $c_{\ell} := c_{\ell}(x) \in L^{\infty}(\Omega)$ is compactly supported in $\Omega$, that is, $\mathop{\mathrm{supp}}(c_{\ell}) \subset \Omega$, and $k^{2}$ is not a Dirichlet eigenvalue of $-\Delta$ in $\Omega$.
Then, the well-posedness of the original problem \eqref{eqn:probl_I} with small Dirichlet boundary data $f$ and small coefficients $\{ c_{\ell} \}_{\ell=1}^{m}$ follows from the implicit function theorem \cite{RR2004}, and the solution $u \in W^{2,p}(\Omega)$, the Dirichlet boundary data $u |_{\partial\Omega} \in \mathcal{D} := W^{2-\frac{1}{p},p}(\partial\Omega)$, and the Neumann boundary data $\partial_{\nu} u |_{\partial\Omega} \in \mathcal{N} := W^{1-\frac{1}{p},p}(\partial\Omega)$, where $p := \frac{n}{2}+1$.
This well-posedness result is similar to that in \cite[Proposition 2.1]{ZLX2024} for a single nonlinearity term $c_{m} u^{m}$, and that in \cite{LLLS2021, N2023} for $k = 0$.


The multi-coefficient identification problem that we consider in current work is to recover the multiple coefficients $\{ c_{\ell} \}_{\ell=1}^{m}$ from the Dirichlet-to-Neumann (DtN) map on the boundary $\partial\Omega$.
For convenience, we denote $\bm{c} := (c_{1},c_{2},\dots,c_{m}) \in (L^{\infty}(\Omega))^{m}$.
The definition of the DtN map $\Lambda_{\bm{c}}$ of the original problem ($I$) in \eqref{eqn:probl_I} is given by
\begin{align}\label{eqn:DtN_I}
\Lambda_{\bm{c}} : u |_{\partial\Omega} \mapsto \partial_{\nu} u |_{\partial\Omega}.
\end{align}
Here, $u |_{\partial\Omega} \in \mathcal{D}$ and $\partial_{\nu} u |_{\partial\Omega} \in \mathcal{N}$ are the Dirichlet boundary data and the corresponding Neumann boundary data of the original problem ($I$) in \eqref{eqn:probl_I}, respectively.
While it is theoretically possible to establish the uniqueness of multi-coefficients by extending the results in \cite{LLLS2021} using the full DtN maps, a numerical approach necessitates an appropriate linearization strategy with respect to the DtN map $\Lambda_{\bm{c}}$.
Specifically, we will apply the first-order linearization method to the DtN map $\Lambda_{\bm{c}}$ associated with the nonlinear Helmholtz equation \eqref{eqn:probl_I}, as outlined below.

More precisely, we consider a linearization approach with respect to the coefficient functions of the polynomial $P_{m}(x,u)$, leading to the following system
\begin{align}
\text{($I^{\sss{0}}$)}~
& \left\{~
\begin{aligned}
\Delta u^{\sss{0}} + k^{2} u^{\sss{0}} &= 0
& & \text{in\ } \Omega, \\
u^{\sss{0}} &= f
& & \text{on\ } \partial\Omega, \\
\end{aligned}
\right. \label{eqn:probl_I0} \\[1ex]
\text{($I^{\sss{1}}$)}~
& \left\{~
\begin{aligned}
\Delta u^{\sss{1}} + k^{2} u^{\sss{1}} &= P_{m}(x,u^{\sss{0}})
& & \text{in\ } \Omega, \\
u^{\sss{1}} &= 0
& & \text{on\ } \partial\Omega. \\
\end{aligned}
\right. \label{eqn:probl_I1}
\end{align}
Then, given any Dirichlet boundary data $u^{\sss{0}} |_{\partial\Omega} = f$ for the unperturbed problem $(I^{\sss{0}})$, the solution $u^{\sss{0}}$ in $\Omega$ is obtained by solving \eqref{eqn:probl_I0}.
Similarly, the corresponding Neumann boundary data $\partial_{\nu} u^{\sss{1}} |_{\partial\Omega}$ for the linearized problem $(I^{\sss{1}})$ can be computed by solving \eqref{eqn:probl_I1} with the polynomial $P_{m}(x,u^{\sss{0}})$ defined in $\Omega$.

Thus, the linearized DtN map $\Lambda^{\prime}_{\bm{c}}$ of the original DtN map $\Lambda_{\bm{c}}$ in \eqref{eqn:DtN_I} is given by
\begin{align}\label{eqn:linearDtN_I1}
\Lambda^{\prime}_{\bm{c}} : u^{\sss{0}} |_{\partial\Omega} \mapsto \partial_{\nu} u^{\sss{1}} |_{\partial\Omega}.
\end{align}
Here, $u^{\sss{0}} |_{\partial\Omega} \in \mathcal{D}$ is the Dirichlet boundary data of the unperturbed problem ($I^{\sss{0}}$) in \eqref{eqn:probl_I0}, and $\partial_{\nu} u^{\sss{1}} |_{\partial\Omega} \in \mathcal{N}$ is the Neumann boundary data of the linearized problem ($I^{\sss{1}}$) in \eqref{eqn:probl_I1}.
The linearized DtN map $\Lambda^{\prime}_{\bm{c}}$ is well-defined for small $f$ and is linear with respect to $\bm{c} \in (L^{\infty}(\Omega))^{m}$.
Indeed, $\Lambda^{\prime}_{\bm{c}}$ is the first-order linearization of DtN map $\Lambda_{\bm{c}}$ for small $\bm{c}$.
A similar result has been shown in \cite[Proposition 2.3]{ZLX2024} for a single nonlinearity term $c_{m} u^{m}$ and we skip these detailed discussion.

Multiplying the linearized problem ($I^{\sss{1}}$) in \eqref{eqn:probl_I1} from both sides with a test solution $\varphi$, which satifies $\Delta \varphi + k^{2} \varphi = 0$ in $\Omega$, we obtain the Alessandrini-type identity for the linearized DtN map $\Lambda^{\prime}_{\bm{c}}$,
\begin{align}\label{eqn:identity_I1}
\int_{\Omega} \varphi \, P_{m}(x,u^{\sss{0}}) \,\rmd x
&= \int_{\partial\Omega} \varphi \, \partial_{\nu} u^{\sss{1}} \,\rmd s_{x}
= \int_{\partial\Omega} \varphi \, \Lambda^{\prime}_{\bm{c}} (u^{\sss{0}}|_{\partial\Omega}) \,\rmd s_{x}
= \int_{\partial\Omega} \varphi \, \Lambda^{\prime}_{\bm{c}} (f) \,\rmd s_{x}.
\end{align}
Notice that in the above equality, the Neumann boundary data $\partial_{\nu} u^{\sss{1}} |_{\partial\Omega}$ cannot be directly measured.
Nevertheless, as mentioned earlier, $\Lambda^{\prime}_{\bm{c}}$ represents the first-order linearization of DtN map $\Lambda_{\bm{c}}$ for sufficient small $\bm{c}$.
Consequently, $\partial_{\nu} u^{\sss{1}} |_{\partial\Omega}$ can be approximated as follows
\begin{align}\label{eqn:Neumann_approximate}
\partial_{\nu} u^{\sss{1}} |_{\partial\Omega}
~\approx~ \partial_{\nu} u |_{\partial\Omega} - \partial_{\nu} u^{\sss{0}} |_{\partial\Omega},
\end{align}
where a similar approach has proven effective in single-coefficient identification problems; see \cite{DS1994, ILX2020, ILX2022, ZLX2024}.

We emphasize that, to numerically recover the entire coefficient vector $\bm{c}$, the aforementioned equality \eqref{eqn:identity_I1} serves as the foundational starting point for our recursive algorithm.
The primary challenge then lies in recursively reconstructing each coefficient $c_{\ell}$ ($\ell = 1,2,\dots,m$).
In the following section, we will systematically present the core methodology by applying the principle of inclusion-exclusion.

\section{The recursive reconstruction formulas}\label{sec:scheme}

In this section, we outline and formulate our recursive reconstruction scheme for recovering all coefficients $\{ c_{\ell} \}_{\ell=1}^{m}$ (or $\bm{c}$) from the linearized DtN map $\Lambda^{\prime}_{\bm{c}}$ given in \eqref{eqn:linearDtN_I1}, along with the corresponding identity \eqref{eqn:identity_I1}.

To successfully recover the entire coefficient vector $\bm{c} = (c_{1},c_{2},\dots,c_{m})$, we must recursively identify each component $c_{\ell}$ for $\ell = 1,2,\dots,m$.
For this purpose, the principle of inclusion-exclusion plays a crucial role and we now detail its application to our multi-coefficient identification problem below.

\subsection{The principle of inclusion-exclusion}

Let $\ell \in \mathbb{N}_{+}$ be a postive integer, we first introduce some notations.
The $\ell$-dimensional multi-index is an $\ell$-tuple of non-negative integers, denoted by
\begin{align*}
\bm{\alpha} := (\alpha_{1},\alpha_{2},\dots,\alpha_{\ell}) \in \mathbb{N}_{0}^{\ell}.
\end{align*}
In particular, $\bm{1} := (1,1,\dots,1) \in \mathbb{N}_{0}^{\ell}$.
We further denote the absolute value of $\bm{\alpha}$ as
\begin{align*}
|\bm{\alpha}| := \alpha_{1} + \alpha_{2} + \cdots + \alpha_{\ell},
\end{align*}
and denote the multinomial coefficient as
\begin{align*}
{\binom{\ell}{\bm{\alpha}}} := \frac{\ell\,!}{\,\alpha_{1}!\,\alpha_{2}!\,\cdots\,\alpha_{\ell}!\,}.
\end{align*}

Then, given a non-negative integer $a \in \mathbb{N}_{0}$, and a finite set of multi-index
\begin{align*}
\mathcal{A}^{\ell}_{a} := \big\{ \bm{\alpha} \in \mathbb{N}_{0}^{\ell} : |\bm{\alpha}| = a \big\},
\end{align*}
the following identities can be easily derived from the principle of inclusion-exclusion in combinatorics.
These identities will later play a key role in the stable reconstruction method.

\begin{lemma}\label{lmm:PIE}
Let $\ell \in \mathbb{N}_{+}$, we denote the index set
\begin{align*}
\mathcal{U}_{\ell} := \{ 1, 2, \dots, \ell \},
\end{align*}
and denote $\mathcal{S}$ to be any non-empty subset of $\mathcal{U}_{\ell}$.
Then, given a set of functions $\big\{ w_{j} : j \in \mathcal{U}_{\ell} \big\}$, for any $\ell' \in \mathbb{N}_{+}$, there holds that
\begin{align}\label{eqn:identity_PIE}
\frac{1}{\,\ell\,!\,} \, {\textstyle\sum\limits_{\emptyset \subsetneqq \mathcal{S} \subseteq \mathcal{U}_{\ell}}} (-1)^{|\mathcal{U}_{\ell} \setminus \mathcal{S}|} \Big( {\textstyle\sum\limits_{j \in \mathcal{S}}} w_{j} \Big)^{\ell'}
= \left\{~
\begin{aligned}
\\[-1.5ex]
& ~ 0, & 0 < \ell' < \ell, \\[2.0ex]
& {\textstyle\prod\limits_{j \in \mathcal{U}_{\ell}}} w_{j}, & \ell' = \ell, \\
& {\textstyle\prod\limits_{j \in \mathcal{U}_{\ell}}} w_{j} \, Q_{\ell,\ell'\!-\ell}(w_{1},w_{2},\dots,w_{\ell}), & \ell' > \ell, \\
\end{aligned}
\right.
\end{align}
where $\mathcal{U}_{\ell} \setminus \mathcal{S}$ is the relative complement of $\mathcal{S}$ in $\mathcal{U}_{\ell}$, $|\cdot|$ is the cardinality of a set, and $Q_{\ell,\ell'\!-\ell}(w_{1},w_{2},\dots,w_{\ell})$ is a multivariate polynomial in $\big\{ w_{j} : j \in \mathcal{U}_{\ell} \big\}$ of degree $\ell'-\ell$.
Here, the multivariate polynomial $Q_{\ell,a}(w_{1},w_{2},\dots,w_{\ell})$ of degree $a$ is defined by
\begin{align}\label{eqn:poly_Qa}
Q_{\ell,a}(w_{1},w_{2},\dots,w_{\ell})
:= \frac{1}{\,\ell\,!\,} {\textstyle\sum\limits_{\bm{\alpha} \in \mathcal{A}^{\ell}_{a}}} {\textstyle\binom{\ell+a}{\bm{1}+\bm{\alpha}}} \, w_{1}^{\alpha_{1}} w_{2}^{\alpha_{2}} \cdots w_{\ell}^{\alpha_{\ell}},
\quad a \in \mathbb{N}_{0}.
\end{align}
Notice that, $Q_{\ell,a} \equiv 1$ when $a = 0$.
\end{lemma}

The proof of this lemma employs a standard combinatorial technique for handling subtraction operations; we omit the details and refer the reader to \cite{GGL1995} for a complete treatment.
A related result addressing single-coefficient inversion in the nonlinear Helmholtz equation can be found in our joint work \cite[Lemma 2.4]{ZLX2024}.
To provide some explicit representations of the polynomial $Q_{\ell,a}$, we state the following remark.

\begin{remark}
we exhibit three specific cases of the polynomial $Q_{\ell,a}$.
\begin{align*}
Q_{\ell,1}(w_{1},w_{2},\dots,w_{\ell})
&= \frac{(\ell+1)!}{\,\ell\,!\,2!\,} {\textstyle\sum\limits_{1 \leqslant j \leqslant \ell}} w_{j}, \\
Q_{\ell,2}(w_{1},w_{2},\dots,w_{\ell})
&= \frac{(\ell+2)!}{\,\ell\,!\,3!\,} {\textstyle\sum\limits_{1 \leqslant j \leqslant \ell}} w_{j}^{2}
+ \frac{(\ell+2)!}{\,\ell\,!\,2!\,2!\,} {\textstyle\sum\limits_{1 \leqslant j < j' \leqslant \ell}} w_{j} w_{j'}, \\
Q_{\ell,3}(w_{1},w_{2},\dots,w_{\ell})
&= \frac{(\ell+3)!}{\,\ell\,!\,4!\,} {\textstyle\sum\limits_{1 \leqslant j \leqslant \ell}} w_{j}^{3}
+ \frac{(\ell+3)!}{\,\ell\,!\,3!\,2!\,} {\textstyle\sum\limits_{1 \leqslant j < j' \leqslant \ell}} ( w_{j}^{2} w_{j'} + w_{j} w_{j'}^{2} ) \\
&\qquad + \frac{(\ell+3)!}{\,\ell\,!\,2!\,2!\,2!\,} {\textstyle\sum\limits_{1 \leqslant j < j' < j'' \leqslant \ell}} w_{j} w_{j'} w_{j''}.
\end{align*}
\end{remark}

To incorporate \textbf{Lemma \ref{lmm:PIE}} into our multi-coefficient inversion framework for the nonlinear Helmholtz equation \eqref{eqn:probl_I}, we begin by defining $u_{\ell,j}^{\sss{0}}$ as the solution to the linear Helmholtz equation \eqref{eqn:probl_I0} with Dirichlet boundary condition $f_{\ell,j}$, for each index $j \in \mathcal{U}_{\ell} = \{ 1, 2, \dots, \ell \}$, such that
\begin{align*}
\text{($I_{\ell,j}^{\sss{0}}$)}~
& \left\{~
\begin{aligned}
\Delta u_{\ell,j}^{\sss{0}} + k^{2} u_{\ell,j}^{\sss{0}} &= 0
& & \text{in\ } \Omega, \\
u_{\ell,j}^{\sss{0}} &= f_{\ell,j}
& & \text{on\ } \partial\Omega. \\
\end{aligned}
\right.
\end{align*}
Then, let $\mathcal{S}$ be any non-empty subset of $\mathcal{U}_{\ell}$ (i.e., $\emptyset \subsetneqq \mathcal{S} \subseteq \mathcal{U}_{\ell}$), since the problems ($I^{\sss{0}}$) and ($I_{\ell,j}^{\sss{0}}$) are linear, we also have that
\begin{align*}
u_{\ell,\mathcal{S}}^{\sss{0}}
:= {\textstyle\sum\limits_{j \in \mathcal{S}}} u_{\ell,j}^{\sss{0}}
\quad \text{in\ } \Omega
\qquad \text{and} \qquad
f_{\ell,\mathcal{S}}
:= {\textstyle\sum\limits_{j \in \mathcal{S}}} f_{\ell,j}
\quad \text{on\ } \partial\Omega
\end{align*}
are the combined solution and the combined Dirichlet boundary date of the problem
\begin{align*}
\text{($I_{\ell,\mathcal{S}}^{\sss{0}}$)}~
& \left\{~
\begin{aligned}
\Delta u_{\ell,\mathcal{S}}^{\sss{0}} + k^{2} u_{\ell,\mathcal{S}}^{\sss{0}} &= 0
& & \text{in\ } \Omega, \\
u_{\ell,\mathcal{S}}^{\sss{0}} &= f_{\ell,\mathcal{S}}
& & \text{on\ } \partial\Omega, \\
\end{aligned}
\right.
\end{align*}
respectively.
Therefore, we can denote $u_{\ell,\mathcal{S}}^{\sss{1}}$ as the corresponding solution of the linearized problem \eqref{eqn:probl_I1} with $P_{m}(x,u_{\ell,\mathcal{S}}^{\sss{0}})$ in $\Omega$,
\begin{align*}
\text{($I_{\ell,\mathcal{S}}^{\sss{1}}$)}~
& \left\{~
\begin{aligned}
\Delta u_{\ell,\mathcal{S}}^{\sss{1}} + k^{2} u_{\ell,\mathcal{S}}^{\sss{1}} &= P_{m}(x,u_{\ell,\mathcal{S}}^{\sss{0}})
& & \text{in\ } \Omega, \\
u_{\ell,\mathcal{S}}^{\sss{1}} &= 0
& & \text{on\ } \partial\Omega. \\
\end{aligned}
\right.
\end{align*}
Thus, for each non-empty subset $\mathcal{S}$ of $\mathcal{U}_{\ell}$, we also derive a similar identity
\begin{align}\label{eqn:identity_I1S}
\int_{\Omega} \varphi_{\ell} \, P_{m}(x,u_{\ell,\mathcal{S}}^{\sss{0}}) \,\rmd x
&= \int_{\partial\Omega} \varphi_{\ell} \, \partial_{\nu} u_{\ell,\mathcal{S}}^{\sss{1}} \,\rmd s_{x}
= \int_{\partial\Omega} \varphi_{\ell} \, \Lambda^{\prime}_{\bm{c}} (f_{\ell,\mathcal{S}}) \,\rmd s_{x}
\end{align}
by following the same method as that used for the Alessandrini-type identity \eqref{eqn:identity_I1}.
Here, the test solution $\varphi_{\ell}$ satisfies Helmholtz equation $\Delta \varphi_{\ell} + k^{2} \varphi_{\ell} = 0$ in $\Omega$.

Combining the polynomial $P_{m}(x,u)$ from \eqref{eqn:poly_Pm}, the combinatorial identities \eqref{eqn:identity_PIE} in \textbf{Lemma \ref{lmm:PIE}}, and the relation \eqref{eqn:identity_I1S}, we derive the following Alessandrini-type identities based on the principle of inclusion-exclusion.
These identities constitute the fundamental framework for our recursive inversion formulas.

\begin{lemma}\label{lmm:APIE}
Let the sets $\mathcal{U}_{\ell}$ and $\mathcal{S}$, the solutions $\varphi_{\ell}$ and $u_{\ell,j}^{\sss{0}}$, and the linearized DtN map $\Lambda^{\prime}_{\bm{c}}$ be defined as above.
Given $f_{\ell,0} := \varphi_{\ell} |_{\partial\Omega}$, $\big\{ f_{\ell,j} = u_{\ell,j}^{\sss{0}} |_{\partial\Omega} : j \in \mathcal{U}_{\ell} \big\}$ and $f_{\ell,\mathcal{S}} = {\textstyle\sum\limits_{j \in \mathcal{S}}} f_{\ell,j}$, we obtain that
\begin{enumerate}[label=\textbf{(\arabic*)}]

\item \textbf{Case $\ell = m$:}
\begin{align}\label{eqn:APIE_m}
\int_{\Omega} c_{m}(x) \, \varphi_{m} {\textstyle\prod\limits_{j \in \mathcal{U}_{m}}} u_{m,j}^{\sss{0}} \,\rmd x
&= \frac{1}{\,m!\,} \, {\textstyle\sum\limits_{\emptyset \subsetneqq \mathcal{S} \subseteq \mathcal{U}_{m}}} (-1)^{|\mathcal{U}_{m} \setminus \mathcal{S}|} \int_{\partial\Omega} f_{m,0} \, \Lambda^{\prime}_{\bm{c}} (f_{m,\mathcal{S}}) \,\rmd s_{x},
\end{align}

\item \textbf{Case $0 < \ell < m$:}
\begin{align}\label{eqn:APIE_l}
\begin{aligned}
\int_{\Omega} c_{\ell}(x) \, \varphi_{\ell} {\textstyle\prod\limits_{j \in \mathcal{U}_{\ell}}} u_{\ell,j}^{\sss{0}} \,\rmd x
&+ {\textstyle\sum\limits_{a = 1}^{m-\ell}} \int_{\Omega} c_{\ell+a}(x) \, Q_{\ell,a}( u_{\ell,1}^{\sss{0}}, u_{\ell,2}^{\sss{0}}, \dots, u_{\ell,\ell}^{\sss{0}} ) \, \varphi_{\ell} {\textstyle\prod\limits_{j \in \mathcal{U}_{\ell}}} u_{\ell,j}^{\sss{0}} \,\rmd x \\
&\qquad = \frac{1}{\,\ell\,!\,} \, {\textstyle\sum\limits_{\emptyset \subsetneqq \mathcal{S} \subseteq \mathcal{U}_{\ell}}} (-1)^{|\mathcal{U}_{\ell} \setminus \mathcal{S}|} \int_{\partial\Omega} f_{\ell,0} \, \Lambda^{\prime}_{\bm{c}} (f_{\ell,\mathcal{S}}) \,\rmd s_{x}.
\end{aligned}
\end{align}

\end{enumerate}
Here, $Q_{\ell,a}( u_{\ell,1}^{\sss{0}}, u_{\ell,2}^{\sss{0}}, \dots, u_{\ell,\ell}^{\sss{0}} )$ is the multivariate polynomial in $\big\{ u_{\ell,j}^{\sss{0}} : j \in \mathcal{U}_{\ell} \big\}$ of degree $a$ defined by \eqref{eqn:poly_Qa}.
\end{lemma}

\begin{proof}
The proof follows by applying the combinatorial identities \eqref{eqn:identity_PIE} from \textbf{Lemma \ref{lmm:PIE}} and systematically verifying the cancellation of each term in the polynomial expansion of $P_{m}(x,u)$ given in \eqref{eqn:poly_Pm}.
The key steps are summarized as follows:
\begin{enumerate}[label=\textbf{(\arabic*)}]

\item For \textbf{Case $\ell = m$}, we have
\begin{align*}
&\frac{1}{\,m!\,} \, {\textstyle\sum\limits_{\emptyset \subsetneqq \mathcal{S} \subseteq \mathcal{U}_{m}}} (-1)^{|\mathcal{U}_{m} \setminus \mathcal{S}|} P_{m}(x,u_{m,\mathcal{S}}^{\sss{0}}) \\
&\qquad = {\textstyle\sum\limits_{\ell' = 1}^{m}} c_{\ell'}(x) \, \frac{1}{\,m!\,} \, {\textstyle\sum\limits_{\emptyset \subsetneqq \mathcal{S} \subseteq \mathcal{U}_{m}}} (-1)^{|\mathcal{U}_{m} \setminus \mathcal{S}|} (u_{m,\mathcal{S}}^{\sss{0}})^{\ell'} \\
&\qquad = c_{m}(x) {\textstyle\prod\limits_{j \in \mathcal{U}_{m}}} u_{m,j}^{\sss{0}}.
\end{align*}

\item  For \textbf{Case $0 < \ell < m$}, we have
\begin{align*}
&\frac{1}{\,\ell\,!\,} \, {\textstyle\sum\limits_{\emptyset \subsetneqq \mathcal{S} \subseteq \mathcal{U}_{\ell}}} (-1)^{|\mathcal{U}_{\ell} \setminus \mathcal{S}|} P_{m}(x,u_{\ell,\mathcal{S}}^{\sss{0}}) \\
&\qquad = {\textstyle\sum\limits_{\ell' = 1}^{m}} c_{\ell'}(x) \, \frac{1}{\,\ell\,!\,} \, {\textstyle\sum\limits_{\emptyset \subsetneqq \mathcal{S} \subseteq \mathcal{U}_{\ell}}} (-1)^{|\mathcal{U}_{\ell} \setminus \mathcal{S}|} (u_{\ell,\mathcal{S}}^{\sss{0}})^{\ell'} \\
&\qquad = \Big( c_{\ell}(x) + {\textstyle\sum\limits_{a = 1}^{m-\ell}} c_{\ell+a}(x) \, Q_{\ell,a}( u_{\ell,1}^{\sss{0}}, u_{\ell,2}^{\sss{0}}, \dots, u_{\ell,\ell}^{\sss{0}} ) \Big) {\textstyle\prod\limits_{j \in \mathcal{U}_{\ell}}} u_{\ell,j}^{\sss{0}}.
\end{align*}

\end{enumerate}
The proof is completed by substituting the above equations into the identity \eqref{eqn:identity_I1S}.
\end{proof}

We shall show that the first terms on the left hand side of identities \eqref{eqn:APIE_m} and \eqref{eqn:APIE_l} will prove to be useful.
In particular, by carefully choosing the test functions and solutions for the Helmholtz equation \eqref{eqn:probl_I0}, a strategy that will be discussed in the next two subsections, we can transform the integral
\begin{align*}
\int_{\Omega} c_{m}(x) \, \varphi_{m} {\textstyle\prod\limits_{j \in \mathcal{U}_{m}}} u_{m,j}^{\sss{0}} \,\rmd x,
\end{align*}
for instance, into a Fourier mode of the unknown coefficient $c_{m}$.
Then, as \textbf{Lemma \ref{lmm:APIE}} demonstrates, recovering the lower-order coefficient $c_{\ell}$ of polynomial $P_{m}(x,u)$ requires accounting for the coupling with higher-order coefficient $c_{\ell'}$, where $\ell < \ell' \leqslant m$.
The decoupling of these coefficients, along with the construction of complex exponential solutions for the Helmholtz equations \eqref{eqn:probl_I0}, will be addressed in the following subsections.

\subsection{The reconstruction formulas via Fourier transform}

Let $\xi \in \mathbb{R}^{n} \setminus \{0\}$ be a vector in the Fourier frequency space.
We set the solutions $\varphi_{\ell}$ and $\big\{ u_{\ell,j}^{\sss{0}} : j \in \mathcal{U}_{\ell} \big\}$ of Helmholtz equation in $\Omega$ to be the complex exponential solutions as below,
\begin{align}\label{eqn:CE_solution}
\varphi_{\ell}(x) := \rme^{\img x \cdot \zeta_{\ell,0}},
\qquad
u_{\ell,j}^{\sss{0}}(x) := \rme^{\img x \cdot \zeta_{\ell,j}},
\quad j \in \mathcal{U}_{\ell},
\end{align}
where $\img = \sqrt{-1}$ and the auxiliary complex vectors $\zeta_{\ell,0} \in \mathbb{C}^{n}$ and $\big\{ \zeta_{\ell,j} \in \mathbb{C}^{n} : j \in \mathcal{U}_{\ell} \big\}$ satisfy
\begin{align}\label{eqn:zeta_condition}
\left\{
\begin{aligned}
\zeta_{\ell,0} \cdot \zeta_{\ell,0} &= k^{2}, \\
\zeta_{\ell,j} \cdot \zeta_{\ell,j} &= k^{2}, \quad j \in \mathcal{U}_{\ell}, \\
\zeta_{\ell,0} + {\textstyle\sum\limits_{j \in \mathcal{U}_{\ell}}} \zeta_{\ell,j} &= \xi.
\end{aligned}
\right.
\end{align}

Noticing that $\mathop{\mathrm{supp}}(c_{\ell}) \subset \Omega$ ($\ell = 1,2,\dots,m$), we can denote
\begin{align*}
\int_{\Omega} c_{\ell}(x) \, \varphi_{\ell} {\textstyle\prod\limits_{j \in \mathcal{U}_{\ell}}} u_{\ell,j}^{\sss{0}} \,\rmd x
= \int_{\Omega} c_{\ell}(x) \, \rme^{\img x \cdot \xi} \,\rmd x
= \int_{\mathbb{R}^{n}} c_{\ell}(x) \, \rme^{\img x \cdot \xi} \,\rmd x
=: \mathcal{F}[c_{\ell}](\xi)
\end{align*}
as the Fourier transform of the coefficient $c_{\ell}$ at frequency $\xi$.
Therefore, according to \textbf{Lemma \ref{lmm:APIE}}, the reconstruction formulas via the Fourier transform $\mathcal{F}$ towards the unknown coefficients $\{ c_{\ell} \}_{\ell=1}^{m}$ are now given below.

\begin{lemma}\label{lmm:formulas}
The reconstruction formulas for the coefficients $\{ c_{\ell} \}_{\ell=1}^{m}$ of the polynomial $P_{m}(x,u)$ using the Fourier transform $\mathcal{F}$ and the linearized DtN map $\Lambda^{\prime}_{\bm{c}}$ are given by
\begin{enumerate}[label=\textbf{(\arabic*)}]

\item \textbf{Case $c_{m}$:}
\begin{align}\label{eqn:formula_m}
\mathcal{F}[c_{m}]
&= d_{m},
\end{align}

\item \textbf{Case $c_{\ell}$ with $0 < \ell < m$:}
\begin{align}\label{eqn:formula_l}
\begin{aligned}
\mathcal{F}[c_{\ell}]
+ {\textstyle\sum\limits_{a = 1}^{m-\ell}} \mathcal{F}[c_{\ell+a}Q_{\ell,a}]
&= d_{\ell}.
\end{aligned}
\end{align}

\end{enumerate}
Here, for convenience, we denote ``the data'' $d_{\ell}$ by
\begin{align}\label{eqn:data_d_l}
d_{\ell}
:= \frac{1}{\,\ell\,!\,} \, {\textstyle\sum\limits_{\emptyset \subsetneqq \mathcal{S} \subseteq \mathcal{U}_{\ell}}} (-1)^{|\mathcal{U}_{\ell} \setminus \mathcal{S}|} \int_{\partial\Omega} f_{\ell,0} \, \Lambda^{\prime}_{\bm{c}} (f_{\ell,\mathcal{S}}) \,\rmd s_{x},
\end{align}
where $f_{\ell,0} = \rme^{\img x \cdot \zeta_{\ell,0}} |_{\partial\Omega}$, $\big\{ f_{\ell,j} = \rme^{\img x \cdot \zeta_{\ell,j}} |_{\partial\Omega} : j \in \mathcal{U}_{\ell} \big\}$, and $f_{\ell,\mathcal{S}} = {\textstyle\sum\limits_{j \in \mathcal{S}}} f_{\ell,j}$ are the Dirichlet boundary data.
\end{lemma}

\begin{proof}
The proof follows from the identities \eqref{eqn:APIE_m} and \eqref{eqn:APIE_l} in \textbf{Lemma \ref{lmm:APIE}}, by observing that
\begin{align*}
\varphi_{\ell} {\textstyle\prod\limits_{j \in \mathcal{U}_{\ell}}} u_{\ell,j}^{\sss{0}}
= \rme^{\img x \cdot \xi},
\end{align*}
where the solutions are given in \eqref{eqn:CE_solution} and the auxiliary complex vectors satisfy the conditions \eqref{eqn:zeta_condition}.
\end{proof}



It is emphasized that the lower-order coefficient $c_{\ell}$ of polynomial $P_{m}(x,u)$ can be recovered recursively using the identity \eqref{eqn:APIE_l} or the formula \eqref{eqn:formula_l} for \textbf{Case $0 < \ell < m$}.
For specific choices of solutions \eqref{eqn:CE_solution} and complex vectors \eqref{eqn:zeta_condition}, the polynomial $Q_{\ell,a}$ \eqref{eqn:poly_Qa} simplifies to
\begin{align*}
Q_{\ell,a}( u_{\ell,1}^{\sss{0}}, u_{\ell,2}^{\sss{0}}, \dots, u_{\ell,\ell}^{\sss{0}} )
&= \frac{1}{\,\ell\,!\,} {\textstyle\sum\limits_{\bm{\alpha} \in \mathcal{A}^{\ell}_{a}}} {\textstyle\binom{\ell+a}{\bm{1}+\bm{\alpha}}} \, (u_{\ell,1}^{\sss{0}})^{\alpha_{1}} (u_{\ell,2}^{\sss{0}})^{\alpha_{2}} \cdots (u_{\ell,\ell}^{\sss{0}})^{\alpha_{\ell}} \\
&= \frac{1}{\,\ell\,!\,} {\textstyle\sum\limits_{\bm{\alpha} \in \mathcal{A}^{\ell}_{a}}} {\textstyle\binom{\ell+a}{\bm{1}+\bm{\alpha}}} \, \exp\big\{\img x \cdot ( {\textstyle\sum\limits_{j \in \mathcal{U}_{\ell}}} \alpha_{j} \zeta_{\ell,j} )\big\}.
\end{align*}
Then, we have that
\begin{align*}
\mathcal{F}[c_{\ell+a}Q_{\ell,a}](\xi)
= \frac{1}{\,\ell\,!\,} {\textstyle\sum\limits_{\bm{\alpha} \in \mathcal{A}^{\ell}_{a}}} {\textstyle\binom{\ell+a}{\bm{1}+\bm{\alpha}}} \, \mathcal{F}[c_{\ell+a}]( \xi + {\textstyle\sum\limits_{j \in \mathcal{U}_{\ell}}} \alpha_{j} \zeta_{\ell,j} ).
\end{align*}
Defining the linear operator $\mathcal{Q}_{\ell,a}$ by
\begin{align}
\mathcal{Q}_{\ell,a}(\mathcal{F}[c_{\ell+a}]) := \mathcal{F}[c_{\ell+a}Q_{\ell,a}],
\end{align}
we conclude that the operator $\mathcal{Q}_{\ell,a}$ is a linear combination of several translation operators acting on $\mathcal{F}[c_{\ell+a}]$. 
Thus, denote the identity operator by $\mathcal{I}$, the formulas \eqref{eqn:formula_m} and \eqref{eqn:formula_l} admit the matrix formulation
\begin{align}\label{eqn:formula_matrix}
\begin{pmatrix}
~ \mathcal{I} & \mathcal{Q}_{1,1} & \mathcal{Q}_{1,2} & \cdots & \mathcal{Q}_{1,m-1} \\
& \mathcal{I} & \mathcal{Q}_{2,1} & \cdots & \mathcal{Q}_{2,m-2} \\
& & \ddots & \ddots & \vdots \\
& & & \mathcal{I} & \mathcal{Q}_{m-1,1} \\
& & & & \mathcal{I} \\ 
\end{pmatrix}
\begin{pmatrix}
\mathcal{F}[c_{1}] \\
\mathcal{F}[c_{2}] \\
\vdots \\
\mathcal{F}[c_{m-1}] \\
\mathcal{F}[c_{m}] \\
\end{pmatrix}
&= \begin{pmatrix}
d_{1} \\
d_{2} \\
\vdots \\
d_{m-1} \\
d_{m} \\
\end{pmatrix},
\end{align}
which can be solved by back substitution.

\begin{remark}
Take $m = 4$, for example, we have
\begin{align*}
\left\{
\begin{aligned}
\mathcal{F}[c_{1}] + \mathcal{F}[c_{2}Q_{1,1}] + \mathcal{F}[c_{3}Q_{1,2}] + \mathcal{F}[c_{4}Q_{1,3}] &= d_{1}, \\
\mathcal{F}[c_{2}] + \mathcal{F}[c_{3}Q_{2,1}] + \mathcal{F}[c_{4}Q_{2,2}] &= d_{2}, \\
\mathcal{F}[c_{3}] + \mathcal{F}[c_{4}Q_{3,1}] &= d_{3}, \\
\mathcal{F}[c_{4}] &= d_{4}, \\
\end{aligned}
\right.
\end{align*}
by the formulas \eqref{eqn:formula_m} and \eqref{eqn:formula_l}.
Then, $\mathcal{F}[c_{1}]$, $\mathcal{F}[c_{2}]$, $\mathcal{F}[c_{3}]$ and $\mathcal{F}[c_{4}]$ can be recovered recursively by
\begin{align*}
\left\{
\begin{aligned}
\mathcal{F}[c_{4}] &= d_{4}, \\
\mathcal{F}[c_{3}] &= d_{3}
- \mathcal{Q}_{3,1}(\mathcal{F}[c_{4}]), \\
\mathcal{F}[c_{2}] &= d_{2}
- \mathcal{Q}_{2,2}(\mathcal{F}[c_{4}])
- \mathcal{Q}_{2,1}(\mathcal{F}[c_{3}]), \\
\mathcal{F}[c_{1}] &= d_{1}
- \mathcal{Q}_{1,3}(\mathcal{F}[c_{4}])
- \mathcal{Q}_{1,2}(\mathcal{F}[c_{3}])
- \mathcal{Q}_{1,1}(\mathcal{F}[c_{2}]). \\
\end{aligned}
\right.
\end{align*}
Moreover, we also present several examples of the polynomial $Q_{\ell,a}$,
\begin{align*}
Q_{1,1}(u_{1,1}^{\sss{0}})
&= u_{1,1}^{\sss{0}}
= \rme^{\img x \cdot \zeta_{1,1}}, \\
Q_{1,2}(u_{1,1}^{\sss{0}})
&= (u_{1,1}^{\sss{0}})^{2}
= \rme^{\img x \cdot (2\zeta_{1,1})}, \\
Q_{1,3}(u_{1,1}^{\sss{0}})
&= (u_{1,1}^{\sss{0}})^{3}
= \rme^{\img x \cdot (3\zeta_{1,1})}, \\
Q_{2,1}(u_{2,1}^{\sss{0}},u_{2,2}^{\sss{0}})
&= \tfrac{3}{2} (u_{2,1}^{\sss{0}}+u_{2,2}^{\sss{0}})
= \tfrac{3}{2} ( \rme^{\img x \cdot \zeta_{2,1}} + \rme^{\img x \cdot \zeta_{2,2}} ), \\[1ex]
Q_{2,2}(u_{2,1}^{\sss{0}},u_{2,2}^{\sss{0}})
&= 2 \big( (u_{2,1}^{\sss{0}})^{2} + (u_{2,2}^{\sss{0}})^{2} \big) + 3 u_{2,1}^{\sss{0}} u_{2,2}^{\sss{0}} \\
&= 2 ( \rme^{\img x \cdot (2\zeta_{2,1})} + \rme^{\img x \cdot (2\zeta_{2,2})} ) + 3 \rme^{\img x \cdot (\zeta_{2,1}+\zeta_{2,2})}, \\[1ex]
Q_{3,1}(u_{3,1}^{\sss{0}},u_{3,2}^{\sss{0}},u_{3,3}^{\sss{0}})
&= 2 ( u_{3,1}^{\sss{0}} + u_{3,2}^{\sss{0}} + u_{3,3}^{\sss{0}} ) \\
&= 2 ( \rme^{\img x \cdot \zeta_{3,1}} + \rme^{\img x \cdot \zeta_{3,2}} + \rme^{\img x \cdot \zeta_{3,3}} ).
\end{align*}
Then, we have that
\begin{align*}
\mathcal{Q}_{1,1}(\mathcal{F}[c_{2}])(\xi)
&= \mathcal{F}[c_{2}Q_{1,1}](\xi)
= \mathcal{F}[c_{2}](\xi+\zeta_{1,1}), \\
\mathcal{Q}_{1,2}(\mathcal{F}[c_{3}])(\xi)
&= \mathcal{F}[c_{3}Q_{1,2}](\xi)
= \mathcal{F}[c_{3}](\xi+2\zeta_{1,1}), \\
\mathcal{Q}_{1,3}(\mathcal{F}[c_{4}])(\xi)
&= \mathcal{F}[c_{4}Q_{1,3}](\xi)
= \mathcal{F}[c_{4}](\xi+3\zeta_{1,1}), \\
\mathcal{Q}_{2,1}(\mathcal{F}[c_{3}])(\xi)
&= \mathcal{F}[c_{3}Q_{2,1}](\xi)
= \tfrac{3}{2} \big\{ \mathcal{F}[c_{3}](\xi+\zeta_{2,1}) + \mathcal{F}[c_{3}](\xi+\zeta_{2,2}) \big\}, \\[1ex]
\mathcal{Q}_{2,2}(\mathcal{F}[c_{4}])(\xi)
&= \mathcal{F}[c_{4}Q_{2,2}](\xi) \\
&= 2 \big\{ \mathcal{F}[c_{4}](\xi+2\zeta_{2,1}) + \mathcal{F}[c_{4}](\xi+2\zeta_{2,2}) \big\} + 3 \mathcal{F}[c_{4}](\xi+\zeta_{2,1}+\zeta_{2,2}), \\[1ex]
\mathcal{Q}_{3,1}(\mathcal{F}[c_{4}])(\xi)
&= \mathcal{F}[c_{4}Q_{3,1}](\xi) \\
&= 2 \big\{ \mathcal{F}[c_{4}](\xi+\zeta_{3,1}) + \mathcal{F}[c_{4}](\xi+\zeta_{3,2}) + \mathcal{F}[c_{4}](\xi+\zeta_{3,3}) \big\}.
\end{align*}
\end{remark}

\subsection{The construction of auxiliary complex vectors}

We now proceed to properly construct the auxiliary complex vectors $\zeta_{\ell,0}$ and $\big\{ \zeta_{\ell,j} : j \in \mathcal{U}_{\ell} \big\}$ mentioned in the previous subsection.
Since $\xi \neq 0$, we first define an orthonormal basis of $\mathbb{R}^{n}$ for $n \geqslant 2$,
\begin{align*}
\big\{ e_{1} := \tfrac{\xi}{|\xi|}, e_{2}, \dots, e_{n} \big\},
\end{align*}
and we also define two complex vectors $\zeta_{\ell}^{+} \in \mathbb{C}^{n}$ and $\zeta_{\ell}^{-} \in \mathbb{C}^{n}$ for different $\ell \in \mathbb{N}_{+}$,
\begin{align}\label{eqn:zeta_pm}
\zeta_{\ell}^{\pm} := \left\{~
\begin{aligned}
\tfrac{|\xi|}{\ell+1} \, e_{1} &\pm \sqrt{ k^{2} - ( \tfrac{|\xi|}{\ell+1} )^{2} } \, e_{2}
& & \text{for odd $\ell$}, \\
\tfrac{|\xi|-k}{\ell} \, e_{1} &\pm \sqrt{ k^{2} - ( \tfrac{|\xi|-k}{\ell} )^{2} } \, e_{2}
& & \text{for even $\ell$}. \\
\end{aligned}
\right.
\end{align}
Then, the auxiliary complex vectors $\zeta_{\ell,0}$ and $\big\{ \zeta_{\ell,j} : j \in \mathcal{U}_{\ell} \big\}$ can be constructed as follows.
\begin{enumerate}[label=\textbf{(\roman*)}]

\item For odd $\ell \in \mathbb{N}_{+}$,
\begin{align}\label{eqn:zeta_odd}
\zeta_{\ell,0} = \zeta_{\ell}^{-},
\qquad
\zeta_{\ell,1} = \zeta_{\ell}^{+},
\quad
\zeta_{\ell,2} = \zeta_{\ell}^{-},
\quad
\dots,
\quad
\zeta_{\ell,\ell} = \zeta_{\ell}^{+}.
\end{align}

\item For even $\ell \in \mathbb{N}_{+}$,
\begin{align}\label{eqn:zeta_even}
\zeta_{\ell,0} = k e_{1},
\qquad
\zeta_{\ell,1} = \zeta_{\ell}^{+},
\quad
\zeta_{\ell,2} = \zeta_{\ell}^{-},
\quad
\dots,
\quad
\zeta_{\ell,\ell-1} = \zeta_{\ell}^{+},
\quad
\zeta_{\ell,\ell} = \zeta_{\ell}^{-}.
\end{align}

\end{enumerate}
One can easily verify that these complex vectors satisfy the required conditions \eqref{eqn:zeta_condition}.

\begin{remark}\label{rmk:real_vec}
For each $\ell = 1,2,\dots,m$, we notice that, if $0 < |\xi| \leqslant (\ell+1) k$, then the vectors $\zeta_{\ell}^{\pm} \in \mathbb{R}^{n} \setminus \{0\}$.
\end{remark}

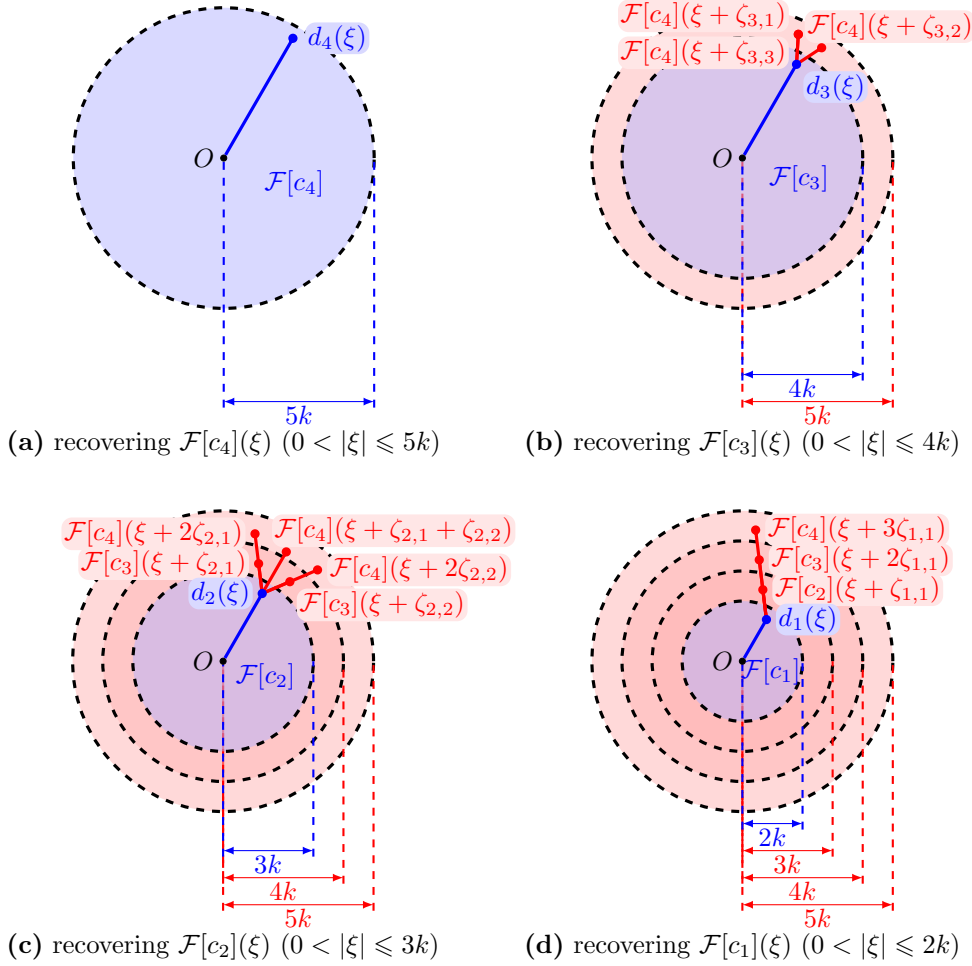
\begin{figure}[htbp]
\centering
\begin{tabular}{@{}l@{}l@{}l@{}}
\begin{tikzpicture}[>=latex,scale=0.4]
\node at (0,-9.5) {\textbf{(a)} recovering $\mathcal{F}[c_{4}](\xi)$ ($0 < |\xi| \leqslant 5 k$)};
\node[inner sep=1] (O) at (0,0) {};
\foreach \l/\c in {5/blue} {%
    \draw[very thick,dashed,fill={\c!30!white},fill opacity=0.5] (O) circle ({\l});
    \draw[thick,dashed,\c] ({\l},{-4.7-0.9*(\l-1)}) -- ({\l},0);
    \draw[<->,\c] (0,{-4.5-0.9*(\l-1)}) -- ++({\l},0) node[pos=0.5,below=-2] {${\l}k$};
    \draw[thick,dashed,\c] (0,{-4.7-0.9*(\l-1)}) -- (O);
}%
\draw[very thick,blue,fill=blue] (O) -- (60:4.6) circle (0.1) %
 node[pos=1,right=4,blue,fill=white!85!blue,inner sep=1,rounded corners] {$d_{4}(\xi)$};
\draw[fill=black] (O) circle (0.1) node[left] {$O$};
\node[blue] at (-20:2.5) {$\mathcal{F}[c_{4}]$};
\end{tikzpicture}
& &
\begin{tikzpicture}[>=latex,scale=0.4]
\node at (0,-9.5) {\textbf{(b)} recovering $\mathcal{F}[c_{3}](\xi)$ ($0 < |\xi| \leqslant 4 k$)};
\node[inner sep=1] (O) at (0,0) {};
\foreach \l/\c in {5/red,4/blue} {%
    \draw[very thick,dashed,fill={\c!30!white},fill opacity=0.5] (O) circle ({\l});
    \draw[thick,dashed,\c] ({\l},{-4.7-0.9*(\l-1)}) -- ({\l},0);
    \draw[<->,\c] (0,{-4.5-0.9*(\l-1)}) -- ++({\l},0) node[pos=0.5,below=-2] {${\l}k$};
    \draw[thick,dashed,\c] (0,{-4.7-0.9*(\l-1)}) -- (O);
}%
\draw[very thick,red,fill=red] (60:3.6) -- +({60+atan(1/2)}:1) circle (0.1) %
 node[pos=1, left=3,yshift=7,red,fill=white!90!red,inner sep=1,rounded corners] {$\mathcal{F}[c_{4}](\xi+\zeta_{3,1})$};
\draw[very thick,red,fill=red] (60:3.6) -- +({60-atan(1/2)}:1) circle (0.1) %
 node[pos=1,right,xshift=-7,yshift=8,red,fill=white!90!red,inner sep=1,rounded corners] {$\mathcal{F}[c_{4}](\xi+\zeta_{3,2})$};
\draw[very thick,red,fill=red] (60:3.6) -- +({60+atan(1/2)}:1) circle (0.1) %
 node[pos=1, left=3,yshift=-7,red,fill=white!90!red,inner sep=1,rounded corners] {$\mathcal{F}[c_{4}](\xi+\zeta_{3,3})$};
\draw[very thick,blue,fill=blue] (O) -- (60:3.6) circle (0.1) %
 node[pos=1,below right=2,blue,fill=white!85!blue,inner sep=1,rounded corners] {$d_{3}(\xi)$};
\draw[fill=black] (O) circle (0.1) node[left] {$O$};
\node[blue] at (-20:2.0) {$\mathcal{F}[c_{3}]$};
\end{tikzpicture}
\\[3ex]
\begin{tikzpicture}[>=latex,scale=0.4]
\node at (0,-9.5) {\textbf{(c)} recovering $\mathcal{F}[c_{2}](\xi)$ ($0 < |\xi| \leqslant 3 k$)};
\node[inner sep=1] (O) at (0,0) {};
\foreach \l/\c in {5/red,4/red,3/blue} {%
    \draw[very thick,dashed,fill={\c!30!white},fill opacity=0.5] (O) circle ({\l});
    \draw[thick,dashed,\c] ({\l},{-4.7-0.9*(\l-1)}) -- ({\l},0);
    \draw[<->,\c] (0,{-4.5-0.9*(\l-1)}) -- ++({\l},0) node[pos=0.5,below=-2] {${\l}k$};
    \draw[thick,dashed,\c] (0,{-4.7-0.9*(\l-1)}) -- (O);
}%
\draw[very thick,red,fill=red] (60:2.6) -- +({60}:1.6) circle (0.1) %
 node[pos=1,right,xshift=-6,yshift=8,red,fill=white!90!red,inner sep=1,rounded corners] {$\mathcal{F}[c_{4}](\xi+\zeta_{2,1}+\zeta_{2,2})$};
\draw[very thick,red,fill=red] (60:2.6) -- +({60-atan(3/4)}:2) circle (0.1) %
 node[pos=1,right=3,red,fill=white!90!red,inner sep=1,rounded corners] {$\mathcal{F}[c_{4}](\xi+2\zeta_{2,2})$};
\draw[very thick,red,fill=red] (60:2.6) -- +({60+atan(3/4)}:2) circle (0.1) %
 node[pos=1, left=3,red,fill=white!90!red,inner sep=1,rounded corners] {$\mathcal{F}[c_{4}](\xi+2\zeta_{2,1})$};
\draw[very thick,red,fill=red] (60:2.6) -- +({60-atan(3/4)}:1) circle (0.1) %
 node[pos=1,below right=2,red,fill=white!90!red,inner sep=1,rounded corners] {$\mathcal{F}[c_{3}](\xi+\zeta_{2,2})$};
\draw[very thick,red,fill=red] (60:2.6) -- +({60+atan(3/4)}:1) circle (0.1) %
 node[pos=1, left=3,red,fill=white!90!red,inner sep=1,rounded corners] {$\mathcal{F}[c_{3}](\xi+\zeta_{2,1})$};
\draw[very thick,blue,fill=blue] (O) -- (60:2.6) circle (0.1) %
 node[pos=1, left=3,blue,fill=white!85!blue,inner sep=1,rounded corners] {$d_{2}(\xi)$};
\draw[fill=black] (O) circle (0.1) node[left] {$O$};
\node[blue] at (-20:1.5) {$\mathcal{F}[c_{2}]$};
\end{tikzpicture}
& &
\begin{tikzpicture}[>=latex,scale=0.4]
\node at (0,-9.5) {\textbf{(d)} recovering $\mathcal{F}[c_{1}](\xi)$ ($0 < |\xi| \leqslant 2 k$)};
\node[inner sep=1] (O) at (0,0) {};
\foreach \l/\c in {5/red,4/red,3/red,2/blue} {%
    \draw[very thick,dashed,fill={\c!30!white},fill opacity=0.5] (O) circle ({\l});
    \draw[thick,dashed,\c] ({\l},{-4.7-0.9*(\l-1)}) -- ({\l},0);
    \draw[<->,\c] (0,{-4.5-0.9*(\l-1)}) -- ++({\l},0) node[pos=0.5,below=-2] {${\l}k$};
    \draw[thick,dashed,\c] (0,{-4.7-0.9*(\l-1)}) -- (O);
}%
\draw[very thick,red,fill=red] (60:1.6) -- +({60+atan(3/4)}:3) circle (0.1) %
 node[pos=1,right=3,red,fill=white!90!red,inner sep=1,rounded corners] {$\mathcal{F}[c_{4}](\xi+3\zeta_{1,1})$};
\draw[very thick,red,fill=red] (60:1.6) -- +({60+atan(3/4)}:2) circle (0.1) %
 node[pos=1,right=3,red,fill=white!90!red,inner sep=1,rounded corners] {$\mathcal{F}[c_{3}](\xi+2\zeta_{1,1})$};
\draw[very thick,red,fill=red] (60:1.6) -- +({60+atan(3/4)}:1) circle (0.1) %
 node[pos=1,right=3,red,fill=white!90!red,inner sep=1,rounded corners] {$\mathcal{F}[c_{2}](\xi+\zeta_{1,1})$};
\draw[very thick,blue,fill=blue] (O) -- (60:1.6) circle (0.1) %
 node[pos=1,right=3,blue,fill=white!85!blue,inner sep=1,rounded corners] {$d_{1}(\xi)$};
\draw[fill=black] (O) circle (0.1) node[left] {$O$};
\node[blue] at (-20:1.0) {$\mathcal{F}[c_{1}]$};
\end{tikzpicture}
\end{tabular}
\caption{The recursive reconstruction scheme \eqref{eqn:formula_matrix} for recovering the Fourier modes $\mathcal{F}[c_{\ell}](\xi)$ ($0 < |\xi| \leqslant (\ell+1) k$ and $\ell = 1,2,\dots,m$), with $m = 4$ and $n = 2$ as a representative case.}
\label{fig:operator_Q}
\end{figure}

\begin{remark}\label{rmk:scheme_fig}
Given the auxiliary vectors $\big\{ \zeta_{\ell,j} : j \in \mathcal{U}_{\ell} \big\}$ by \eqref{eqn:zeta_odd} or \eqref{eqn:zeta_even}, let
\begin{align*}
0 < |\xi| \leqslant (\ell+1)k
\quad \text{for\ } \mathcal{F}[c_{\ell}](\xi),
\end{align*}
then we can also obtain that the vectors $\zeta_{\ell,j} \in \mathbb{R}^{n} \setminus \{0\}$, and
\begin{align*}
0 < |\xi|
&\leqslant \big| \xi + {\textstyle\sum\limits_{j \in \mathcal{U}_{\ell}}} \alpha_{j} \zeta_{\ell,j} \big| \\
&\leqslant |\xi| + |\bm{\alpha}| k
= |\xi| + a k
\leqslant (\ell+a+1) k
\quad \text{for\ } \mathcal{F}[c_{\ell+a}]( \xi + {\textstyle\sum\limits_{j \in \mathcal{U}_{\ell}}} \alpha_{j} \zeta_{\ell,j} ),
\end{align*}
since $\zeta_{\ell,j} \cdot \zeta_{\ell,j} = k^{2}$ ($j \in \mathcal{U}_{\ell}$) in \eqref{eqn:zeta_condition} and $\bm{\alpha} = (\alpha_{1},\alpha_{2},\dots,\alpha_{\ell}) \in \mathcal{A}^{\ell}_{a}$ ($a \in \mathbb{N}_{0}$).

Hence, \textbf{Figure \ref{fig:operator_Q}} illustrates the recursive reconstruction scheme \eqref{eqn:formula_matrix} for recovering the Fourier modes $\mathcal{F}[c_{\ell}](\xi)$ ($0 < |\xi| \leqslant (\ell+1) k$ and $\ell = 1,2,\dots,m$).
Take $m = 4$, for example, we have that
\begin{enumerate}
\item[]
\begin{enumerate}[label=\textbf{Fig.\,\ref{fig:operator_Q}}\textbf{(\alph*)}]

\item $\mathcal{F}[c_{4}](\xi)$ ($0 < |\xi| \leqslant 5 k$) is recovered by
\begin{align*}
\mathcal{F}[c_{4}](\xi) = d_{4}(\xi).
\end{align*}

\item $\mathcal{F}[c_{3}](\xi)$ ($0 < |\xi| \leqslant 4 k$) is recovered by
\begin{align*}
\mathcal{F}[c_{3}](\xi) = d_{3}(\xi)
- 2 \big\{ \mathcal{F}[c_{4}](\xi+\zeta_{3,1}) + \mathcal{F}[c_{4}](\xi+\zeta_{3,2}) + \mathcal{F}[c_{4}](\xi+\zeta_{3,3}) \big\}
\end{align*}
with $\zeta_{3,1} = \zeta_{3,3} = \zeta_{3}^{+} \in \mathbb{R}^{n}$ and $\zeta_{3,2} = \zeta_{3}^{-} \in \mathbb{R}^{n}$.

\item $\mathcal{F}[c_{2}](\xi)$ ($0 < |\xi| \leqslant 3 k$) is recovered by
\begin{align*}
\mathcal{F}[c_{2}](\xi) = d_{2}(\xi)
&- 2 \big\{ \mathcal{F}[c_{4}](\xi+2\zeta_{2,1}) + \mathcal{F}[c_{4}](\xi+2\zeta_{2,2}) \big\} \\
&- 3 \mathcal{F}[c_{4}](\xi+\zeta_{2,1}+\zeta_{2,2}) \\
&- \tfrac{3}{2} \big\{ \mathcal{F}[c_{3}](\xi+\zeta_{2,1}) + \mathcal{F}[c_{3}](\xi+\zeta_{2,2}) \big\}
\end{align*}
with $\zeta_{2,1} = \zeta_{2}^{+} \in \mathbb{R}^{n}$ and $\zeta_{2,2} = \zeta_{2}^{-} \in \mathbb{R}^{n}$.

\item $\mathcal{F}[c_{1}](\xi)$ ($0 < |\xi| \leqslant 2 k$) is recovered by
\begin{align*}
\mathcal{F}[c_{1}](\xi) = d_{1}(\xi)
- \mathcal{F}[c_{4}](\xi+3\zeta_{1,1})
- \mathcal{F}[c_{3}](\xi+2\zeta_{1,1})
- \mathcal{F}[c_{2}](\xi+\zeta_{1,1})
\end{align*}
with $\zeta_{1,1} = \zeta_{1}^{+} \in \mathbb{R}^{n}$.

\end{enumerate}
\end{enumerate}
\end{remark}

As discussed in this section, particularly in \textbf{Lemma \ref{lmm:formulas}}, we demonstrate that by carefully choosing the boundary conditions for the nonlinear Helmholtz equation \eqref{eqn:probl_I} and its linearized system \eqref{eqn:probl_I0}--\eqref{eqn:probl_I1}, we can recursively decouple the multiple coefficients $\{ c_{\ell} \}_{\ell=1}^{m}$ in the polynomial $P_{m}(x,u)$.
This approach necessitates the use of subsequent theoretical analysis and constructive algorithms.

\section{Uniqueness and increasing stability for multi-coefficient inversion}\label{sec:theorem}

Under the definitions and assumptions in Section \ref{sec:setup}, we now present the main theorems establishing the uniqueness and increasing stability of multi-coefficient inversion.

Assume that each coefficient $c_{\ell} \in L^{\infty}(\Omega)$ with $\mathop{\mathrm{supp}}(c_{\ell}) \subset \Omega$ ($\ell = 1,2,\dots,m$), and $k^{2}$ is not a Dirichlet eigenvalue of $-\Delta$ in $\Omega$.
According to the well-posedness of the forward problem, we denote the operator norm of the linearized DtN map $\Lambda^{\prime}_{\bm{c}}$ in \eqref{eqn:linearDtN_I1} by
\begin{align}\label{eqn:op_norm}
\varepsilon := \| \Lambda^{\prime}_{\bm{c}} \|
= \sup\limits_{\| f \|_{\mathcal{D}} = 1} \| \Lambda^{\prime}_{\bm{c}}(f) \|_{\mathcal{N}},
\end{align}
by recalling that $\bm{c} = (c_{1},c_{2},\dots,c_{m})$, $\mathcal{D} = W^{2-\frac{1}{p},p}(\partial\Omega)$ and $\mathcal{N} = W^{1-\frac{1}{p},p}(\partial\Omega)$.

The first theorem verifies the uniqueness of the multi-coefficient identification problem.

\begin{theorem}\label{thm:Uniqueness}
Assume that each coefficient $c_{\ell} \in L^{\infty}(\Omega)$ with $\mathop{\mathrm{supp}}(c_{\ell}) \subset \Omega$ ($\ell = 1,2,\dots,m$), and $k^{2}$ is not a Dirichlet eigenvalue of $-\Delta$ in $\Omega$.
If $\Lambda^{\prime}_{\bm{c}_1} = \Lambda^{\prime}_{\bm{c}_2}$, then $\bm{c}_1 = \bm{c}_2$.
\end{theorem}

\begin{proof}
The proof follows directly from \textbf{Lemma \ref{lmm:formulas}}.
Specifically, if $\Lambda^{\prime}_{\bm{c}_1} = \Lambda^{\prime}_{\bm{c}_2}$, we deduce that
$\mathcal{F}[c_{m,1} - c_{m,2}](\xi) = 0$
for all $\xi \in \mathbb{R}^{n} \setminus \{0\}$, which implies $c_{m,1} = c_{m,2}$ due to the smallness assumption on the coefficients.
To complete the argument, we iteratively apply \eqref{eqn:formula_l} to show that $c_{\ell,1} = c_{\ell,2}$ for each $\ell \in \{ 1, 2, \dots, m-1 \}$.
\end{proof}

In numerical inversion algorithms, while uniqueness is important, stability is also crucial for achieving high-resolution reconstructions.
In the current work, by leveraging the wavenumber $k$, we establish an increasing stability estimate when $k$ is chosen sufficiently large.

The second theorem proves a Lipschitz-type stability result for recovering the Fourier modes $\mathcal{F}[c_{\ell}](\xi)$ with $0 < |\xi| \leqslant (\ell+1)k$ ($\ell = 1,2,\dots,m$) from the linearized DtN map $\Lambda^{\prime}_{\bm{c}}$.
Notably, a larger wavenumber $k$ enables stable recovery of more Fourier modes of each coefficient $c_{\ell}$ ($\ell = 1,2,\dots,m$).

\begin{theorem}\label{thm:Lipschitz}
Let $k > 1$, the Fourier transform of the coefficient $c_{\ell}$ at frequency $\xi \neq 0$ satisfies
\begin{enumerate}[label=\textbf{(\arabic*)}]

\item \textbf{Case $c_{m}$:}
\begin{align*}
\big|\mathcal{F}[c_{m}](\xi)\big| \leqslant \mathcal{C}(m,\Omega) k^{2m} \varepsilon
\quad
\text{for\ } |\xi| \leqslant (m+1) k.
\end{align*}

\item \textbf{Case $c_{\ell}$ with $0 < \ell < m$:}
\begin{align*}
\big|\mathcal{F}[c_{\ell}](\xi)\big| \leqslant \mathcal{C}(m,\ell,\Omega) k^{2m} \varepsilon
\quad
\text{for\ } |\xi| \leqslant (\ell+1) k.
\end{align*}

\end{enumerate}
\end{theorem}

\begin{proof}

For each $\ell = 1,2,\dots,m$, let a Fourier frequency be $\xi \in \mathbb{R}^{n}$ with $0 < |\xi| \leqslant (\ell+1)k$, then the solutions
\begin{align*}
\varphi_{\ell}(x) = \rme^{\img x \cdot \zeta_{\ell,0}},
\qquad
u_{\ell,j}^{\sss{0}}(x) = \rme^{\img x \cdot \zeta_{\ell,j}},
\quad j \in \mathcal{U}_{\ell},
\end{align*}
defined in \eqref{eqn:CE_solution} are indeed the plane waves, since $\zeta_{\ell,0}$ and $\zeta_{\ell,j}$ are chosen as real vectors; see \textbf{Remark \ref{rmk:real_vec}} and \textbf{Remark \ref{rmk:scheme_fig}}.
Let $\mathcal{C}(\Omega)$ be a constant depending on $\Omega$, we first obtain that
\begin{align*}
\| \varphi_{\ell} |_{\partial\Omega} \|_{L^{\infty}(\partial\Omega)} = 1,
\qquad
\| u_{\ell,j}^{\sss{0}} \|_{W^{2,p}(\Omega)}
\leqslant \mathcal{C}(\Omega) k^{2},
\end{align*}
and
\begin{align*}
\| u_{\ell,\mathcal{S}}^{\sss{0}} \|_{W^{2,p}(\Omega)}
\leqslant {\textstyle\sum\limits_{j \in \mathcal{S}}} \| u_{\ell,j}^{\sss{0}} \|_{W^{2,p}(\Omega)}
\leqslant \mathcal{C}(\Omega) k^{2} |\mathcal{S}|,
\end{align*}
for any $j \in \mathcal{U}_{\ell}$ and non-empty subset $\mathcal{S}$ of $\mathcal{U}_{\ell}$.
Hence, by following the proof in \cite[Section 3.1]{ZLX2024},
we have
\begin{equation}\label{eqn:est_dl}
\begin{aligned}
|d_{\ell}|
&= \left| \frac{1}{\,\ell\,!\,} \, {\textstyle\sum\limits_{\emptyset \subsetneqq \mathcal{S} \subseteq \mathcal{U}_{\ell}}} (-1)^{|\mathcal{U}_{\ell} \setminus \mathcal{S}|} \int_{\partial\Omega} f_{\ell,0} \, \Lambda^{\prime}_{\bm{c}} (f_{\ell,\mathcal{S}}) \,\rmd s_{x} \right| \\
&\leqslant \frac{1}{\,\ell\,!\,} \, {\textstyle\sum\limits_{\emptyset \subsetneqq \mathcal{S} \subseteq \mathcal{U}_{\ell}}} \left| \int_{\partial\Omega} f_{\ell,0} \, \Lambda^{\prime}_{\bm{c}} (f_{\ell,\mathcal{S}}) \,\rmd s_{x} \right| \\
&\leqslant \mathcal{C}(\Omega) \| f_{\ell,0} \|_{L^{\infty}(\partial\Omega)} \frac{1}{\,\ell\,!\,} \, {\textstyle\sum\limits_{\emptyset \subsetneqq \mathcal{S} \subseteq \mathcal{U}_{\ell}}} \| \Lambda^{\prime}_{\bm{c}}(f_{\ell,\mathcal{S}}) \|_{\mathcal{N}} \\
&\leqslant \mathcal{C}(\Omega) \| f_{\ell,0} \|_{L^{\infty}(\partial\Omega)} \Big( {\textstyle\sum\limits_{\ell'=1}^{m}} \frac{1}{\,\ell\,!\,} \, {\textstyle\sum\limits_{\emptyset \subsetneqq \mathcal{S} \subseteq \mathcal{U}_{\ell}}} \| f_{\ell,\mathcal{S}} \|_{\mathcal{D}}^{\ell'} \Big) \varepsilon \\
&\leqslant \mathcal{C}(m,\ell,\Omega) \| \varphi_{\ell} |_{\partial\Omega} \|_{L^{\infty}(\partial\Omega)} \Big( {\textstyle\sum\limits_{\ell'=1}^{m}} \frac{1}{\,\ell\,!\,} \, {\textstyle\sum\limits_{\emptyset \subsetneqq \mathcal{S} \subseteq \mathcal{U}_{\ell}}} \| u_{\ell,\mathcal{S}}^{\sss{0}} \|_{W^{2,p}(\Omega)}^{\ell'} \Big) \varepsilon.
\end{aligned}
\end{equation}
Here, $\mathcal{C}(m,\ell,\Omega)$ is a constant depending on $m$, $\ell$ and $\Omega$.
Moreover, while $0 < |\xi| \leqslant (\ell+1)k$ and $k > 1$,
\begin{align*}
|d_{\ell}|
&\leqslant \mathcal{C}(m,\ell,\Omega) \Big( {\textstyle\sum\limits_{\ell'=1}^{m}} \frac{1}{\,\ell\,!\,} \, {\textstyle\sum\limits_{\emptyset \subsetneqq \mathcal{S} \subseteq \mathcal{U}_{\ell}}} k^{2\ell'} |\mathcal{S}|^{\ell'} \Big) \varepsilon \\
&\leqslant \mathcal{C}(m,\ell,\Omega) k^{2m} \Big( {\textstyle\sum\limits_{\ell'=1}^{m}} \frac{1}{\,\ell\,!\,} \, {\textstyle\sum\limits_{j=1}^{\ell}} {\textstyle\binom{\ell}{j}} j^{\ell'} \Big) \varepsilon \\
&\leqslant \mathcal{C}(m,\ell,\Omega) k^{2m} \varepsilon.
\end{align*}

\begin{enumerate}[label=\textbf{(\arabic*)}]

\item Case $\ell = m$: By \textbf{Lemma \ref{lmm:formulas}}, we then obtain that
\begin{align*}
\big|\mathcal{F}[c_{m}](\xi)\big|
= |d_{m}|
\leqslant \mathcal{C}(m,\Omega) k^{2m} \varepsilon
\quad
\text{for\ } 0 < |\xi| \leqslant (m+1) k.
\end{align*}

\item Case $0 < \ell < m$: The proof will be completed by induction.
Assuming that, for each $a = 1,2,\dots,m-\ell$, it holds
\begin{align*}
\big|\mathcal{F}[c_{\ell+a}](\widetilde{\xi})\big|
\leqslant \mathcal{C}(m,\ell+a,\Omega) k^{2m} \varepsilon
\quad
\text{for\ } 0 < |\widetilde{\xi}| \leqslant (\ell+a+1)k.
\end{align*}
Thus, for $0 < |\xi| \leqslant (\ell+1)k$ and $\bm{\alpha} = (\alpha_{1},\alpha_{2},\dots,\alpha_{\ell}) \in \mathcal{A}^{\ell}_{a}$, by noticing
\begin{align*}
0 < \big| \xi + {\textstyle\sum\limits_{j \in \mathcal{U}_{\ell}}} \alpha_{j} \zeta_{\ell,j} \big|
\leqslant (\ell+a+1)k;
\end{align*}
see \textbf{Remark \ref{rmk:scheme_fig}}, we finally obtain that
\begin{align*}
\big|\mathcal{F}[c_{\ell}](\xi)\big|
&\leqslant |d_{\ell}|
+ {\textstyle\sum\limits_{a = 1}^{m-\ell}} \big|\mathcal{F}[c_{\ell+a}Q_{\ell,a}](\xi)\big| \\
&\leqslant |d_{\ell}|
+ {\textstyle\sum\limits_{a = 1}^{m-\ell}} \frac{1}{\,\ell\,!\,} {\textstyle\sum\limits_{\bm{\alpha} \in \mathcal{A}^{\ell}_{a}}} {\textstyle\binom{\ell+a}{\bm{1}+\bm{\alpha}}} \big|\mathcal{F}[c_{\ell+a}]( \xi + {\textstyle\sum\limits_{j \in \mathcal{U}_{\ell}}} \alpha_{j} \zeta_{\ell,j} )\big| \\
&\leqslant \mathcal{C}(m,\ell,\Omega) k^{2m} \varepsilon
+ {\textstyle\sum\limits_{a = 1}^{m-\ell}} \frac{1}{\,\ell\,!\,} {\textstyle\sum\limits_{\bm{\alpha} \in \mathcal{A}^{\ell}_{a}}} {\textstyle\binom{\ell+a}{\bm{1}+\bm{\alpha}}} \, \mathcal{C}(m,\ell+a,\Omega) k^{2m} \varepsilon \\
&\leqslant \mathcal{C}(m,\ell,\Omega) k^{2m} \varepsilon.
\end{align*}

\end{enumerate}
Therefore, we conclude the proof of \textbf{Theorem \ref{thm:Lipschitz}}.
\end{proof}

\begin{remark}
There is a remark for the operator norm of $\Lambda^{\prime}_{\bm{c}}$ in \eqref{eqn:op_norm} and the estimate of the data $d_{\ell}$ in \eqref{eqn:est_dl}.
Let $p = \frac{n}{2}+1$, then $W^{2,p}(\Omega) \subset C^{0,s}(\overline{\Omega}) \subset L^{\infty}(\Omega)$ for some $s \in (0,1)$.
Hence, for the linearized systems ($I^{\sss{0}}$) and ($I^{\sss{1}}$) in \eqref{eqn:probl_I0}--\eqref{eqn:probl_I1}, we have
\begin{align*}
\| u^{\sss{1}} \|_{W^{2,p}(\Omega)}
&\leqslant \mathcal{C}(\Omega,k) \| P_{m}(x,u^{\sss{0}}) \|_{L^{p}(\Omega)} \\
&\leqslant \mathcal{C}(\Omega,k) \| \bm{c} \|_{(L^{\infty}(\Omega))^{m}} {\textstyle\sum\limits_{\ell=1}^{m}} \| u^{\sss{0}} \|_{L^{\infty}(\Omega)}^{\ell} \\
&\leqslant \mathcal{C}(\Omega,k) \| \bm{c} \|_{(L^{\infty}(\Omega))^{m}} {\textstyle\sum\limits_{\ell=1}^{m}} \| u^{\sss{0}} \|_{W^{2,p}(\Omega)}^{\ell} \\
&\leqslant \mathcal{C}(m,\Omega,k) \| \bm{c} \|_{(L^{\infty}(\Omega))^{m}} {\textstyle\sum\limits_{\ell=1}^{m}} \| f \|_{\mathcal{D}}^{\ell},
\end{align*}
which leads to
\begin{align*}
\| \Lambda^{\prime}_{\bm{c}}(f) \|_{\mathcal{N}}
= \| \partial_{\nu} u^{\sss{1}} |_{\partial\Omega} \|_{\mathcal{N}}
\leqslant \mathcal{C}(m,\Omega,k) \| \bm{c} \|_{(L^{\infty}(\Omega))^{m}} {\textstyle\sum\limits_{\ell=1}^{m}} \| f \|_{\mathcal{D}}^{\ell}.
\end{align*}
It also shows the linear dependence of $\Lambda^{\prime}_{\bm{c}}$ on the coefficient vector $\bm{c}$ versus its nonlinear dependence on the boundary data $f$.
\end{remark}

The final theorem establishes an increasing stability estimate for the coefficient vector $\bm{c}$ using the linearized DtN map $\Lambda^{\prime}_{\bm{c}}$. 
This result generalizes earlier works in the linear case \cite[Theorem 2.1]{ILX2020}, the quadratic case \cite[Theorem 3.1]{LSX2022} and the integer power-type case \cite[Theorem 3.3]{ZLX2024}.

\begin{theorem}\label{thm:main}
If $\varepsilon < 1$ and $\| \bm{c} \|_{(H^{1}(\Omega))^{m}} \leqslant M$, the following stability estimate holds
\begin{align*}
\| \bm{c} \|_{(L^{2}(\Omega))^{m}}^{2}
\leqslant \mathcal{C}(m,\Omega) \Big( k^{4m+n} \varepsilon^{2} + \frac{M^{2}}{1+3k^{2}+E^{2}} \Big),
\end{align*}
for sufficiently large wavenumber $k$, that is $k > E := -\ln \varepsilon > 0$.
\end{theorem}

\begin{proof}
Following the proof in \cite[Section 3.1]{ZLX2024}, we outline the proof as follows.

Under a priori assumption $\| \bm{c} \|_{(H^{1}(\Omega))^{m}} \leqslant M$ and $\mathop{\mathrm{supp}}(c_{\ell}) \subset \Omega$ ($\ell = 1,2,\dots,m$), we first derive an $L^{2}$-norm estimate for $\bm{c} = (c_{1},c_{2},\dots,c_{m})$ in high Fourier frequency regime $\big\{ \xi \in \mathbb{R}^{n} : |\xi| > \rho \big\}$ with $\rho > 0$,
\begin{align*}
\int_{|\xi| > \rho} \big|\mathcal{F}[c_{\ell}](\xi)\big|^{2} \,\rmd \xi
&\leqslant \frac{1}{1+\rho^{2}} \int_{\mathbb{R}^{n}} (1+|\xi|^{2}) \big|\mathcal{F}[c_{\ell}](\xi)\big|^{2} \,\rmd \xi \\
&\leqslant \mathcal{C}(\Omega) \frac{\| c_{\ell} \|_{H^{1}(\Omega)}^{2}}{1+\rho^{2}},
\quad
\ell = 1,2,\dots,m.
\end{align*}

While the wavenumber $k$ is large enough (i.e. $k > E$), let $\rho_{\ell} := (\ell+1)k$, we have
\begin{align*}
\| \bm{c} \|_{(L^{2}(\Omega))^{m}}^{2}
&= {\textstyle\sum\limits_{\ell=1}^{m}} \| c_{\ell} \|_{L^{2}(\Omega)}^{2} \\
&= {\textstyle\sum\limits_{\ell=1}^{m}} \int_{|\xi| \leqslant (\ell+1)k} \big|\mathcal{F}[c_{\ell}](\xi)\big|^{2} \,\rmd \xi
+ {\textstyle\sum\limits_{\ell=1}^{m}} \int_{|\xi| > (\ell+1)k} \big|\mathcal{F}[c_{\ell}](\xi)\big|^{2} \,\rmd \xi \\
&\leqslant {\textstyle\sum\limits_{\ell=1}^{m}} \mathcal{C}(m,\ell,\Omega) (\ell+1)^{n} k^{n} k^{4m} \varepsilon^{2}
+ \mathcal{C}(\Omega) {\textstyle\sum\limits_{\ell=1}^{m}} \frac{\| c_{\ell} \|_{H^{1}(\Omega)}^{2}}{1+\rho_{\ell}^{2}} \\
&\leqslant \mathcal{C}(m,\Omega) k^{4m+n} \varepsilon^{2} + \mathcal{C}(\Omega) \frac{\| \bm{c} \|_{(H^{1}(\Omega))^{m}}^{2}}{1+4k^{2}} \\
&\leqslant \mathcal{C}(m,\Omega) \Big( k^{4m+n} \varepsilon^{2} + \frac{M^{2}}{1+3k^{2}+E^{2}} \Big),
\end{align*}
by the estimate of Fourier modes $\mathcal{F}[c_{\ell}](\xi)$ for $|\xi| \leqslant (\ell+1)k$; see \textbf{Theorem \ref{thm:Lipschitz}}.
\end{proof}

\begin{remark}
As yet, it should be noted that no H\"{o}lder-type estimate is established for the multi-coefficient identification problem with a small wavenumber (i.e. $k \leqslant E$), in contrast to the single-coefficient case studied in \cite[Theorem 3.3]{ZLX2024}.
This phenomenon can be understood as representing the additional challenges inherent in multi-coefficient inversion for nonlinear Helmholtz equations.

The primary challenges arise in the intermediate Fourier frequency regime $\big\{ \xi \in \mathbb{R}^{n} : (\ell+1)k < |\xi| \leqslant \rho \big\}$ when $k \leqslant E$.
Notice that the Fourier modes are updated iteratively by
\begin{align*}
\mathcal{F}[c_{\ell}](\xi)
&= d_{\ell}(\xi)
- {\textstyle\sum\limits_{a = 1}^{m-\ell}} \frac{1}{\,\ell\,!\,} {\textstyle\sum\limits_{\bm{\alpha} \in \mathcal{A}^{\ell}_{a}}} {\textstyle\binom{\ell+a}{\bm{1}+\bm{\alpha}}} \, \mathcal{F}[c_{\ell+a}]( \xi + {\textstyle\sum\limits_{j \in \mathcal{U}_{\ell}}} \alpha_{j} \zeta_{\ell,j} ),
\end{align*}
and the vectors $\zeta_{\ell,j}$ and $\xi + {\textstyle\sum_{j \in \mathcal{U}_{\ell}}} \alpha_{j} \zeta_{\ell,j}$ are complex while $(\ell+1)k < |\xi| \leqslant \rho$.
Thus, by an argument analogous to the proof of \textbf{Theorem \ref{thm:Lipschitz}}, the stability result for Fourier mode $\mathcal{F}[c_{\ell}](\xi)$ in this regime is not clear.

Fortunately, for each $\ell$, the truncated Fourier modes $\mathcal{F}[c_{\ell}](\xi)$ restricted to the regime $\big\{ \xi \in \mathbb{R}^{n} : 0 < |\xi| \leqslant (\ell+1)k \big\}$ are sufficient to numerically reconstruct the coefficient $c_{\ell}$ owing to the Lipschitz-type stability established in \textbf{Theorem \ref{thm:Lipschitz}}.
\end{remark}

%

\section{Recursive reconstruction algorithm and numerical examples}\label{sec:algorithm}

In this section, we present an explicit recursive reconstruction algorithm and validate its performance through numerical experiments.

\subsection{The recursive reconstruction algorithm}

We present our recursive reconstruction algorithm in \textbf{Algorithm \ref{alg:main}}, which extends previous approaches from \cite{ILX2020, LSX2022, ZLX2024} to handle polynomial-type nonlinearities with multiple coefficients.

\begin{algorithm}[h]
\caption{Recursive reconstruction algorithm for multi-coefficient inversion in nonlinear Helmholtz equation with polynomial-type nonlinearity}\label{alg:main}
\begin{algorithmic}[1]
\Require the wavenumber $k > 1$, the degree $m \geqslant 2$, the open bounded domain $\Omega$%
\Ensure the coefficients of polynomial $P_{m}(x,u)$, i.e. $\bm{c} = (c_{1},c_{2},\dots,c_{m})$ in $\Omega$
\For{$\ell = m,m-1,\dots,1$}
\For{each $\xi \in \mathbb{R}^{n} \setminus \{0\}$ satisfies $|\xi| \leqslant (\ell+1)k$}
\State Define the auxiliary vectors $\zeta_{\ell}^{+}$ and $\zeta_{\ell}^{-}$ as \eqref{eqn:zeta_pm};
\State Construct the vectors $\zeta_{\ell,0}$ and $\zeta_{\ell,j}$ ($j \in \mathcal{U}_{\ell}$) by \eqref{eqn:zeta_odd} or \eqref{eqn:zeta_even};
\State Choose the solutions $\varphi_{\ell} = \exp\{\img x \cdot \zeta_{\ell,0}\}$ and $u_{\ell,j}^{\sss{0}} = \exp\{\img x \cdot \zeta_{\ell,j}\}$;
\State Give the Dirichlet boundary data $f_{\ell,0} = \varphi_{\ell} |_{\partial\Omega}$ and $f_{\ell,j} = u_{\ell,j}^{\sss{0}} |_{\partial\Omega}$;
\For{each non-empty subset $\mathcal{S} \subset \mathcal{U}_{\ell}$}
\State Give the combined solution $u_{\ell,\mathcal{S}}^{\sss{0}} = {\textstyle\sum_{j \in \mathcal{S}}} u_{\ell,j}^{\sss{0}}$,
\Statex \qquad\qquad\qquad~ and its Dirichlet boundary data $f_{\ell,\mathcal{S}} = {\textstyle\sum_{j \in \mathcal{S}}} f_{\ell,j}$;
\label{alg:Dirichlet}
\State Measure the Neumann boundary data $\partial_{\nu} u_{\ell,\mathcal{S}} |_{\partial\Omega}$ of the problem \eqref{eqn:probl_I},
\Statex \qquad\qquad\qquad~ while the Dirichlet boundary data $u_{\ell,\mathcal{S}} |_{\partial\Omega} = f_{\ell,\mathcal{S}}$ are given;
\label{alg:Neumann}
\State Approximate the linearized Neumann boundary data $\partial_{\nu} u_{\ell,\mathcal{S}}^{\sss{1}} |_{\partial\Omega}$ by
\Statex \qquad\qquad\qquad~ $\Lambda^{\prime}_{\bm{c}} (f_{\ell,\mathcal{S}}) = \partial_{\nu} u_{\ell,\mathcal{S}}^{\sss{1}} |_{\partial\Omega} ~\approx~ ( \partial_{\nu} u_{\ell,\mathcal{S}} - \partial_{\nu} u_{\ell,\mathcal{S}}^{\sss{0}} )|_{\partial\Omega}$;
\label{alg:comment}
\EndFor
\State Compute $d_{\ell}(\xi) = \frac{1}{\,\ell\,!\,} \, {\textstyle\sum_{\emptyset \subsetneqq \mathcal{S} \subseteq \mathcal{U}_{\ell}}} (-1)^{|\mathcal{U}_{\ell} \setminus \mathcal{S}|} \int_{\partial\Omega} f_{\ell,0} \, \Lambda^{\prime}_{\bm{c}} (f_{\ell,\mathcal{S}}) \,\rmd s_{x}$;
\label{alg:Fc_l_start}
\If{$\ell = m$}
\State $\mathcal{F}[c_{m}](\xi) = d_{m}(\xi)$;
%
\Else
\State $\mathcal{F}[c_{\ell}](\xi) = d_{\ell}(\xi) - {\textstyle\sum_{a = 1}^{m-\ell}} \frac{1}{\,\ell\,!\,} {\textstyle\sum_{\bm{\alpha} \in \mathcal{A}^{\ell}_{a}}} {\binom{\ell+a}{\bm{1}+\bm{\alpha}}} \, \mathcal{F}[c_{\ell+a}]( \xi + {\textstyle\sum_{j \in \mathcal{U}_{\ell}}} \alpha_{j} \zeta_{\ell,j} )$;
%
\EndIf
\label{alg:Fc_l_end}
\EndFor
\State Compute $c_{\ell}$ by the inverse Fourier transform of $\mathcal{F}[c_{\ell}]$;
\label{alg:c_l}
\EndFor
\end{algorithmic}
\end{algorithm}

The proposed \textbf{Algorithm \ref{alg:main}} contains three iterations.
\begin{enumerate}[label=\textbf{(\arabic*)}]

\item \textbf{The iteration of $\ell$} (\textsf{Line 1--20}):
The reconstruction proceeds recursively via backward iteration from $\ell = m$ to $\ell = 1$.
At the iteration of $\ell$, we successively recover the coefficient $c_{\ell}$ via the inverse Fourier transform of its Fourier mode $\mathcal{F}[c_{\ell}]$; see \textsf{Line \ref{alg:c_l}}.

\item \textbf{The iteration of $\xi$} (\textsf{Line 2--18}):
For each iteration of $\ell$, we compute the (truncated) Fourier mode $\mathcal{F}[c_{\ell}](\xi)$ at a bounded domain $0 < |\xi| \leqslant (\ell+1)k$, by using previously reconstructed Fourier modes $\big\{ \mathcal{F}[c_{\ell+a}](\widetilde{\xi}) : 0 < |\widetilde{\xi}| \leqslant (\ell+a+1)k,~ a = 1,2,\dots,m-\ell \big\}$ and the data $d_{\ell}(\xi)$; see \textsf{Line \ref{alg:Fc_l_start}--\ref{alg:Fc_l_end}}.

\item \textbf{The iteration of $\mathcal{S}$} (\textsf{Line 7--11}):
For each non-empty subset $\mathcal{S}$ of $\mathcal{U}_{\ell}$, the Neumann boundary data $\partial_{\nu} u_{\ell,\mathcal{S}} |_{\partial\Omega}$ is measured on $\partial\Omega$, by solving the original problem \eqref{eqn:probl_I} while the Dirichlet boundary data $u_{\ell,\mathcal{S}} |_{\partial\Omega} = f_{\ell,\mathcal{S}} = u_{\ell,\mathcal{S}}^{\sss{0}} |_{\partial\Omega}$ are given; see \textsf{Line \ref{alg:Dirichlet}--\ref{alg:Neumann}}.
Then, as mentioned in \eqref{eqn:Neumann_approximate}, we can approximate the linearized Neumann boundary data $\partial_{\nu} u_{\ell,\mathcal{S}}^{\sss{1}} |_{\partial\Omega}$ by
\begin{align}\label{eqn:approx_Neumann}
\Lambda^{\prime}_{\bm{c}} (f_{\ell,\mathcal{S}})
= \partial_{\nu} u_{\ell,\mathcal{S}}^{\sss{1}} |_{\partial\Omega}
~\approx~ ( \partial_{\nu} u_{\ell,\mathcal{S}} - \partial_{\nu} u_{\ell,\mathcal{S}}^{\sss{0}} )|_{\partial\Omega};
\end{align}
see \textsf{Line \ref{alg:comment}}.
So far, the data $d_{\ell}$ can be computed in \textsf{Line \ref{alg:Fc_l_start}}.

\end{enumerate}

To better understand \textbf{Algorithm \ref{alg:main}}, we provide some comments listed as follows.

In \textsf{Line \ref{alg:comment}}, we can NOT utilize the linearized Neumann boundary data $\partial_{\nu} u_{\ell,\mathcal{S}}^{\sss{1}} |_{\partial\Omega}$ directly, since it depends on the unknown $\bm{c} = (c_{1},c_{2},\dots,c_{m})$ referring to the linearized problem \eqref{eqn:probl_I1}.
Fortunately, when the well-posedness of the original problem \eqref{eqn:probl_I} is guaranteed, the linearized Neumann boundary data $\Lambda^{\prime}_{\bm{c}} (f_{\ell,\mathcal{S}}) = \partial_{\nu} u_{\ell,\mathcal{S}}^{\sss{1}} |_{\partial\Omega}$ can be approximated well by the difference between the Neumann boundary data $\partial_{\nu} u_{\ell,\mathcal{S}} |_{\partial\Omega}$ of the original problem \eqref{eqn:probl_I} and $\partial_{\nu} u_{\ell,\mathcal{S}}^{\sss{0}} |_{\partial\Omega}$ of the unperturbed problem \eqref{eqn:probl_I0}; see \eqref{eqn:Neumann_approximate} and \eqref{eqn:approx_Neumann}.

In \textsf{Line \ref{alg:Neumann}}, since $P_{m}(x,u)$ is a polynomial in $u$ of degree $m$, we have to solve the nonlinear Helmholtz equation \eqref{eqn:probl_I} several times at each iteration of $\ell$, which causes considerably high computational costs.
Let $N_{\ell}$ denote the number of mesh points of the Fourier frequency regime $\big\{ \xi \in \mathbb{R}^{n} : |\xi| \leqslant (\ell+1)k \big\}$.
For each iteration of $\ell$, the forward problem must be solved $(2^{\ell}-1) N_{\ell}$ times.
If $\mathcal{F}[c_{\ell+a}]( \xi + {\textstyle\sum_{j \in \mathcal{U}_{\ell}}} \alpha_{j} \zeta_{\ell,j} )$ can be approximated by interpolating $\mathcal{F}[c_{\ell+a}](\xi)$ at mesh points, then our algorithm requires solving the forward problem a total of ${\textstyle\sum_{\ell=1}^{m}} (2^{\ell}-1) N_{\ell}$ times.
We further refer to \cite{FT2005, XB2010, YL2017, WZ2018, JLWZ2022, JWXZ2024} for discussion on varies efficient numerical schemes solving nonlinear Helmholtz equations.

\subsection{The numerical experiments}

Finally, we provide several numerical examples for the two-dimensional case ($n = 2$) to validate our algorithm (\textbf{Algorithm \ref{alg:main}}).
Here, the wavenumber $k = 20$ is fixed for all examples, and the open bounded domain $\Omega \subset \mathbb{R}^{2}$ is given by a disk of radius $0.5$ centered at the origin.
To prevent inverse crime in numerical examples, we employ distinct discretizations on the same computational square domain $[{-0.5},{+0.5}]^{2}$: a fine grid ($201 \times 201$ equidistant points) for the forward problem, and a coarse grid ($191 \times 191$ equidistant points) for the inversion.

\subsubsection{The example for the degree $m = 2$}

\begin{figure}[htbp]
\centering
\begin{tabular}{@{}ccc@{}}
\includegraphics[width=0.33\textwidth]{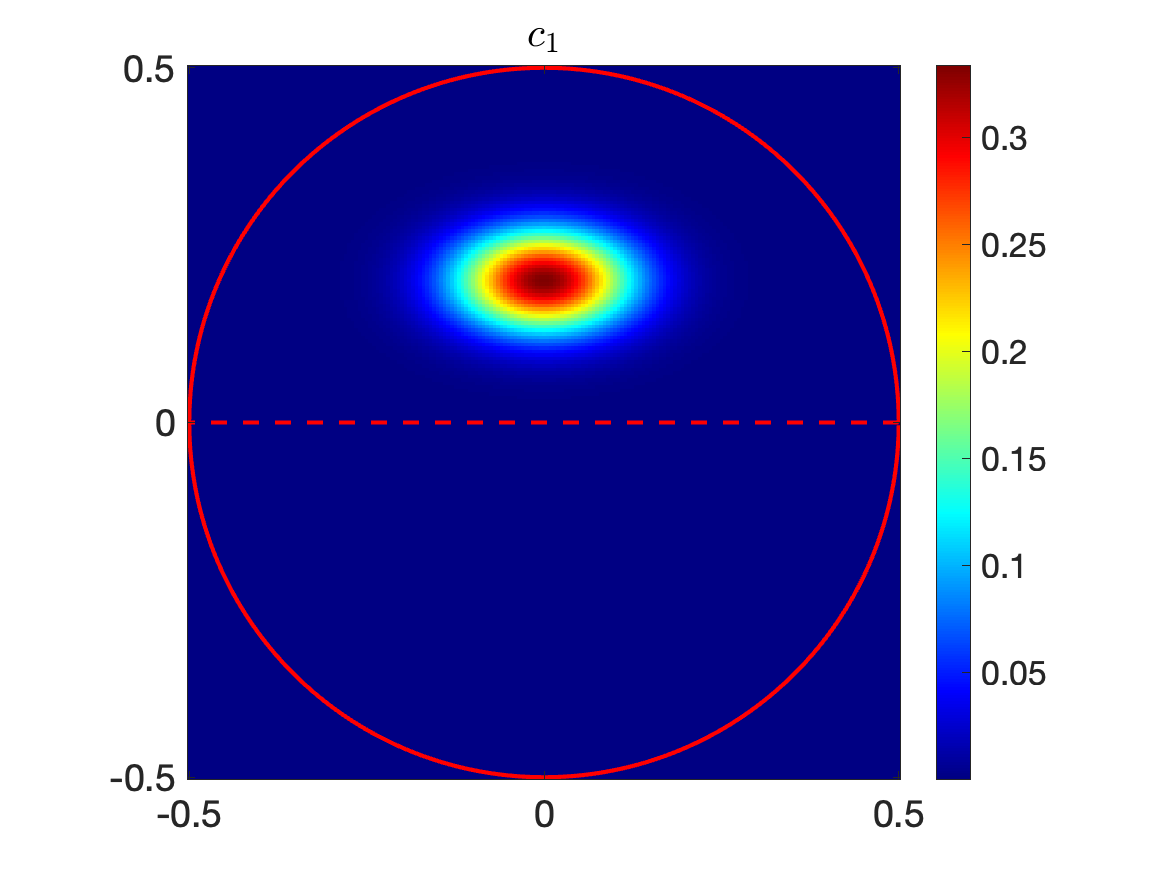} &
\includegraphics[width=0.3\textwidth]{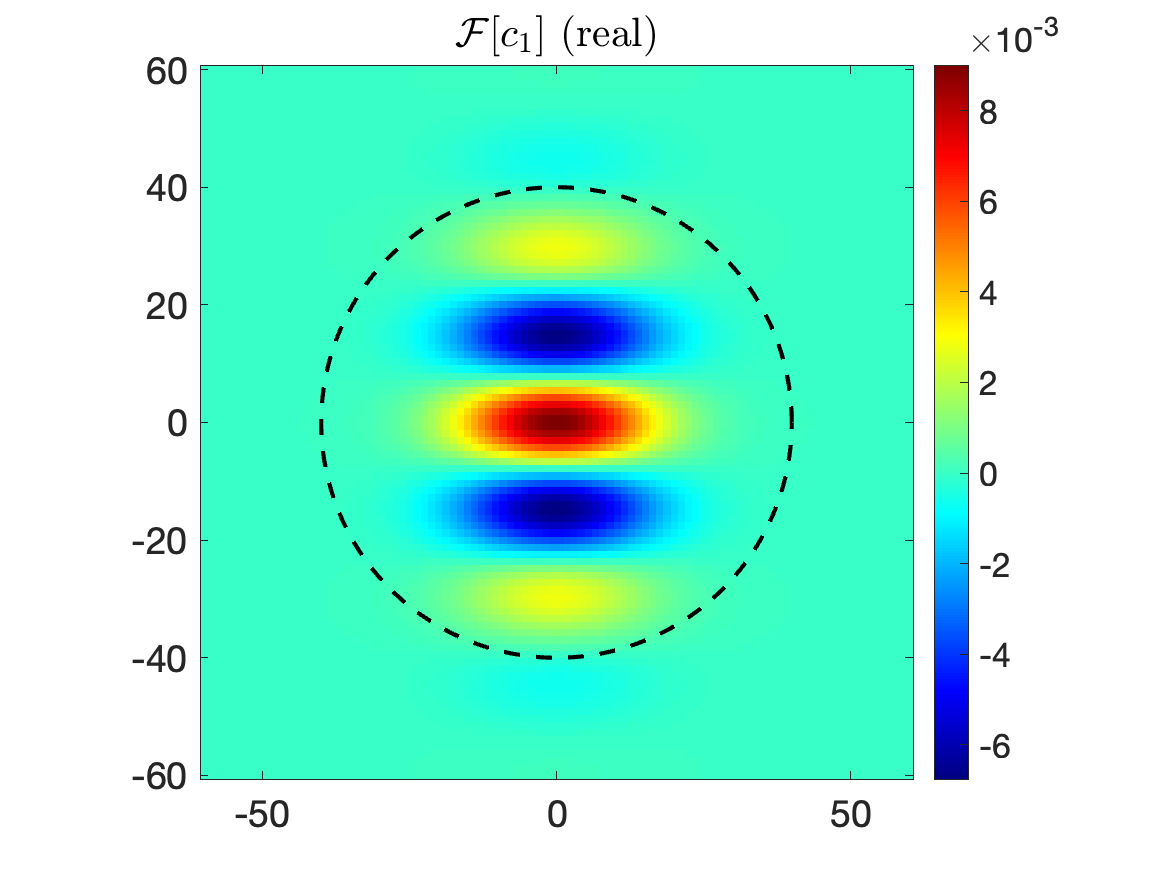} &
\includegraphics[width=0.3\textwidth]{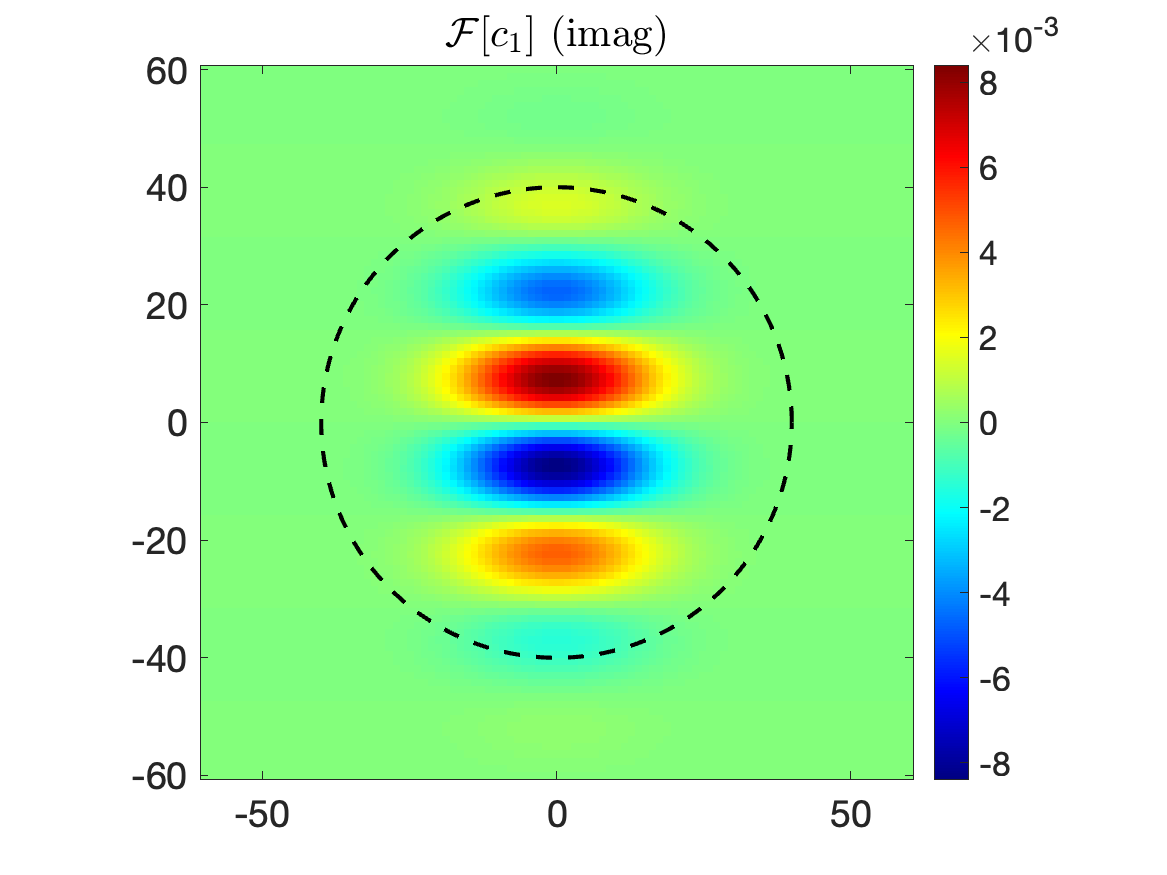} \\
\textbf{(1)} $c_{1}(x)$ &
\textbf{(1-a)} $\mathcal{F}[c_{1}](\xi)$ (real) &
\textbf{(1-b)} $\mathcal{F}[c_{1}](\xi)$ (imaginary) \\[3ex]
\includegraphics[width=0.33\textwidth]{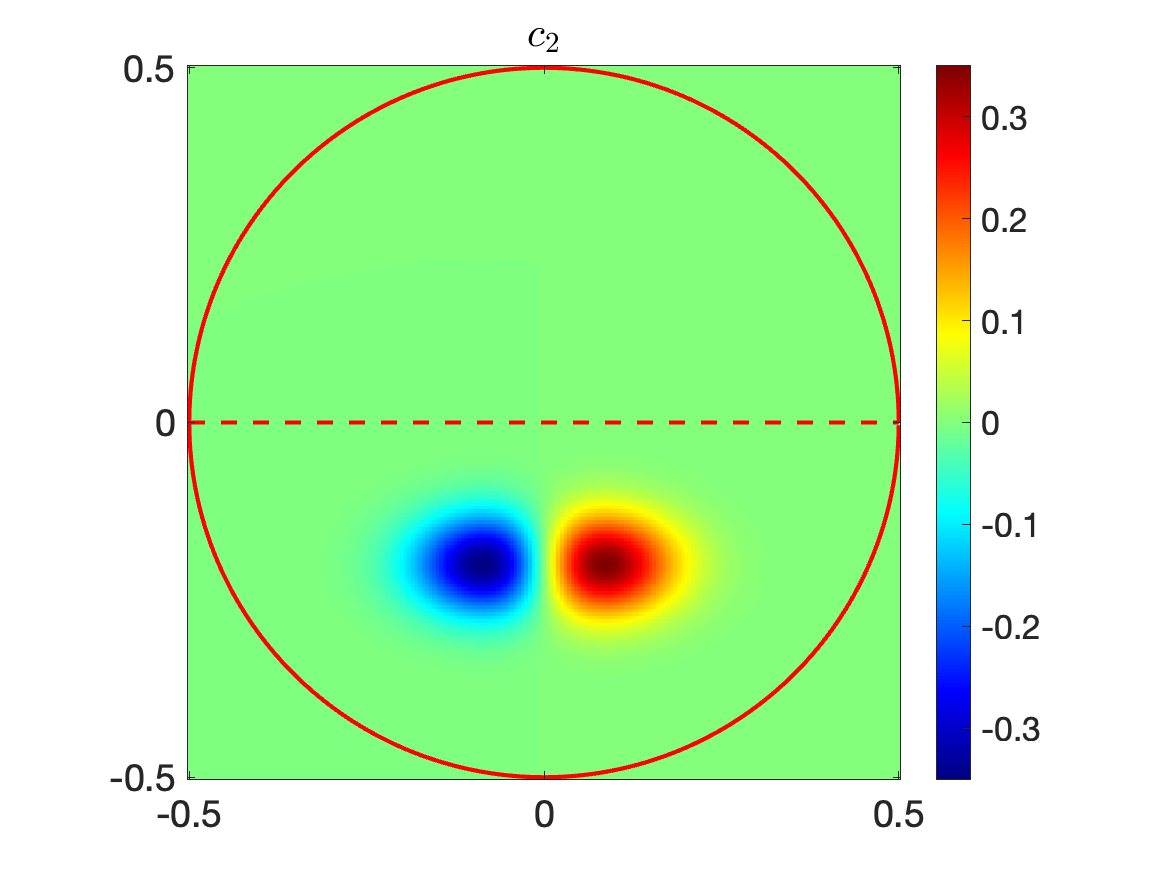} &
\includegraphics[width=0.3\textwidth]{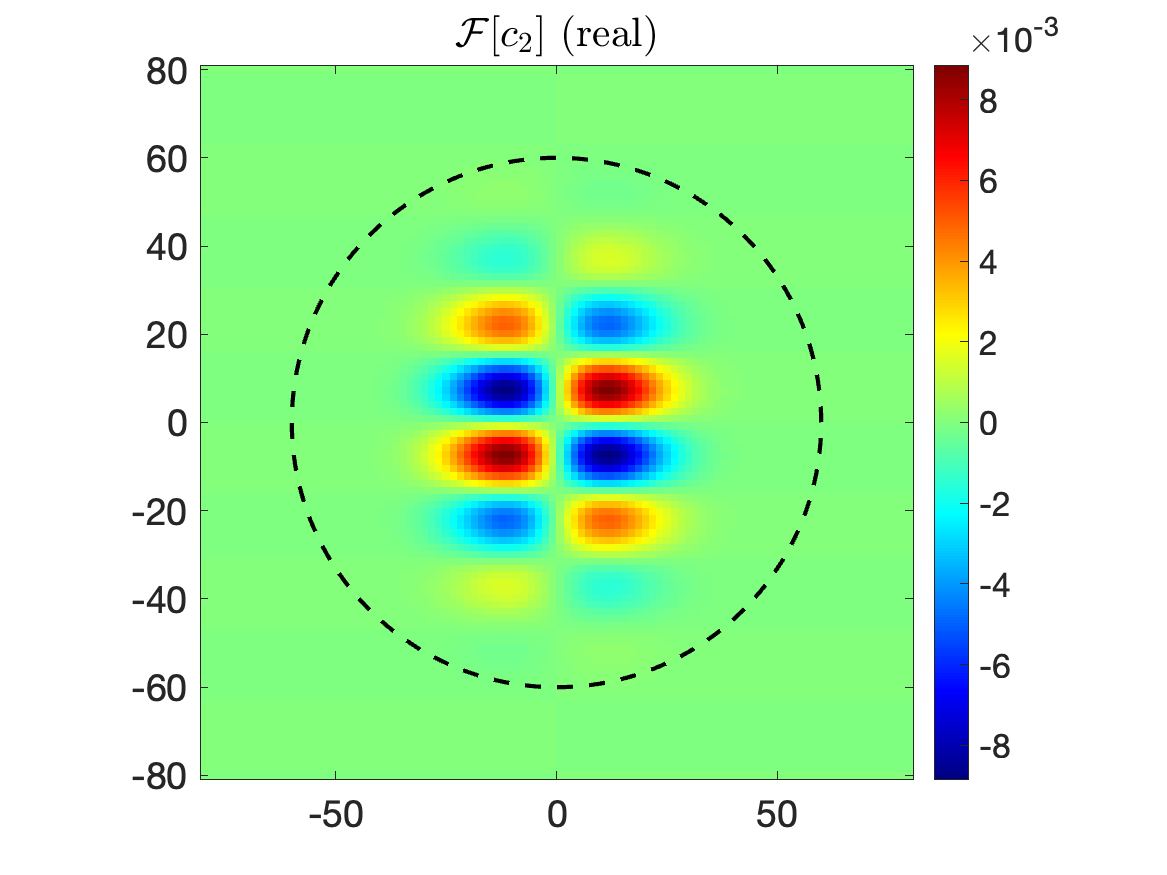} &
\includegraphics[width=0.3\textwidth]{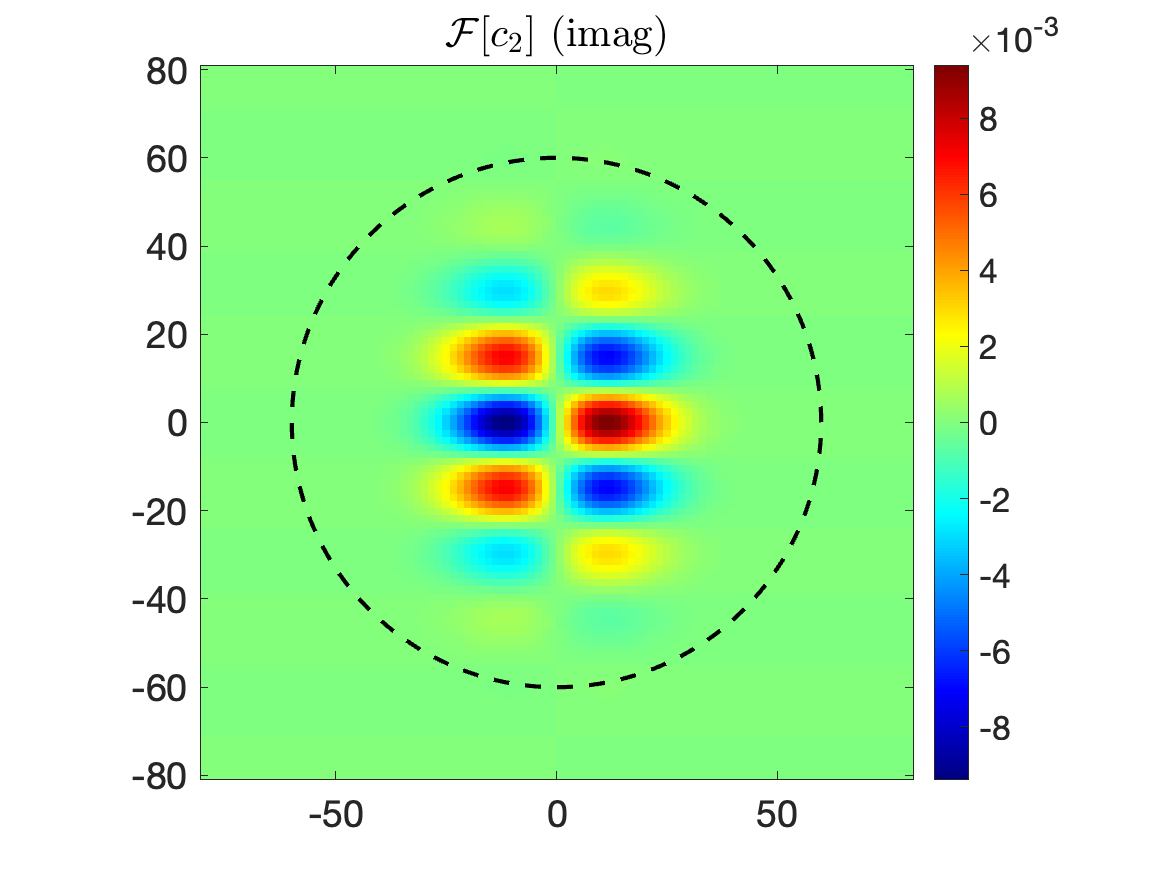} \\
\textbf{(2)} $c_{2}(x)$ &
\textbf{(2-a)} $\mathcal{F}[c_{2}](\xi)$ (real) &
\textbf{(2-b)} $\mathcal{F}[c_{2}](\xi)$ (imaginary) \\
\end{tabular}
\caption{%
\textsf{Left}: The exact coefficients $c_{\ell}(x)$ with the domain boundary $\partial\Omega$ (\textcolor{red}{\bf red} circle), $\ell = 1,2$.
\textsf{Middle \& Right}: The corresponding Fourier modes $\mathcal{F}[c_{\ell}](\xi)$ with the regime boundary $\big\{ \xi : |\xi| = (\ell+1)k \big\}$ (\textcolor{black}{\bf black} dashed circle), $\ell = 1,2$.
\textbf{(1)} The exact coefficient $c_{1}(x)$, \textbf{(1-a)} the real part and \textbf{(1-b)} the imaginary part of Fourier mode $\mathcal{F}[c_{1}](\xi)$.
\textbf{(2)} The exact coefficient $c_{2}(x)$, \textbf{(2-a)} the real part and \textbf{(2-b)} the imaginary part of Fourier mode $\mathcal{F}[c_{2}](\xi)$.
}
\label{figs:poly_m2_exact}
\end{figure}

Let $m = 2$, we first consider the following polynomial in $u$ of degree $2$,
\begin{align*}
P_{2}(x,u) = c_{1}(x) \, u + c_{2}(x) \, u^{2}.
\end{align*}
In this numerical example, the exact coefficients $c_{1}(x)$ and $c_{2}(x)$ are displayed in \textbf{Figure \ref{figs:poly_m2_exact} (1)} and \textbf{(2)} respectively, with the domain boundary $\partial\Omega$ delineated by a red circle.
Recall that, the coefficients $c_{1}$ and $c_{2}$ are compactly supported in $\Omega$.
Furthermore, the real and imaginary parts of the corresponding Fourier modes are also shown in \textbf{Figure \ref{figs:poly_m2_exact}}:
\textbf{Panel (1-a)} and \textbf{(1-b)} for the components of $\mathcal{F}[c_{1}](\xi)$,
\textbf{Panel (2-a)} and \textbf{(2-b)} for the components of $\mathcal{F}[c_{2}](\xi)$.
In these panels, the black dashed circle demarcates the boundary of the Fourier frequency regime $\big\{ \xi \in \mathbb{R}^{2} : |\xi| \leqslant (\ell+1)k \big\}$, where $\ell = 1,2$ and $k = 20$ as specified previously.

\begin{figure}[htbp]
\centering
\begin{tabular}{@{}ccc@{}}
\includegraphics[width=0.33\textwidth]{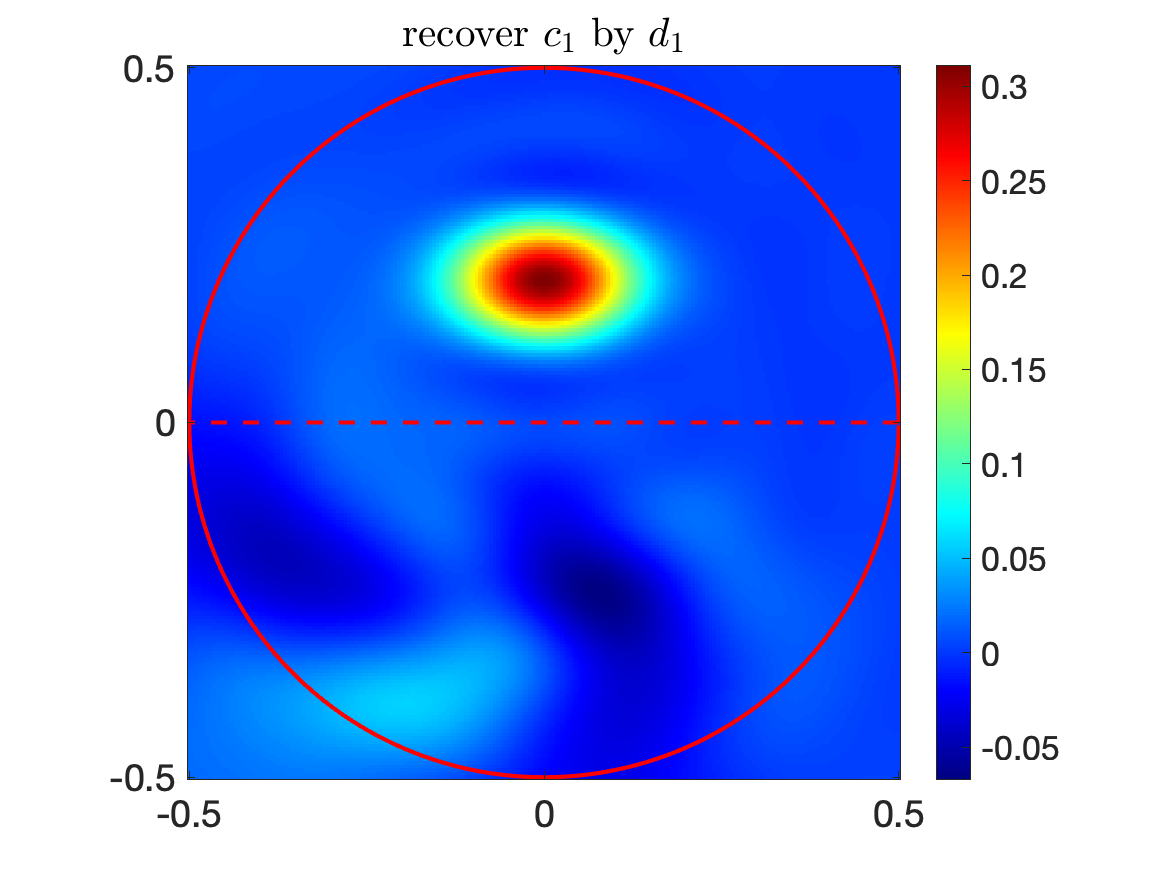} &
\includegraphics[width=0.3\textwidth]{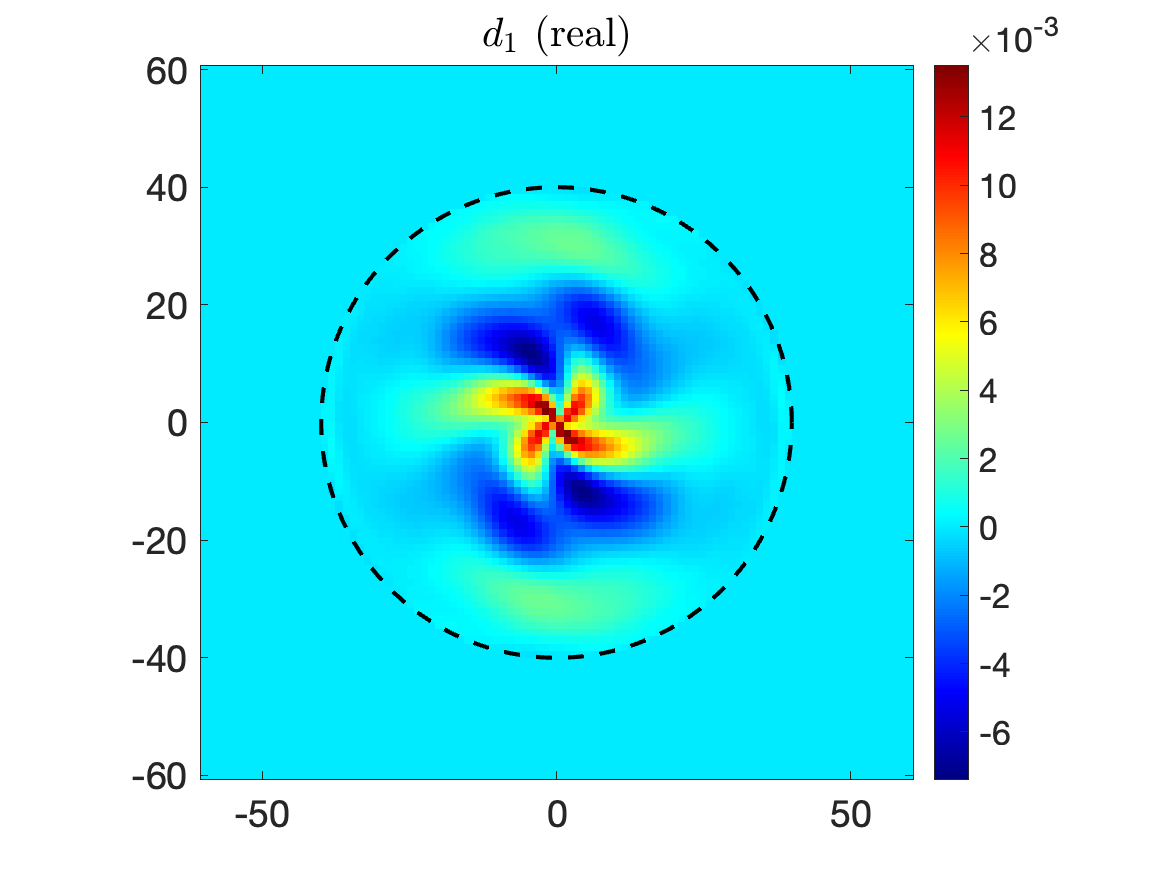} &
\includegraphics[width=0.3\textwidth]{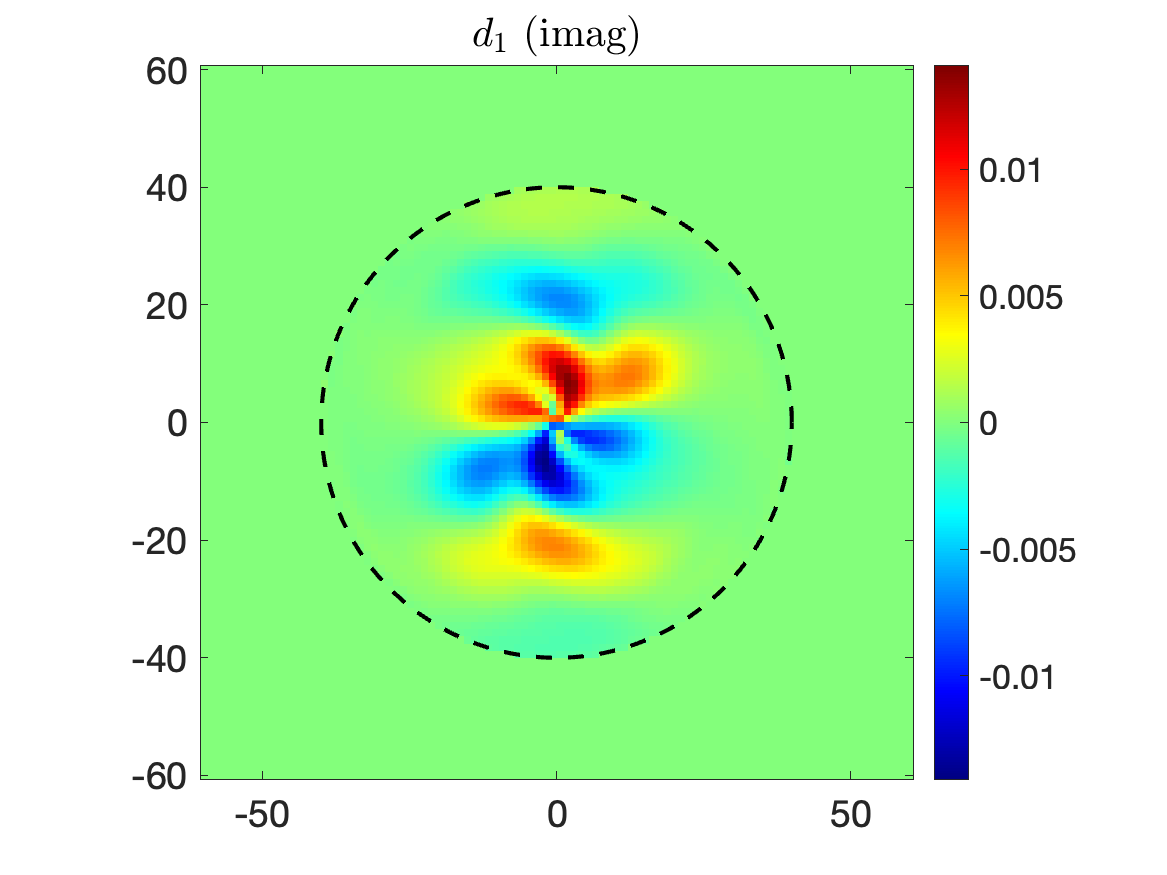} \\
\textbf{(1)} $c^{d}_{1}(x)$ &
\textbf{(1-a)} $d_{1}(\xi)$ (real) &
\textbf{(1-b)} $d_{1}(\xi)$ (imaginary) \\[3ex]
\includegraphics[width=0.33\textwidth]{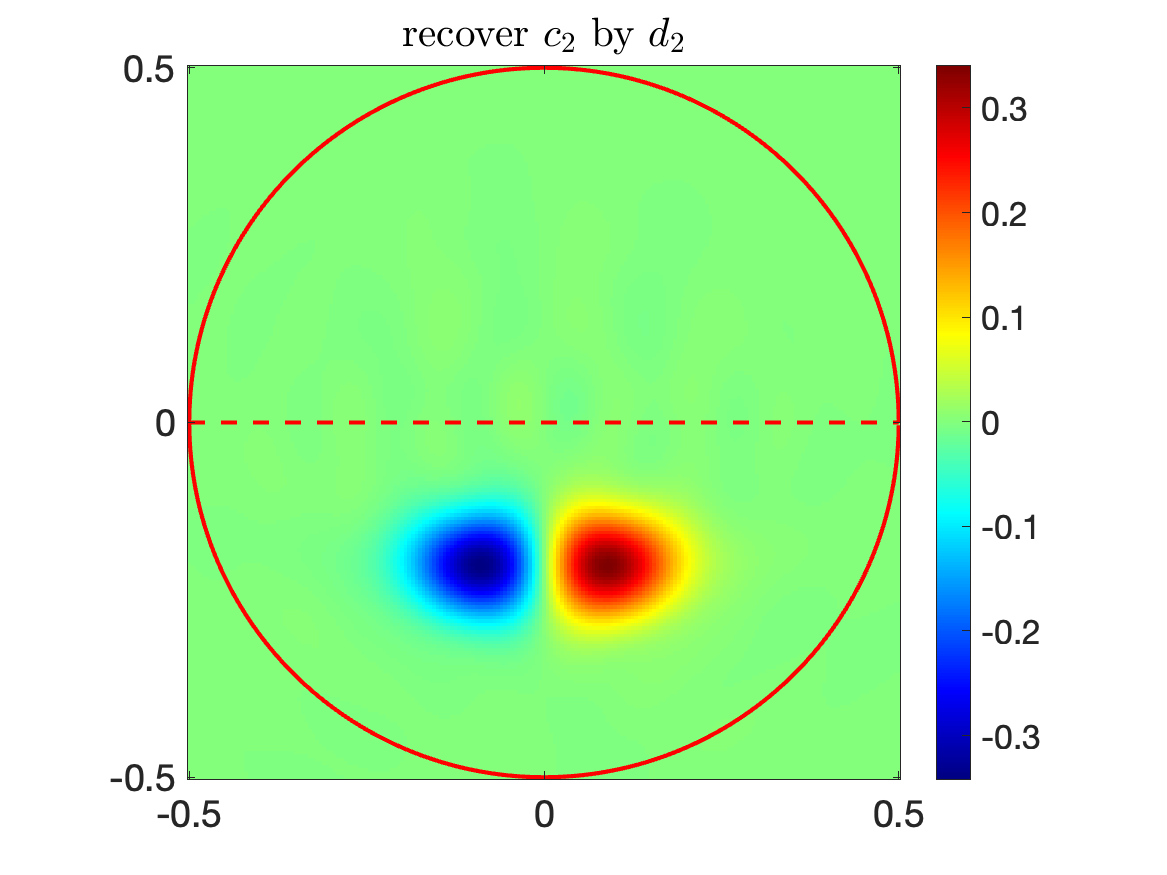} &
\includegraphics[width=0.3\textwidth]{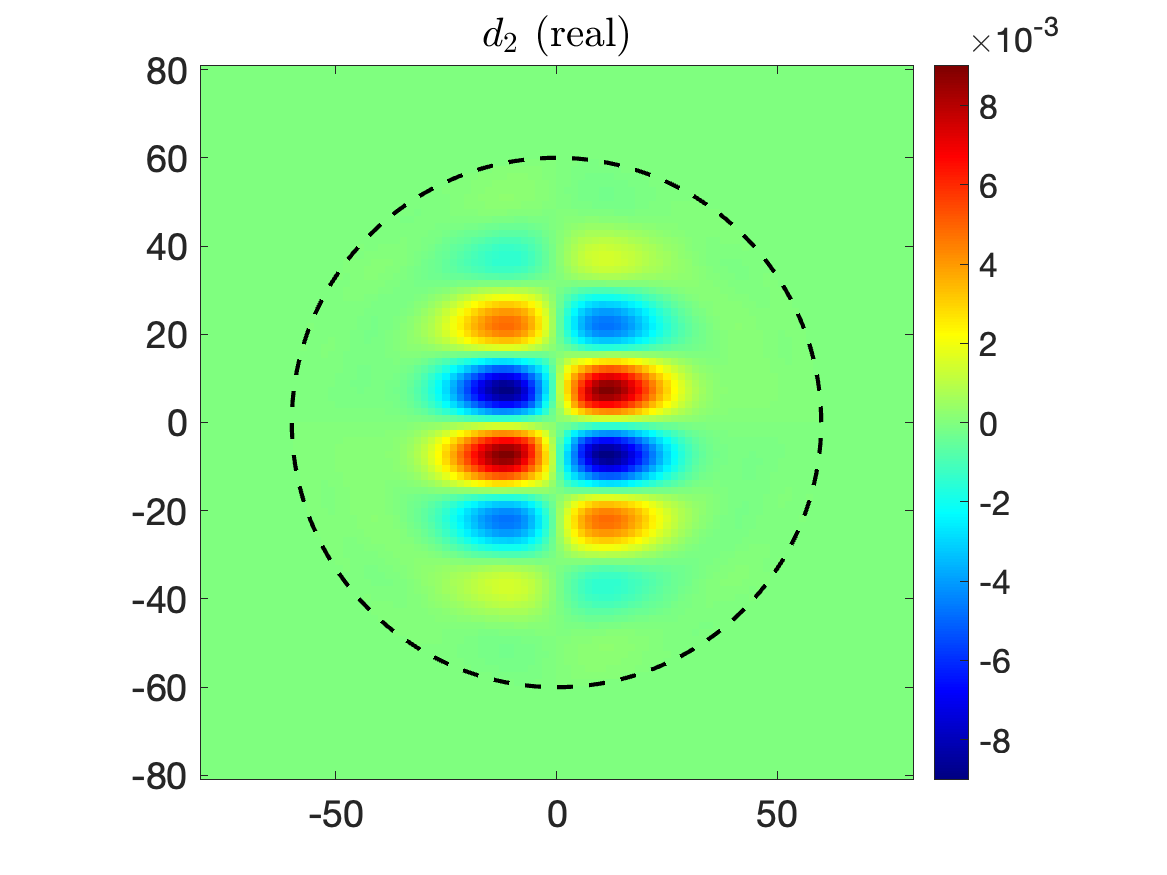} &
\includegraphics[width=0.3\textwidth]{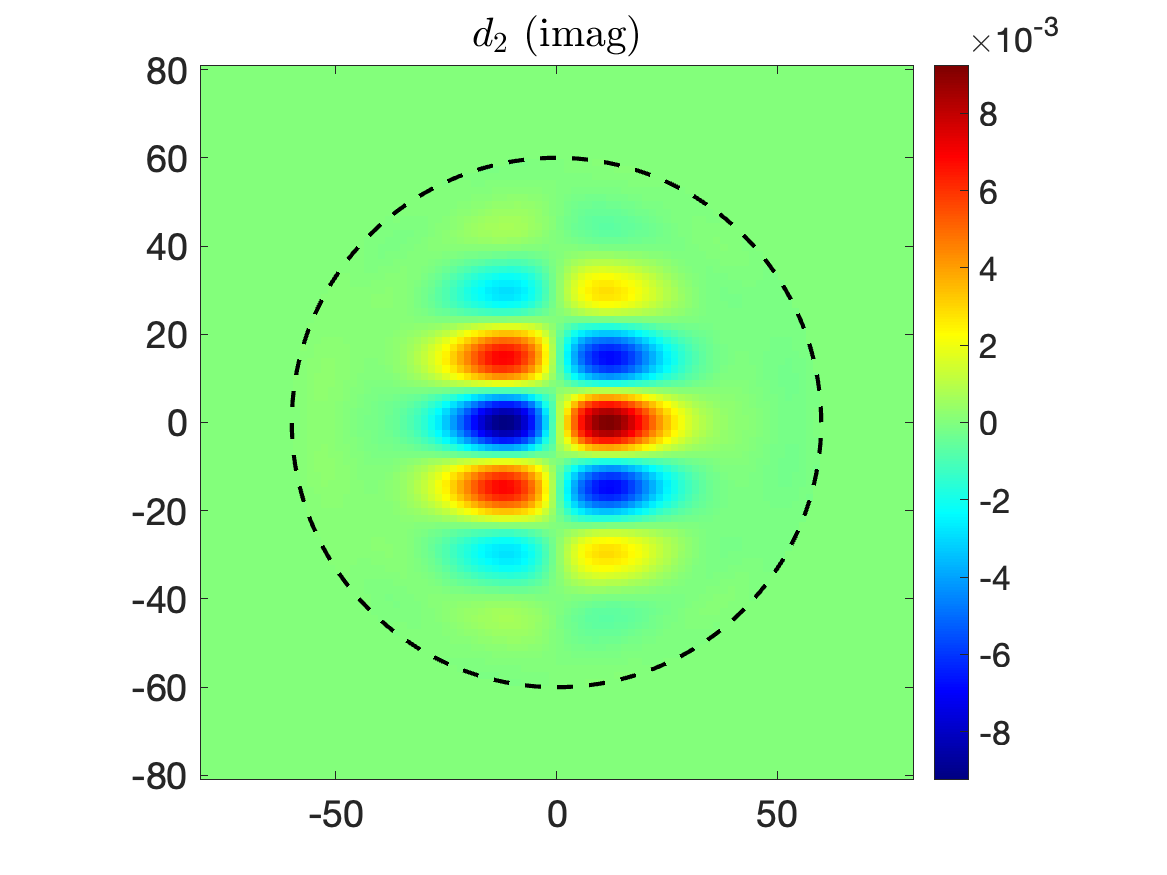} \\
\textbf{(2)} $c^{d}_{2}(x)$ &
\textbf{(2-a)} $d_{2}(\xi)$ (real) &
\textbf{(2-b)} $d_{2}(\xi)$ (imaginary) \\
\end{tabular}
\caption{%
\textsf{Left}: The reconstructed coefficients $c^{d}_{\ell}(x) = \mathcal{F}^{-1}[d_{\ell}](x)$ by using the inverse Fourier transform of the data $d_{\ell}(\xi)$ directly, $\ell = 1,2$.
\textsf{Middle \& Right}: The data $d_{\ell}(\xi)$ computed from measurements on the boundary $\partial\Omega$, $\ell = 1,2$.
\textbf{(1)} The reconstructed coefficient $c^{d}_{1}(x)$, \textbf{(1-a)} the real part and \textbf{(1-b)} the imaginary part of the data $d_{1}(\xi)$.
\textbf{(2)} The reconstructed coefficient $c^{d}_{2}(x)$, \textbf{(2-a)} the real part and \textbf{(2-b)} the imaginary part of the data $d_{2}(\xi)$.
}
\label{figs:poly_m2_d_l}
\end{figure}

Recall that, the data $d_{\ell}(\xi)$ are computed from measurements on the domain boundary $\partial\Omega$; see \eqref{eqn:data_d_l} or \textsf{Line \ref{alg:Fc_l_start}} in \textbf{Algorithm \ref{alg:main}}.
The real and imaginary parts of the data $d_{1}(\xi)$ and $d_{2}(\xi)$ are displayed in \textbf{Figure \ref{figs:poly_m2_d_l}}:
\textbf{Panel (1-a)} and \textbf{(1-b)} for the components of $d_{1}(\xi)$,
\textbf{Panel (2-a)} and \textbf{(2-b)} for the components of $d_{2}(\xi)$.
A direct comparison between \textbf{Figure \ref{figs:poly_m2_exact} (2-a)--(2-b)} and \textbf{Figure \ref{figs:poly_m2_d_l} (2-a)--(2-b)} verifies the equality
\begin{align*}
\mathcal{F}[c_{2}] = d_{2},
\end{align*}
in agreement with \textbf{Lemma \ref{lmm:formulas}}: \textbf{Case $c_{m}$} for $m = 2$.
In contrast, the Fourier mode
\begin{align*}
\mathcal{F}[c_{1}] \neq d_{1},
\end{align*}
as evidenced by comparing between \textbf{Figure \ref{figs:poly_m2_exact} (1-a)--(1-b)} (the exact Fourier mode $\mathcal{F}[c_{1}]$) and \textbf{Figure \ref{figs:poly_m2_d_l} (1-a)--(1-b)} (the computed data $d_{1}$).

For each $\ell = 1,2$, we denote
\begin{align*}
c^{d}_{\ell}(x) := \mathcal{F}^{-1}[d_{\ell}](x)
\end{align*}
as the coefficients reconstructed directly from the inverse Fourier transform of the data $d_{\ell}$.
Thus, these numerical results demonstrate that:
The reconstructed coefficient $c^{d}_{2}$ converges accurately to the exact coefficient $c_{2}$; see \textbf{Figure \ref{figs:poly_m2_exact} (2)} and \textbf{Figure \ref{figs:poly_m2_d_l} (2)}.
The reconstructed coefficient $c^{d}_{1}$ exhibits distortion due to interference from the higher-order coefficient $c_{2}$; see \textbf{Figure \ref{figs:poly_m2_d_l} (1)}.

\begin{figure}[htbp]
\centering
\begin{tabular}{@{}ccc@{}}
\includegraphics[width=0.33\textwidth]{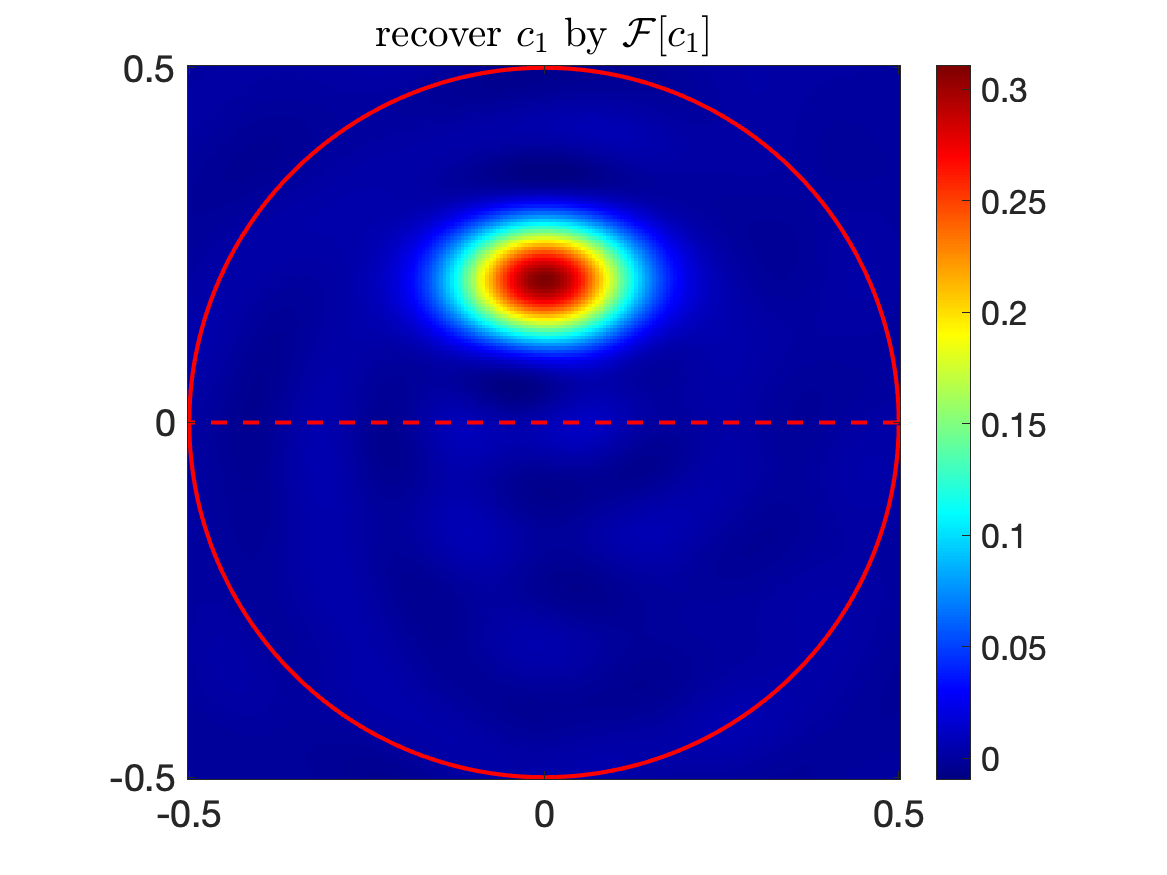} &
\includegraphics[width=0.3\textwidth]{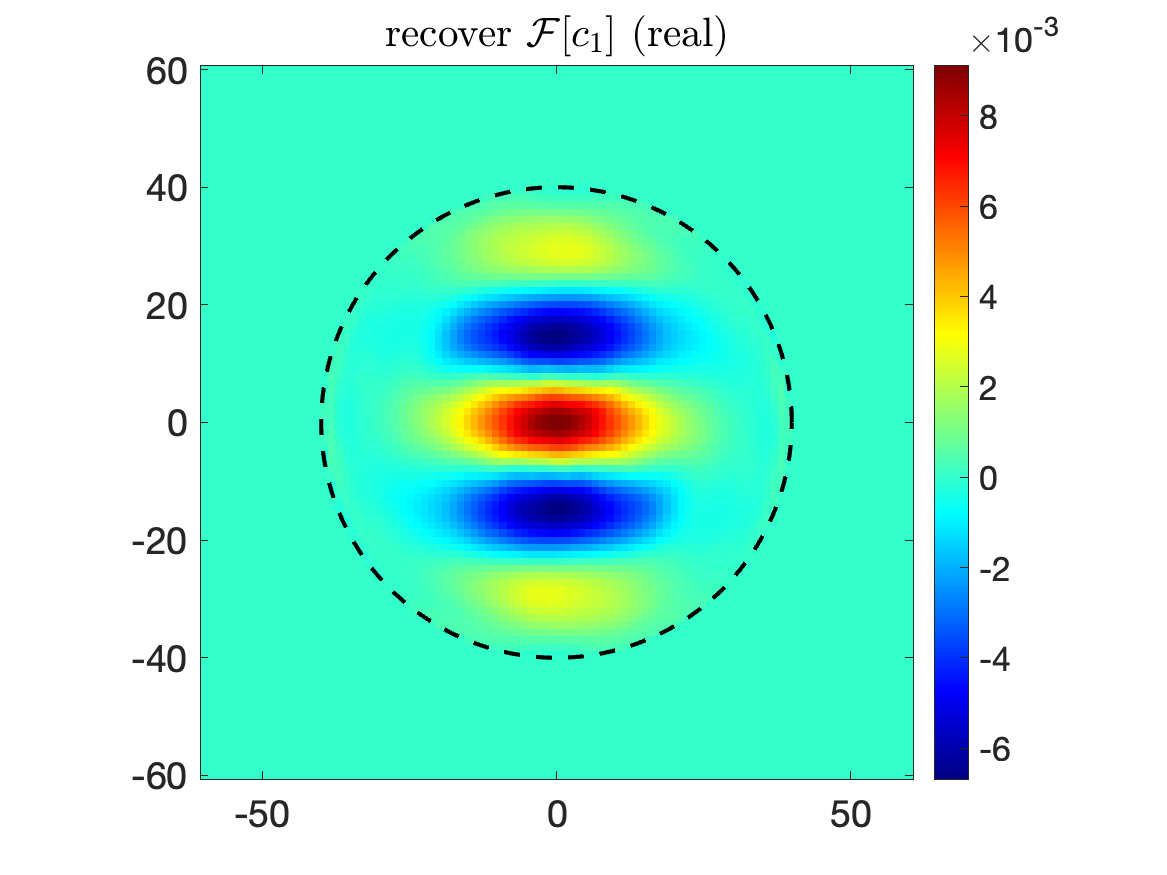} &
\includegraphics[width=0.3\textwidth]{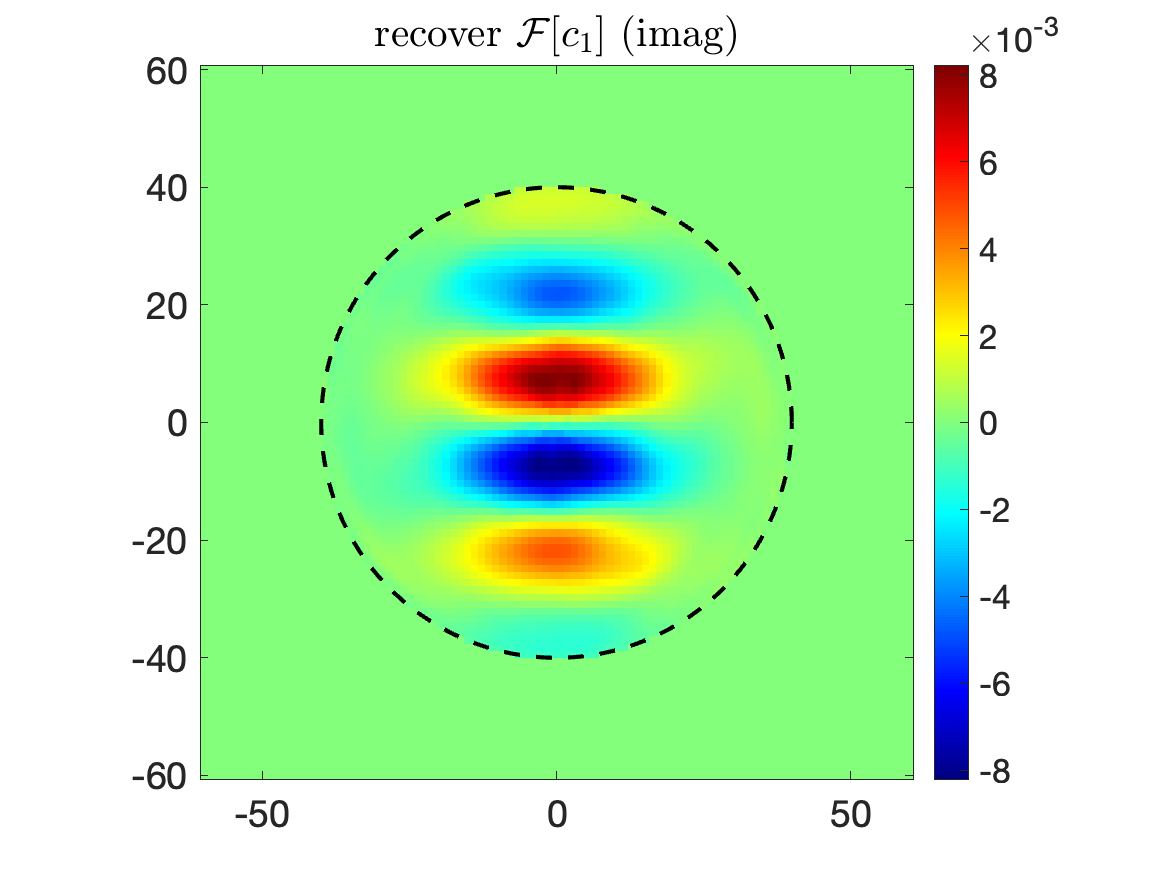} \\
\textbf{(1)} $\widetilde{c}_{1}(x)$ &
\textbf{(1-a)} $\widetilde{d}_{1}(\xi)$ (real) &
\textbf{(1-b)} $\widetilde{d}_{1}(\xi)$ (imaginary) \\[3ex]
\end{tabular}
\caption{%
\textsf{Left}: The recovered coefficients $\widetilde{c}_{\ell}(x) = \mathcal{F}^{-1}[\widetilde{d}_{\ell}](x)$ by using the inverse transform of the modified data $\widetilde{d}_{1}(\xi)$, $\ell = 1$.
\textsf{Middle \& Right}: The modified data $\widetilde{d}_{\ell}(\xi)$ updated by using previously modified data $\big\{ \widetilde{d}_{\ell+a}(\xi) : a = 1,2,\dots,m-\ell \big\}$ and the data $d_{\ell}(\xi)$, $\ell = 1$.
\textbf{(1)} The recovered coefficient $\widetilde{c}_{1}(x)$, \textbf{(1-a)} the real part and \textbf{(1-b)} the imaginary part of the modified data $\widetilde{d}_{1}(\xi)$.
}
\label{figs:poly_m2_Fc_l}
\end{figure}

Furthermore, by \textbf{Lemma \ref{lmm:formulas}}, we then denote
\begin{align*}
\widetilde{d}_{\ell}(\xi)
:= \left\{~
\begin{aligned}
& \\[-1ex]
& d_{m}(\xi) & &~ \text{if\ } \ell = m, \\[1ex]
& d_{\ell}(\xi) - {\textstyle\sum\limits_{a=1}^{m-\ell}} \frac{1}{\,\ell\,!\,} {\textstyle\sum\limits_{\bm{\alpha} \in \mathcal{A}^{\ell}_{a}}} {\textstyle\binom{\ell+a}{\bm{1}+\bm{\alpha}}} \, \widetilde{d}_{\ell+a}( \xi + {\textstyle\sum\limits_{j \in \mathcal{U}_{\ell}}} \alpha_{j} \zeta_{\ell,j} ) & &~ \text{if\ } 0 < \ell < m.
\end{aligned}
\right.
\end{align*}
as ``the modified data'', and denote
\begin{align*}
\widetilde{c}_{\ell}(x) := \mathcal{F}^{-1}[\widetilde{d}_{\ell}](x)
\end{align*}
as the coefficients recovered from the inverse Fourier transform of the modified data $\widetilde{d}_{\ell}$.
Indeed, for $m = 2$, the modified data $\widetilde{d}_{\ell}(\xi)$ ($\ell = 1,2$) is updated by
\begin{align*}
\widetilde{d}_{2}(\xi) = d_{2}(\xi),
\qquad
\widetilde{d}_{1}(\xi) = d_{1}(\xi) - \widetilde{d}_{2}(\xi+\zeta_{1,1}),
\end{align*}
respectively.
Since $\widetilde{d}_{2} = d_{2}$ in this example, we only show the modified data $\widetilde{d}_{1}(\xi)$ in \textbf{Figure \ref{figs:poly_m2_Fc_l} (1-a)--(1-b)}, which verifies the equality
\begin{align*}
\mathcal{F}[c_{1}] = \widetilde{d}_{1},
\end{align*}
in agreement with \textbf{Lemma \ref{lmm:formulas}}: \textbf{Case $c_{\ell}$ with $0 < \ell < m$}.
Thus, the recovered coefficient $\widetilde{c}_{1}(x)$ converges to the exact coefficient $c_{1}$; see \textbf{Figure \ref{figs:poly_m2_Fc_l} (1)}.

\subsubsection{The example for the degree $m = 3$}

\begin{figure}[htbp]
\centering
\begin{tabular}{@{}ccc@{}}
\includegraphics[width=0.33\textwidth]{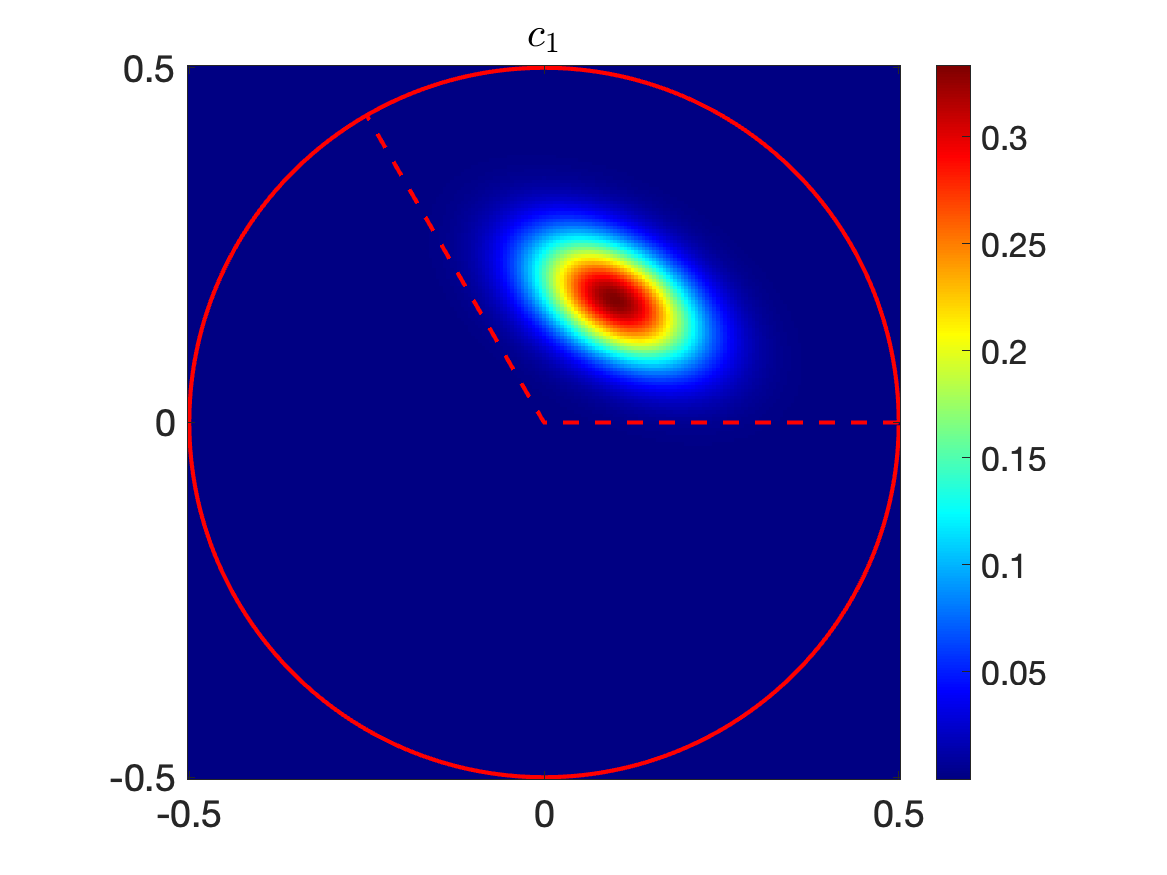} &
\includegraphics[width=0.3\textwidth]{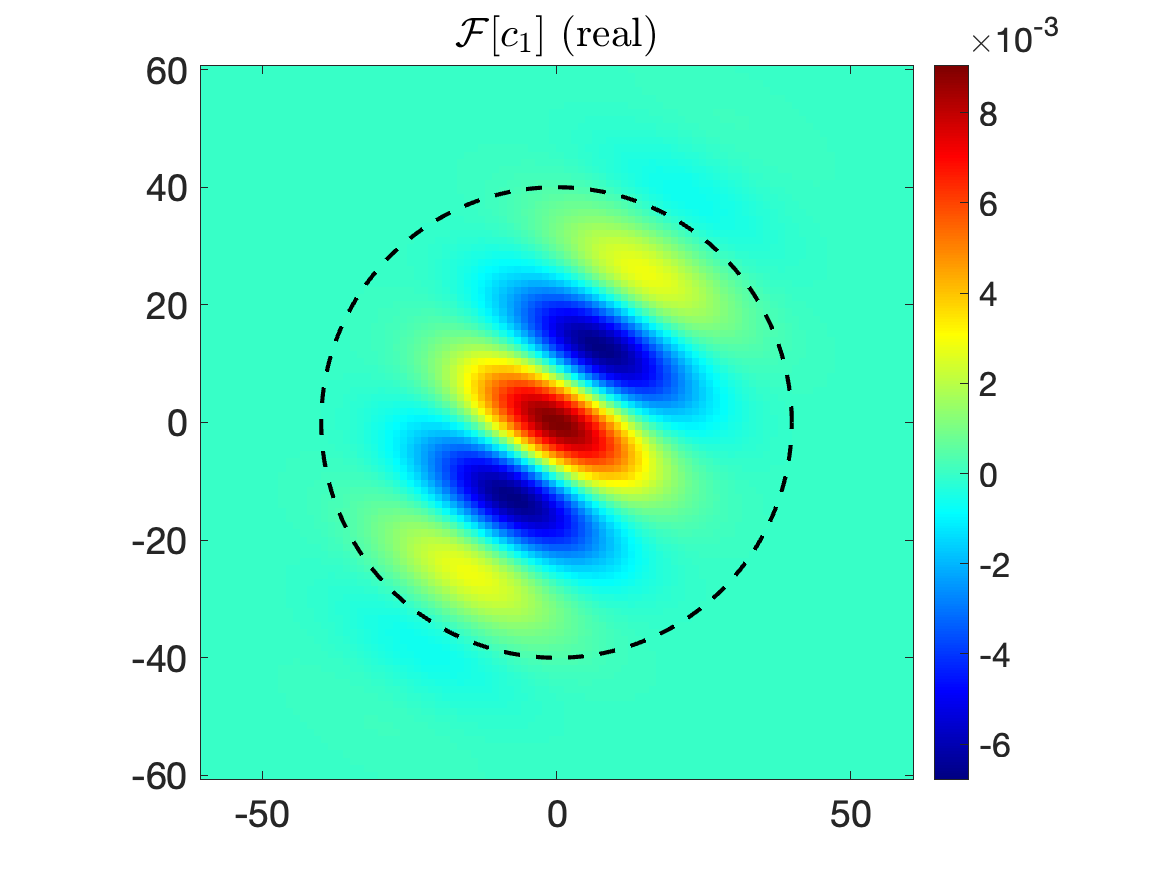} &
\includegraphics[width=0.3\textwidth]{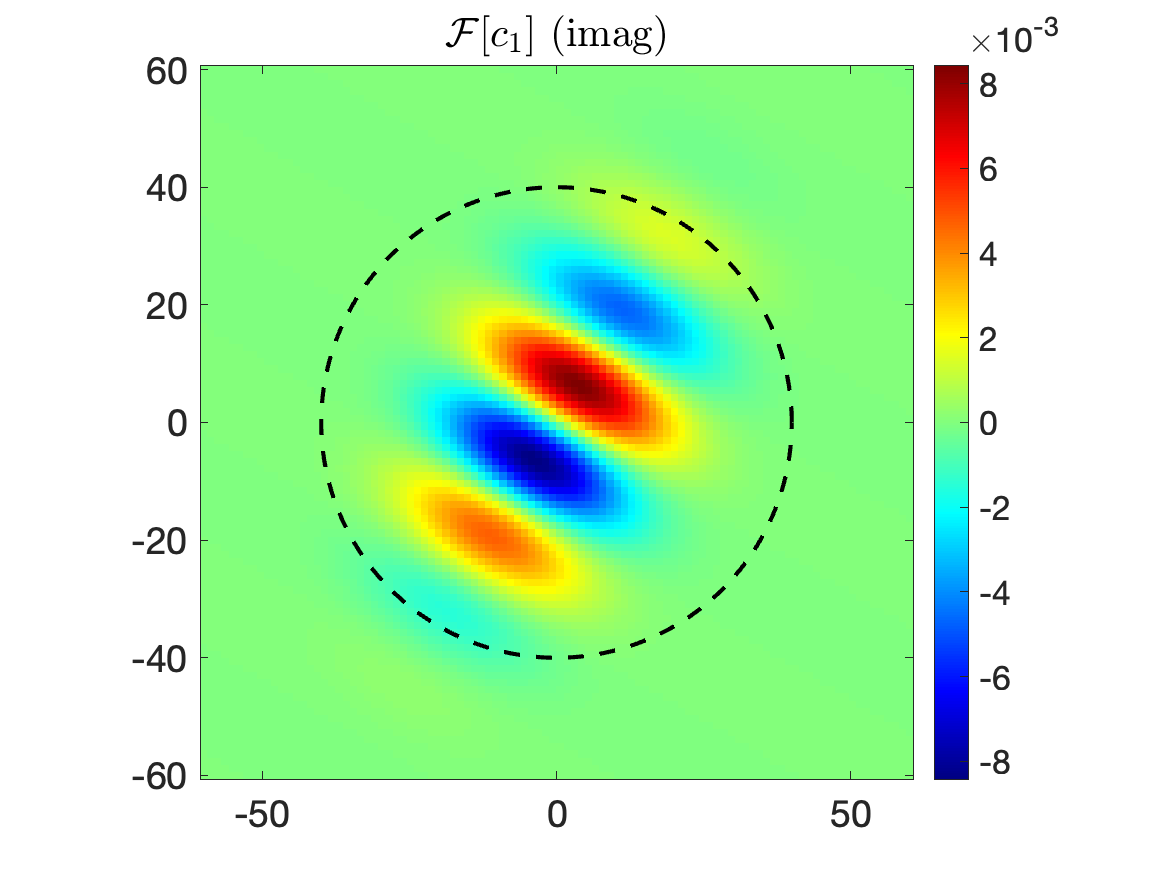} \\
\textbf{(1)} $c_{1}(x)$ &
\textbf{(1-a)} $\mathcal{F}[c_{1}](\xi)$ (real) &
\textbf{(1-b)} $\mathcal{F}[c_{1}](\xi)$ (imaginary) \\[3ex]
\includegraphics[width=0.33\textwidth]{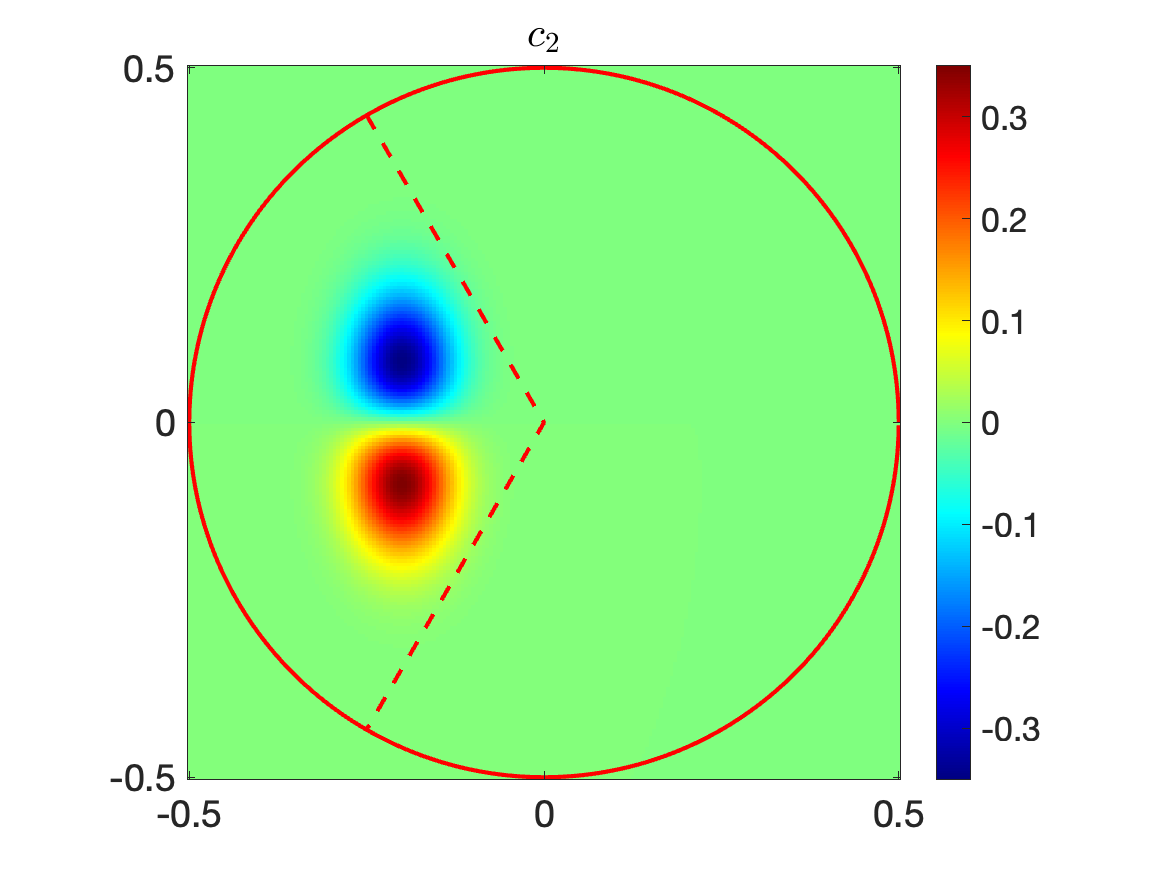} &
\includegraphics[width=0.3\textwidth]{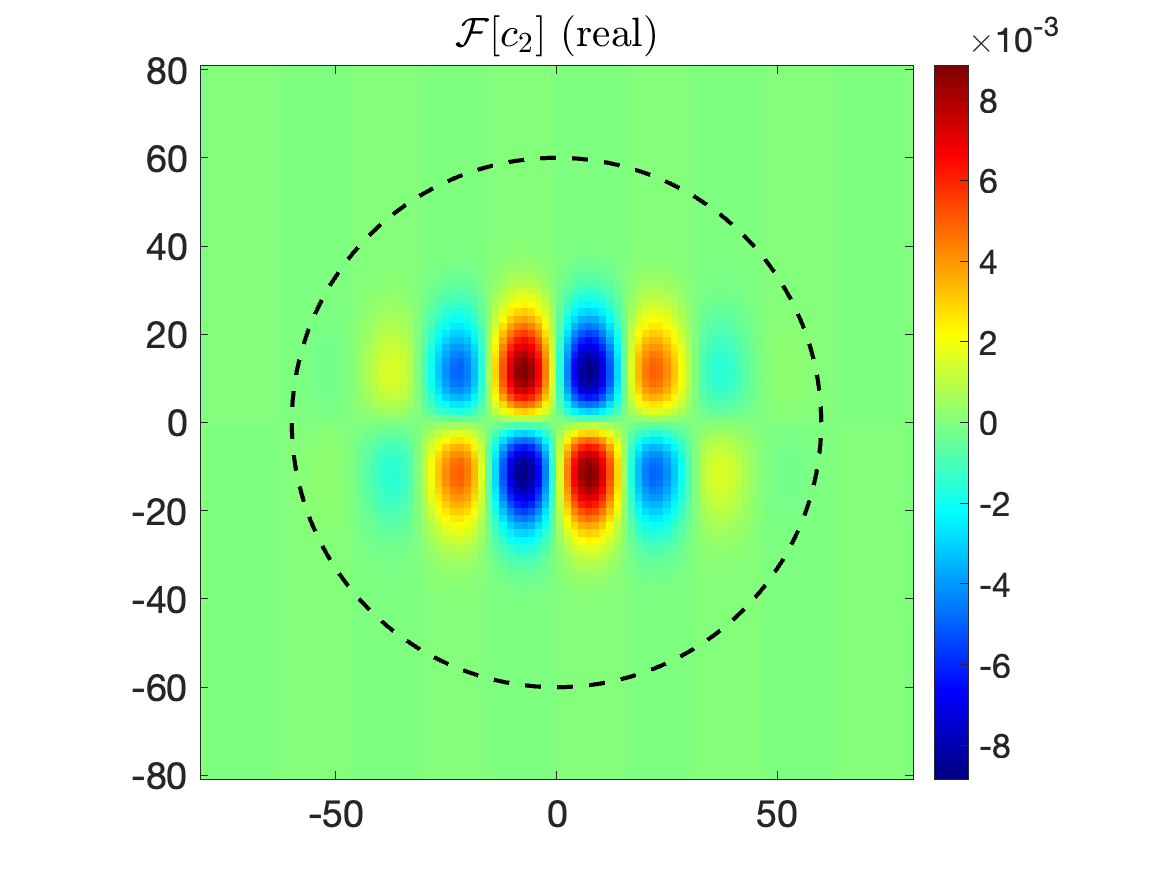} &
\includegraphics[width=0.3\textwidth]{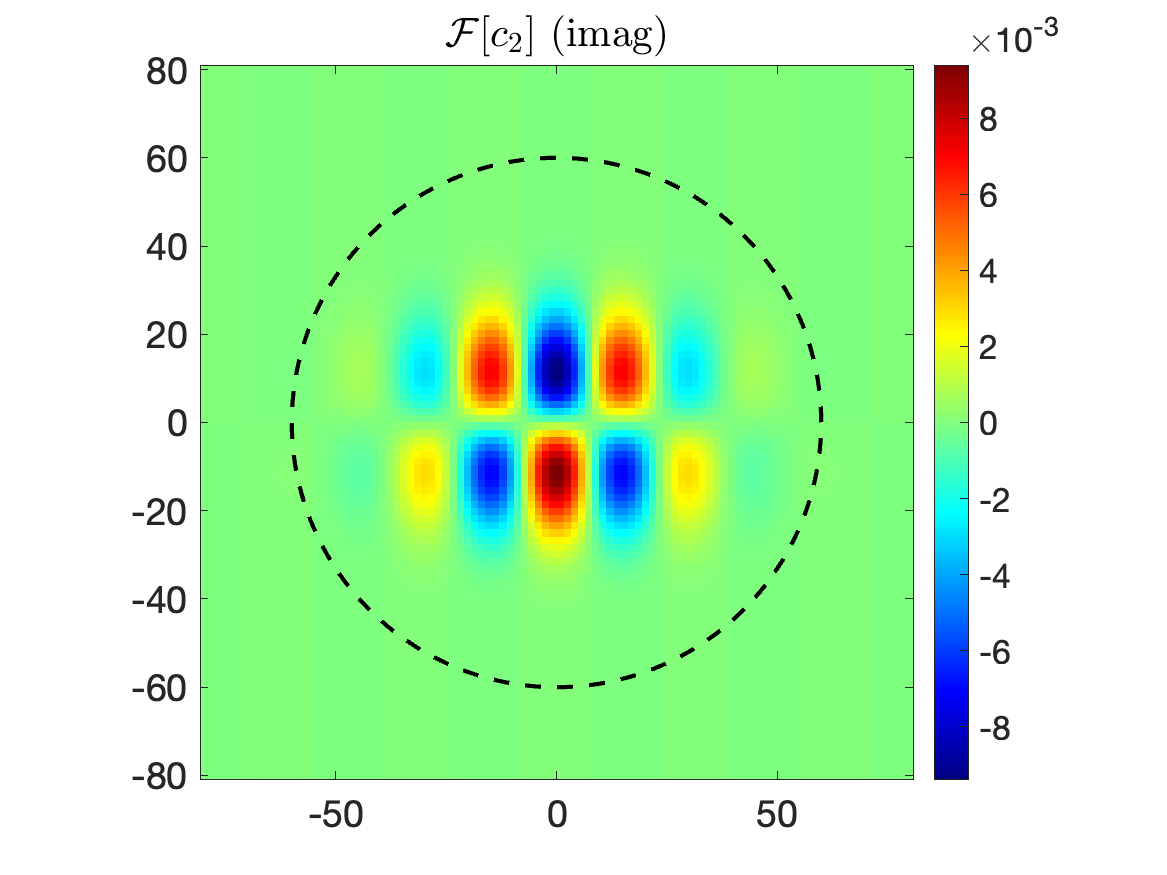} \\
\textbf{(2)} $c_{2}(x)$ &
\textbf{(2-a)} $\mathcal{F}[c_{2}](\xi)$ (real) &
\textbf{(2-b)} $\mathcal{F}[c_{2}](\xi)$ (imaginary) \\[3ex]
\includegraphics[width=0.33\textwidth]{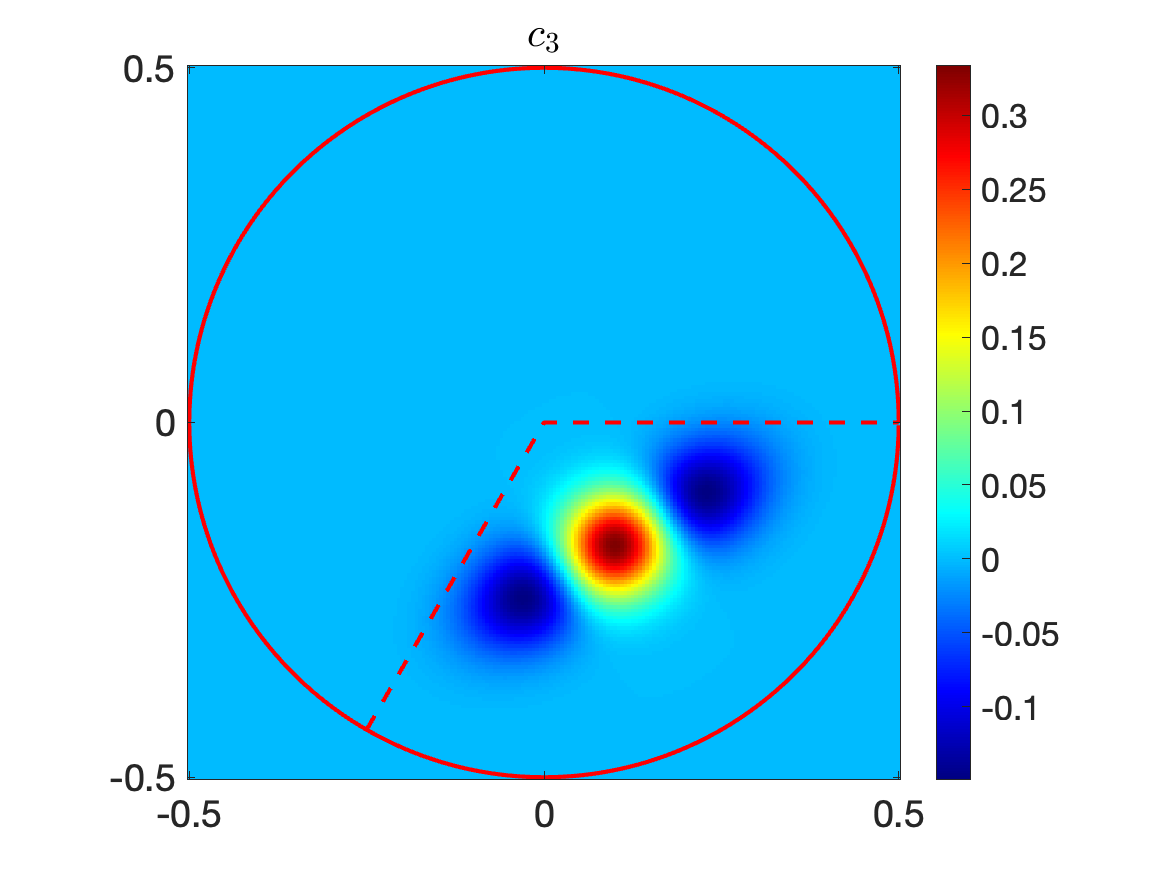} &
\includegraphics[width=0.3\textwidth]{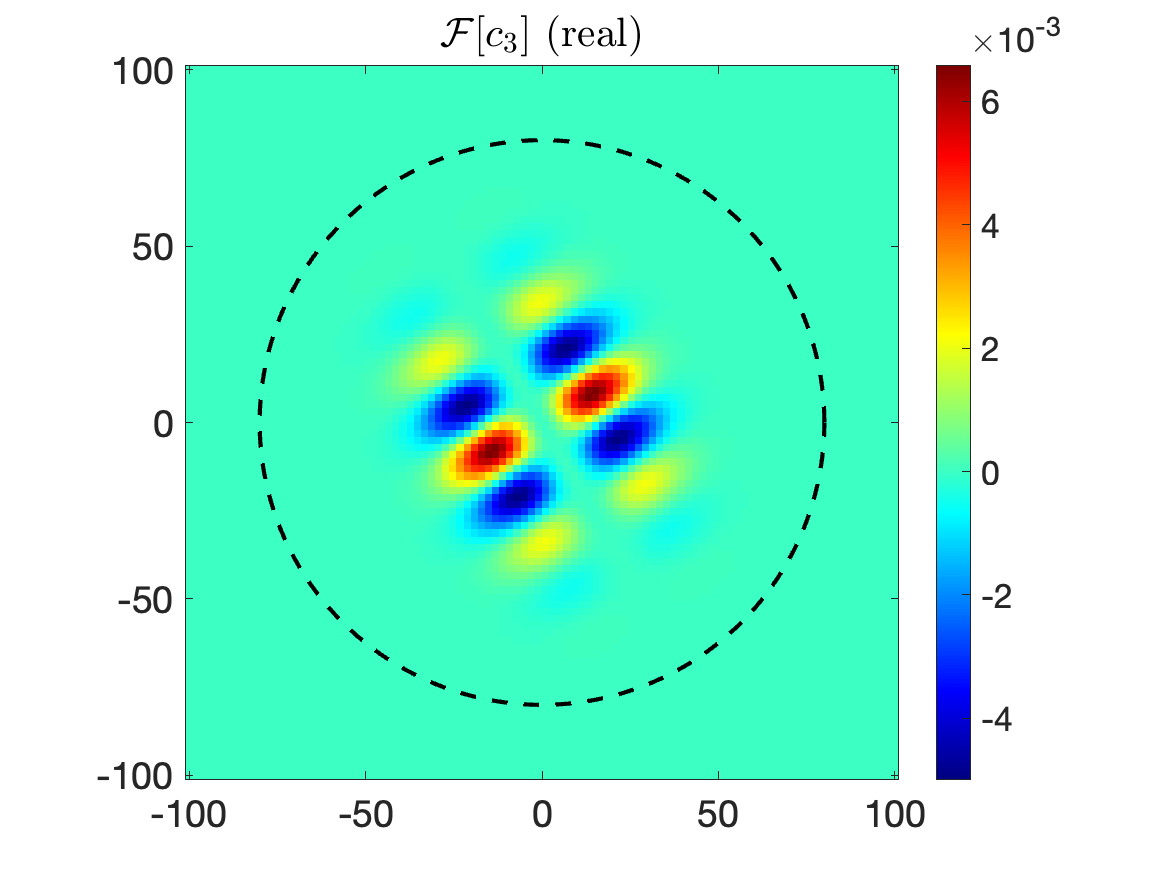} &
\includegraphics[width=0.3\textwidth]{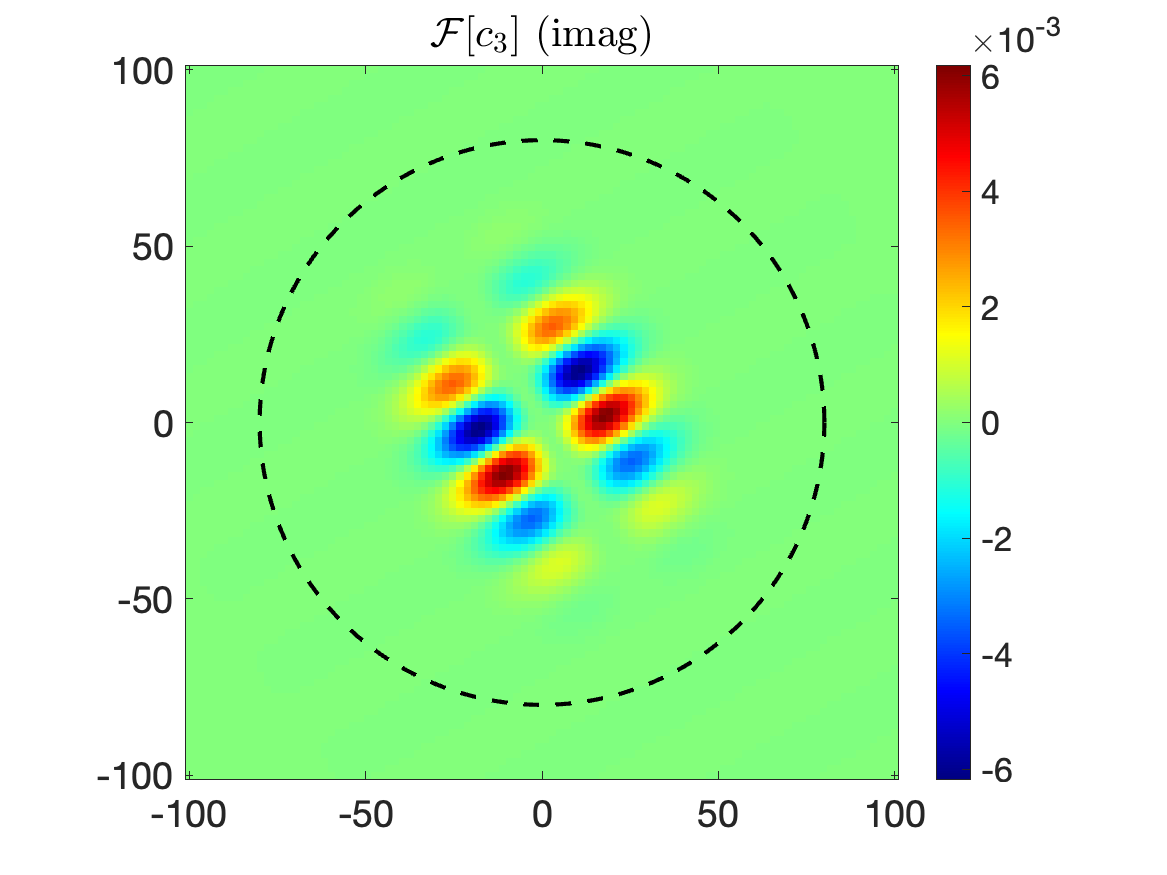} \\
\textbf{(3)} $c_{3}(x)$ &
\textbf{(3-a)} $\mathcal{F}[c_{3}](\xi)$ (real) &
\textbf{(3-b)} $\mathcal{F}[c_{3}](\xi)$ (imaginary) \\
\end{tabular}
\caption{%
\textsf{Left}: The exact coefficients $c_{\ell}(x)$ with the domain boundary $\partial\Omega$ (\textcolor{red}{\bf red} circle), $\ell = 1,2,3$.
\textsf{Middle \& Right}: The corresponding Fourier modes $\mathcal{F}[c_{\ell}](\xi)$ with the regime boundary $\big\{ \xi : |\xi| = (\ell+1)k \big\}$ (\textcolor{black}{\bf black} dashed circle), $\ell = 1,2,3$.
\textbf{(1)} The exact coefficient $c_{1}(x)$, \textbf{(1-a)} the real part and \textbf{(1-b)} the imaginary part of Fourier mode $\mathcal{F}[c_{1}](\xi)$.
\textbf{(2)} The exact coefficient $c_{2}(x)$, \textbf{(2-a)} the real part and \textbf{(2-b)} the imaginary part of Fourier mode $\mathcal{F}[c_{2}](\xi)$.
\textbf{(3)} The exact coefficient $c_{3}(x)$, \textbf{(3-a)} the real part and \textbf{(3-b)} the imaginary part of Fourier mode $\mathcal{F}[c_{3}](\xi)$.
}
\label{figs:poly_m3_exact}
\end{figure}

Let $m = 3$, we then consider the following polynomial in $u$ of degree $3$,
\begin{align*}
P_{3}(x,u) = c_{1}(x) \, u + c_{2}(x) \, u^{2} + c_{3}(x) \, u^{3}.
\end{align*}
In this numerical example, the exact coefficients $c_{1}(x)$, $c_{2}(x)$, and $c_{3}(x)$ are displayed in \textbf{Figure \ref{figs:poly_m3_exact} (1)}, \textbf{(2)}, and \textbf{(3)} respectively, with the domain boundary $\partial\Omega$ delineated by a red circle.
The real and imaginary parts of the corresponding Fourier modes are shown in \textbf{Figure \ref{figs:poly_m3_exact}}:
\textbf{Panel (1-a)} and \textbf{(1-b)} for the components of $\mathcal{F}[c_{1}](\xi)$,
\textbf{Panel (2-a)} and \textbf{(2-b)} for the components of $\mathcal{F}[c_{2}](\xi)$,
\textbf{Panel (3-a)} and \textbf{(3-b)} for the components of $\mathcal{F}[c_{3}](\xi)$.
In these panels, the black dashed circle demarcates the boundary of the Fourier frequency regime $\big\{ \xi \in \mathbb{R}^{2} : |\xi| \leqslant (\ell+1)k \big\}$ with $\ell = 1,2,3$.

\begin{figure}[htbp]
\centering
\begin{tabular}{@{}ccc@{}}
\includegraphics[width=0.33\textwidth]{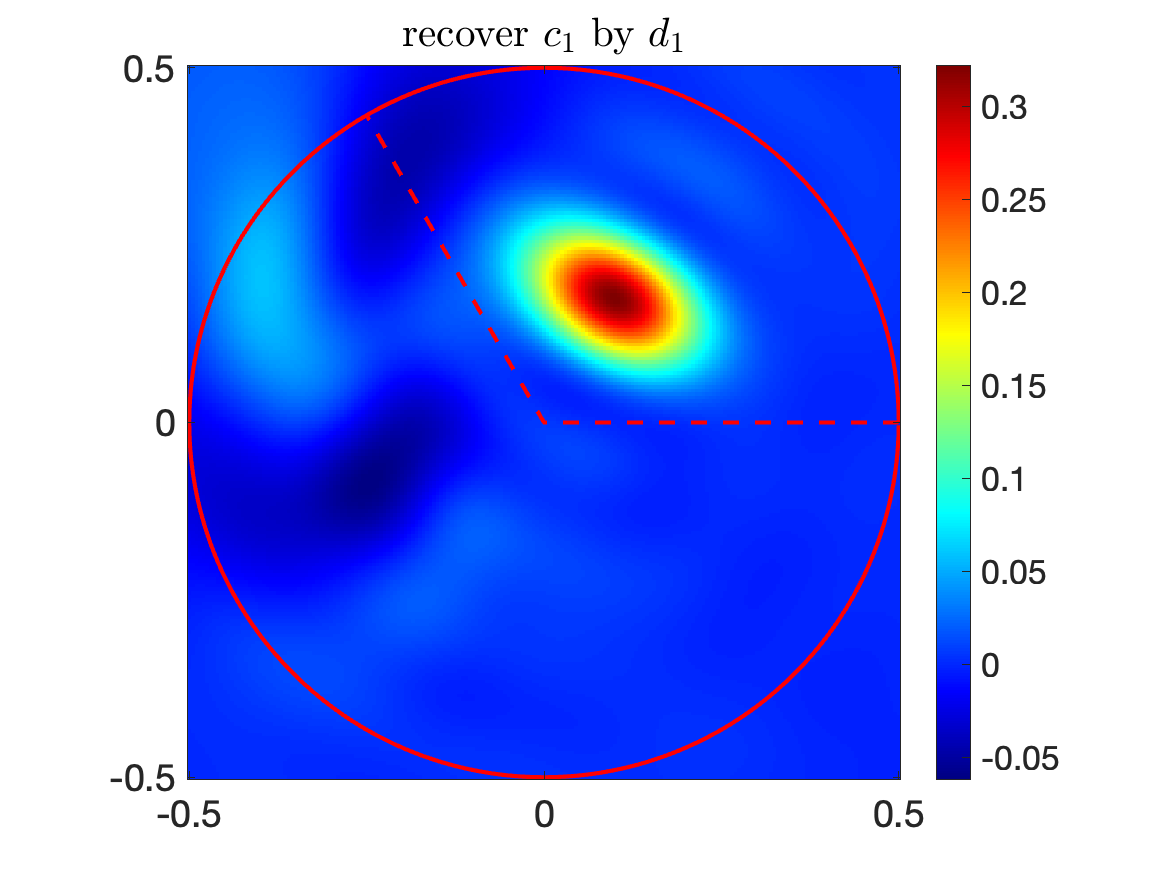} &
\includegraphics[width=0.3\textwidth]{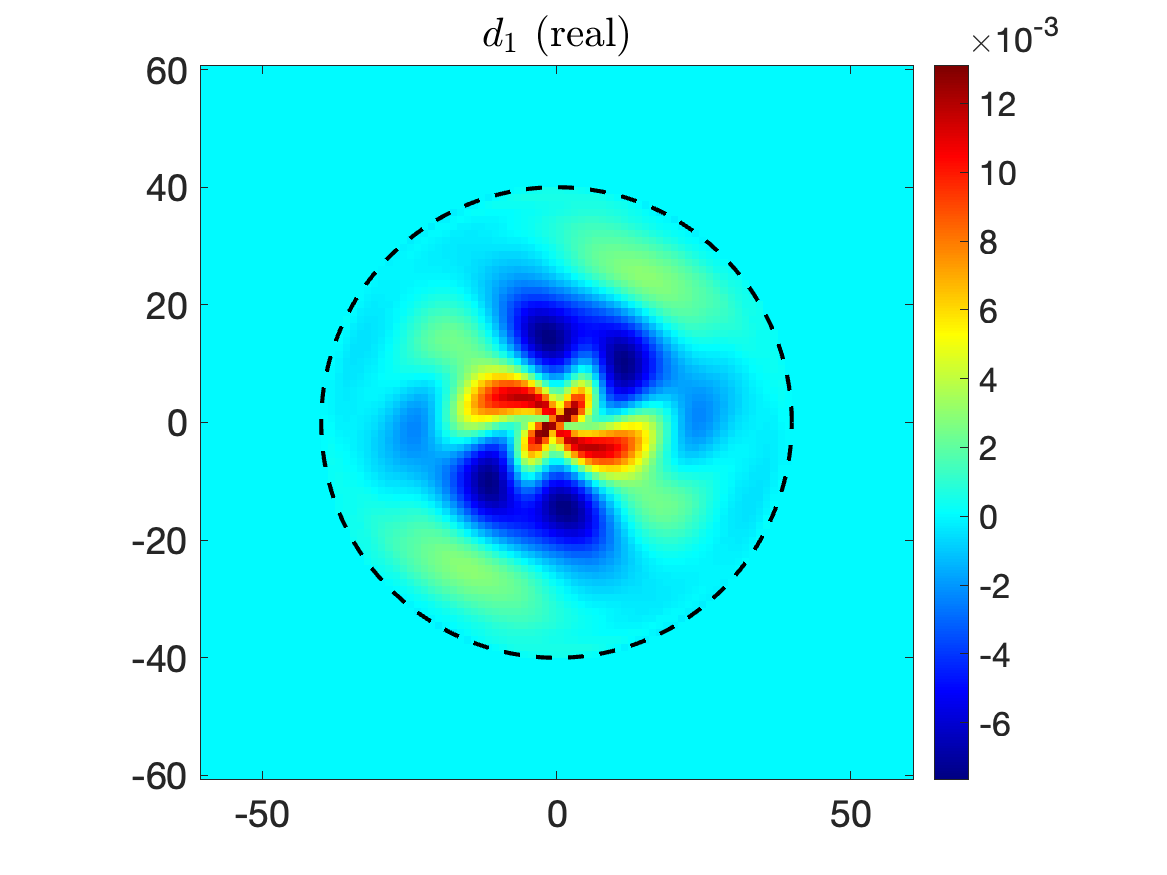} &
\includegraphics[width=0.3\textwidth]{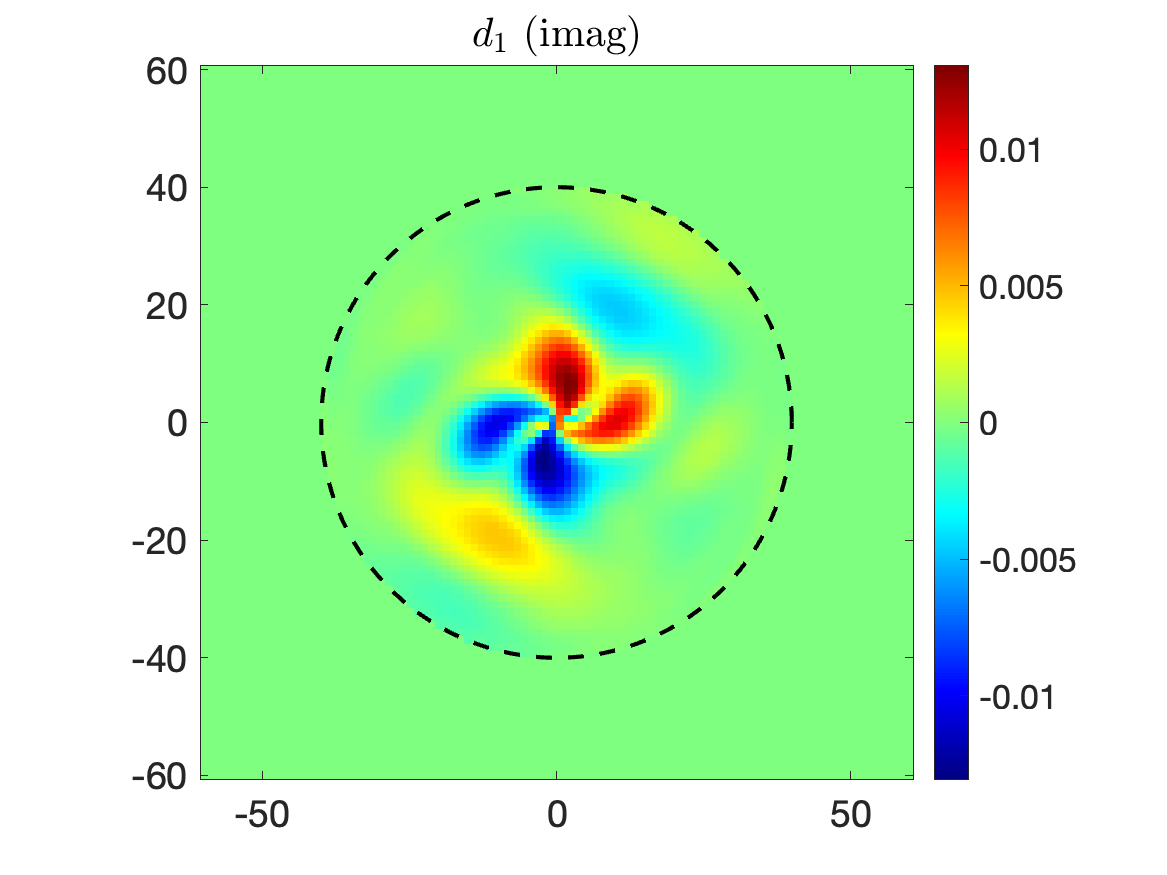} \\
\textbf{(1)} $c^{d}_{1}(x)$ &
\textbf{(1-a)} $d_{1}(\xi)$ (real) &
\textbf{(1-b)} $d_{1}(\xi)$ (imaginary) \\[3ex]
\includegraphics[width=0.33\textwidth]{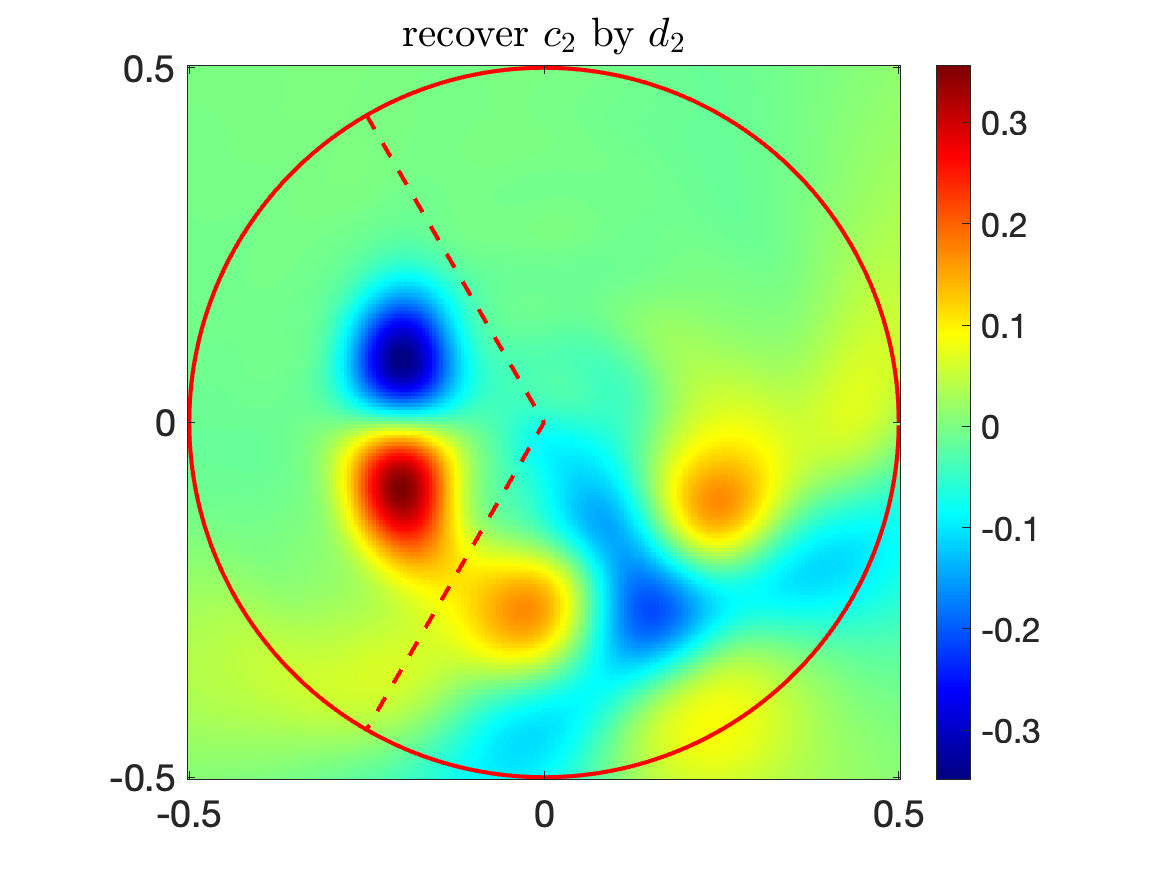} &
\includegraphics[width=0.3\textwidth]{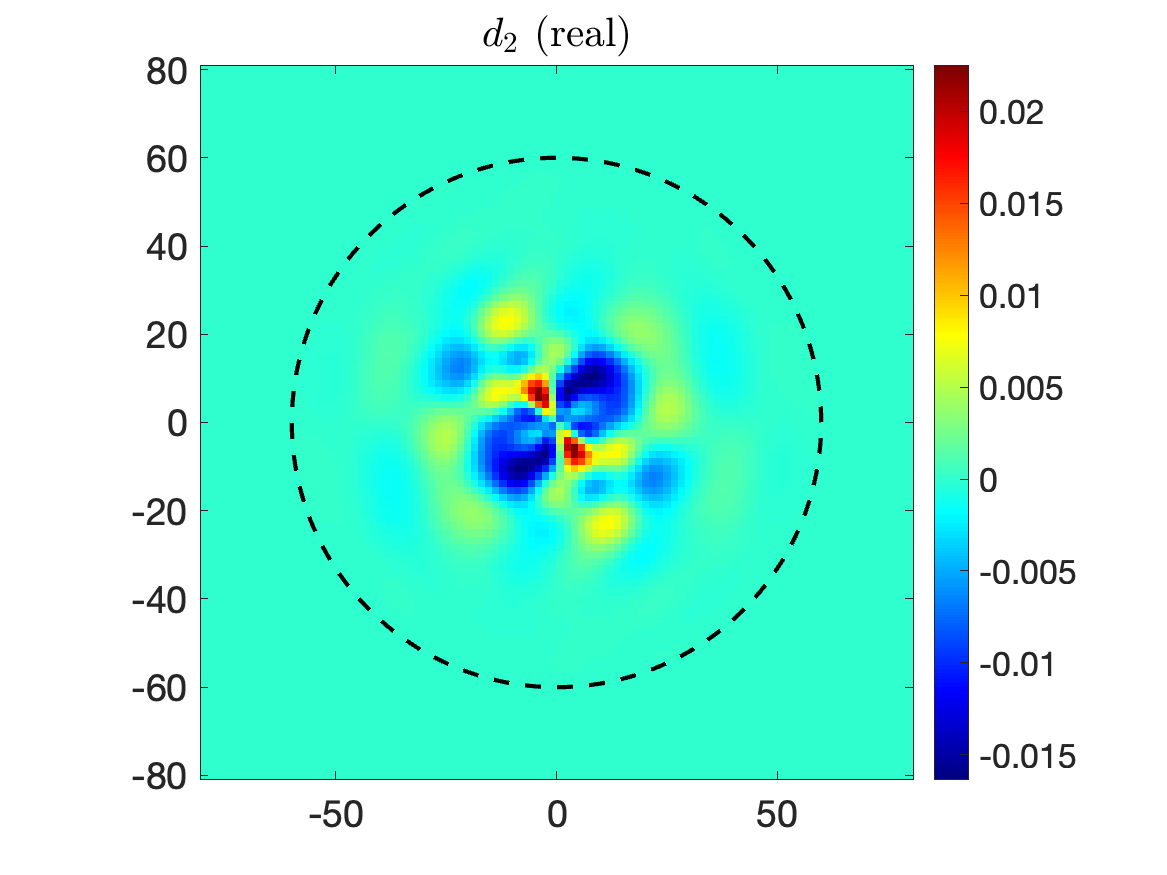} &
\includegraphics[width=0.3\textwidth]{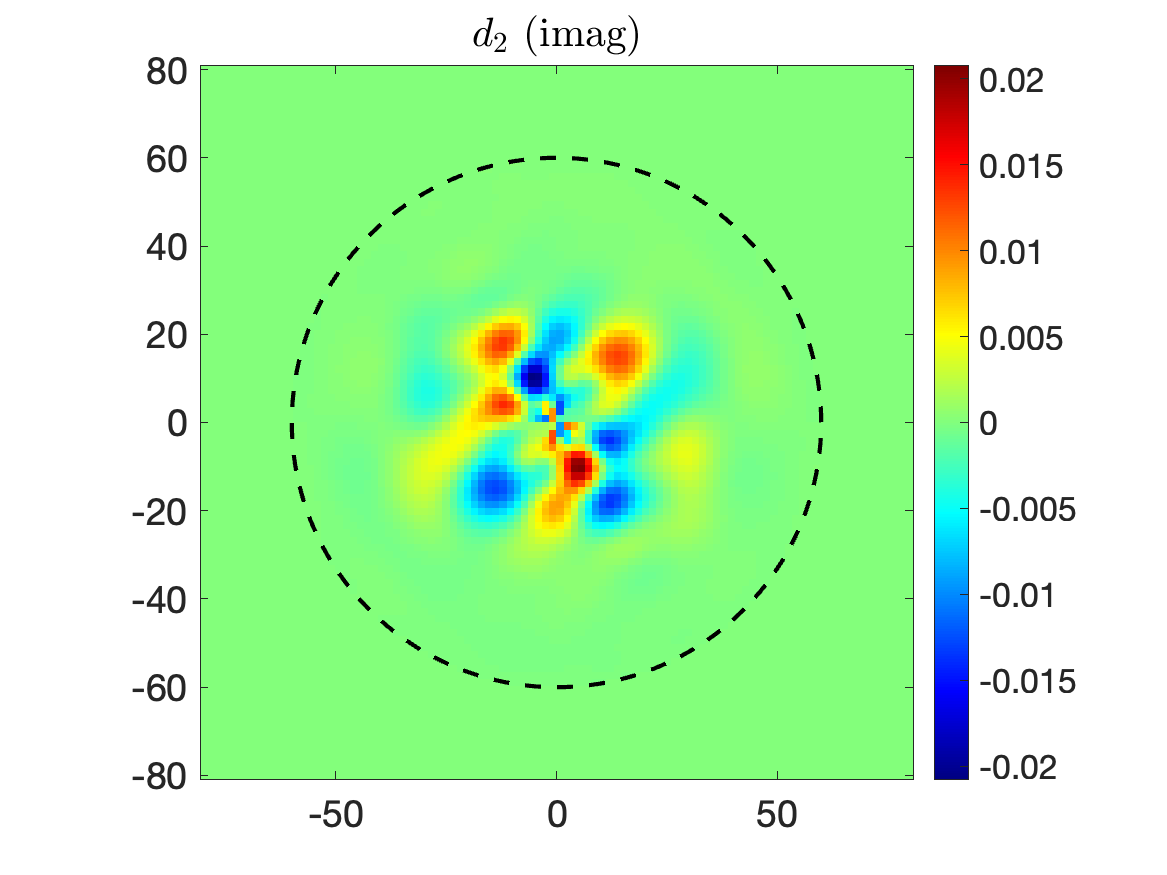} \\
\textbf{(2)} $c^{d}_{2}(x)$ &
\textbf{(2-a)} $d_{2}(\xi)$ (real) &
\textbf{(2-b)} $d_{2}(\xi)$ (imaginary) \\[3ex]
\includegraphics[width=0.33\textwidth]{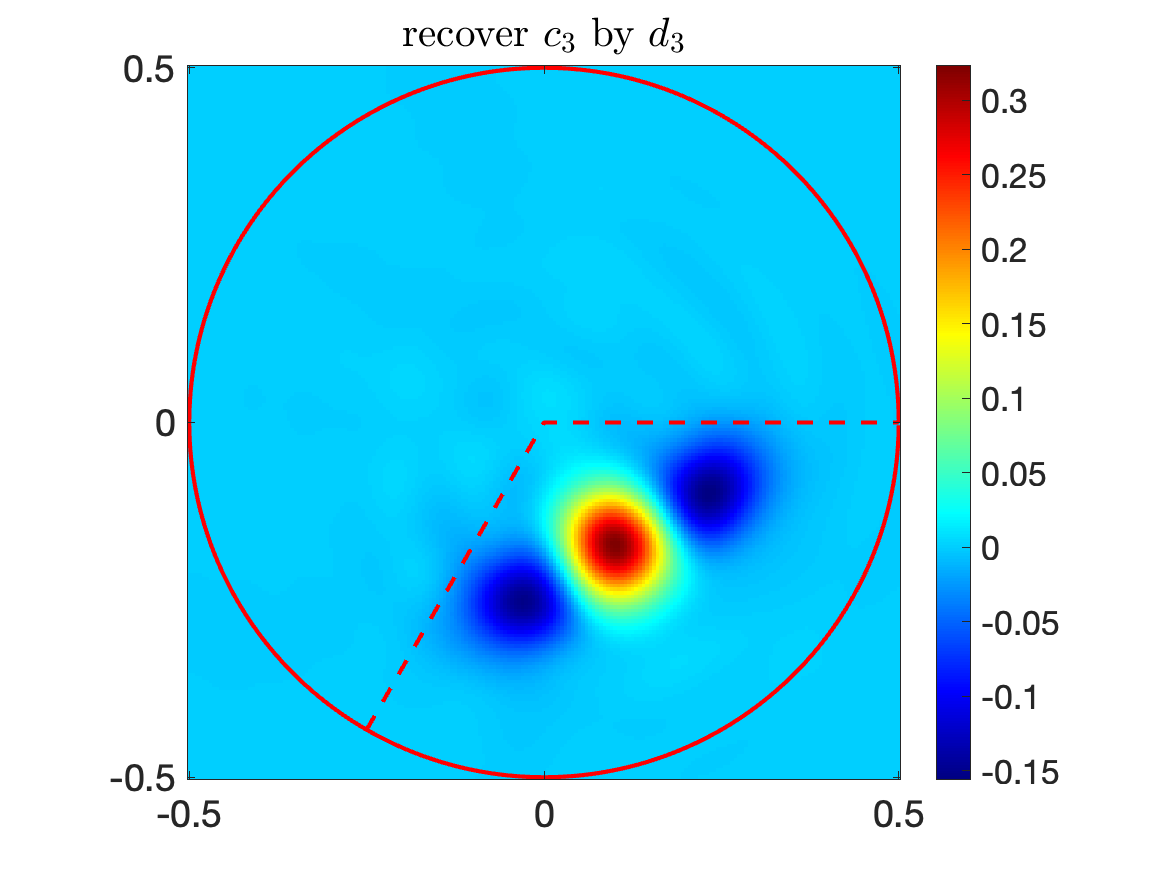} &
\includegraphics[width=0.3\textwidth]{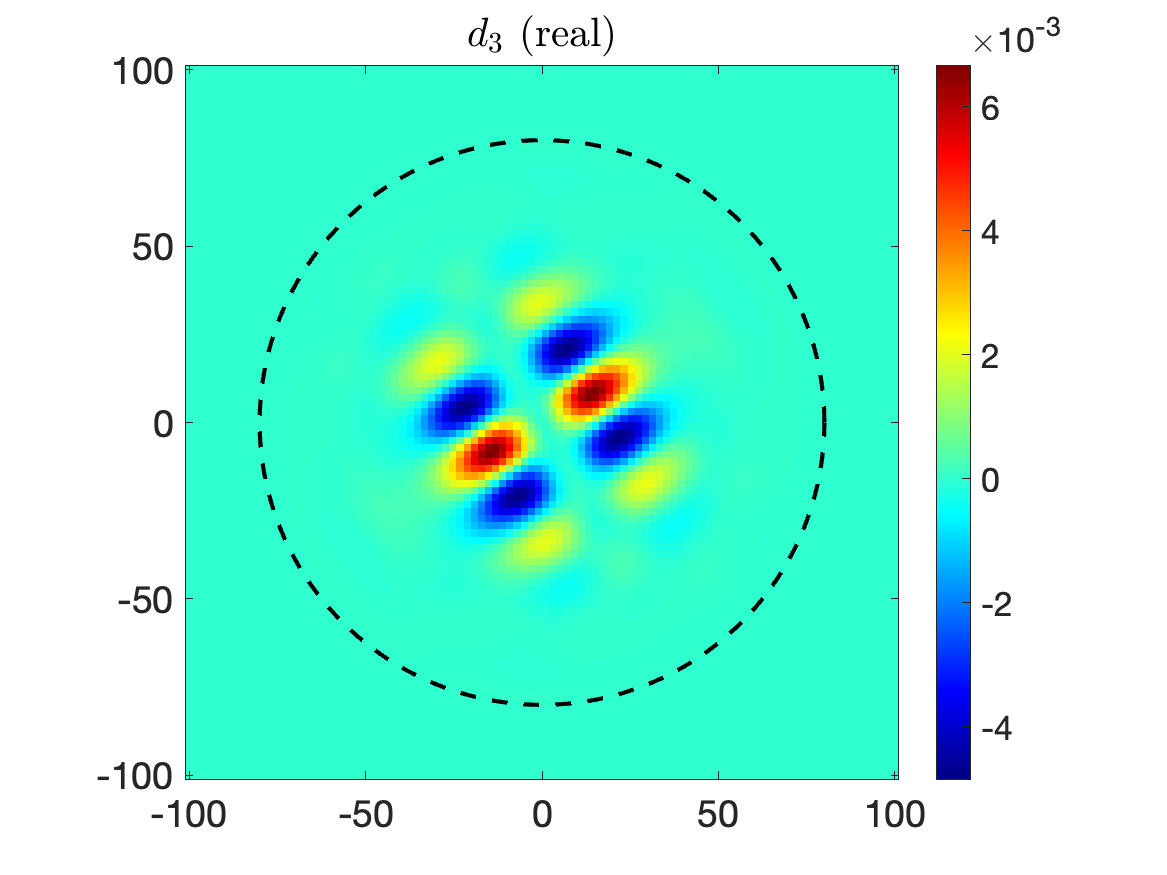} &
\includegraphics[width=0.3\textwidth]{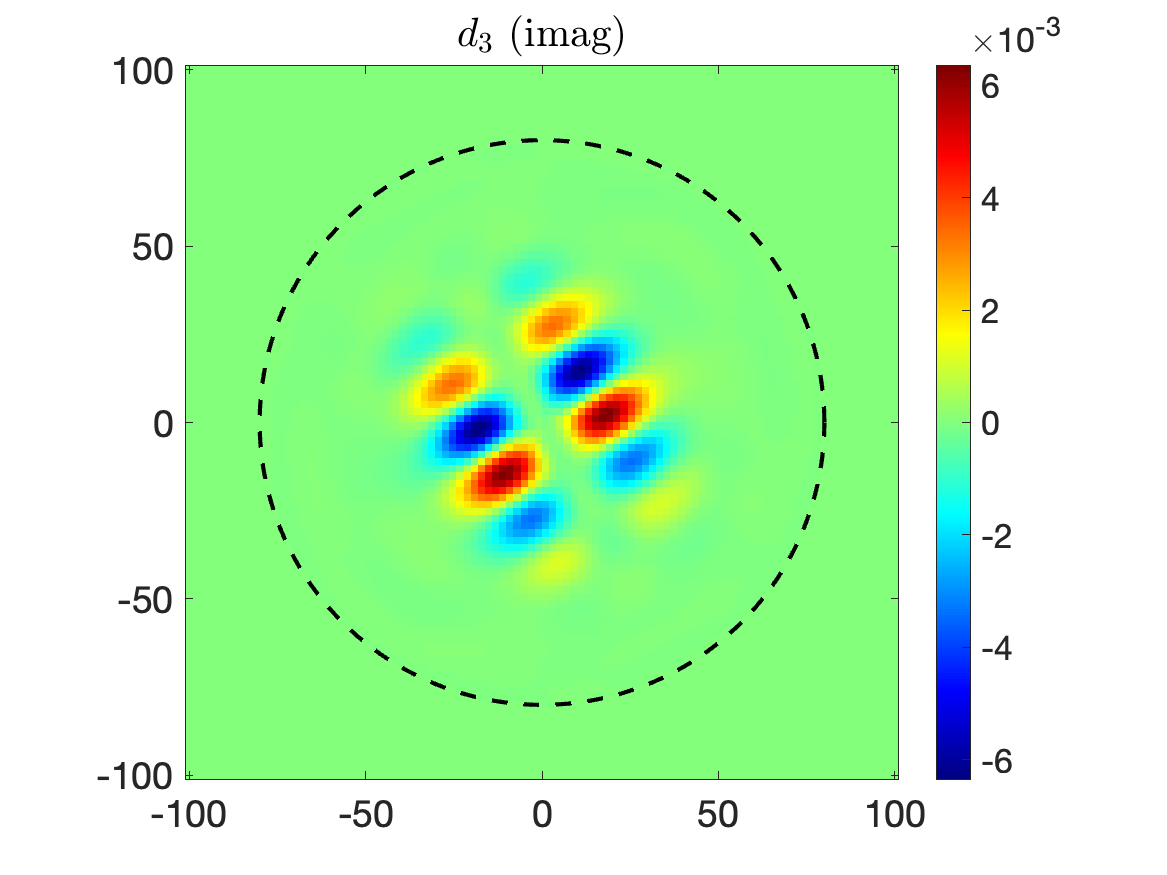} \\
\textbf{(3)} $c^{d}_{3}(x)$ &
\textbf{(3-a)} $d_{3}(\xi)$ (real) &
\textbf{(3-b)} $d_{3}(\xi)$ (imaginary) \\
\end{tabular}
\caption{%
\textsf{Left}: The reconstructed coefficients $c^{d}_{\ell}(x) = \mathcal{F}^{-1}[d_{\ell}](x)$ by using the inverse Fourier transform of the data $d_{\ell}(\xi)$ directly, $\ell = 1,2,3$.
\textsf{Middle \& Right}: The data $d_{\ell}(\xi)$ computed from measurements on the boundary $\partial\Omega$, $\ell = 1,2,3$.
\textbf{(1)} The reconstructed coefficient $c^{d}_{1}(x)$, \textbf{(1-a)} the real part and \textbf{(1-b)} the imaginary part of the data $d_{1}(\xi)$.
\textbf{(2)} The reconstructed coefficient $c^{d}_{2}(x)$, \textbf{(2-a)} the real part and \textbf{(2-b)} the imaginary part of the data $d_{2}(\xi)$.
\textbf{(3)} The reconstructed coefficient $c^{d}_{3}(x)$, \textbf{(3-a)} the real part and \textbf{(3-b)} the imaginary part of the data $d_{3}(\xi)$.
}
\label{figs:poly_m3_d_l}
\end{figure}

Similar to the previous example, the real and imaginary parts of the data $d_{1}(\xi)$, $d_{2}(\xi)$, and $d_{3}(\xi)$ are displayed in \textbf{Figure \ref{figs:poly_m3_d_l}}:
\textbf{Panel (1-a)} and \textbf{(1-b)} for the components of $d_{1}(\xi)$,
\textbf{Panel (2-a)} and \textbf{(2-b)} for the components of $d_{2}(\xi)$,
\textbf{Panel (3-a)} and \textbf{(3-b)} for the components of $d_{3}(\xi)$.
A comparison between \textbf{Figure \ref{figs:poly_m3_exact}} and \textbf{Figure \ref{figs:poly_m3_d_l}} also verifies that
\begin{align*}
\mathcal{F}[c_{3}] = d_{3},
\end{align*}
and
\begin{align*}
\mathcal{F}[c_{2}] \neq d_{2},
\qquad
\mathcal{F}[c_{1}] \neq d_{1},
\end{align*}
in agreement with \textbf{Lemma \ref{lmm:formulas}} for $m = 3$.

Thus, the reconstructed coefficient $c^{d}_{3}$ converges to the exact coefficient $c_{3}$; see \textbf{Figure \ref{figs:poly_m3_exact} (3)} and \textbf{Figure \ref{figs:poly_m3_d_l} (3)}.
But, both the reconstructed coefficients $c^{d}_{1}$ and $c^{d}_{2}$ are distorted due to interference from the higher-order coefficients $c_{2}$ and $c_{3}$; see \textbf{Figure \ref{figs:poly_m3_d_l} (1)} and \textbf{(2)}.

\begin{figure}[htbp]
\centering
\begin{tabular}{@{}ccc@{}}
\includegraphics[width=0.33\textwidth]{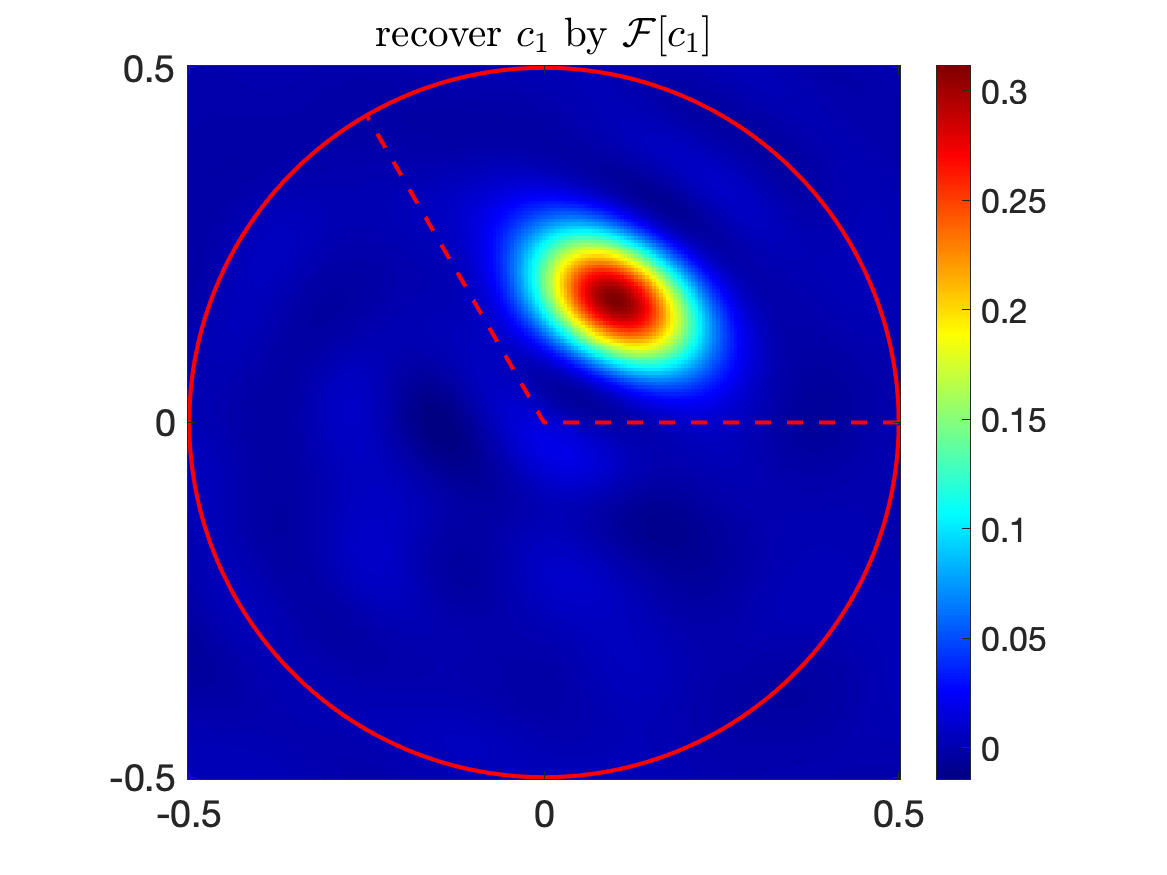} &
\includegraphics[width=0.3\textwidth]{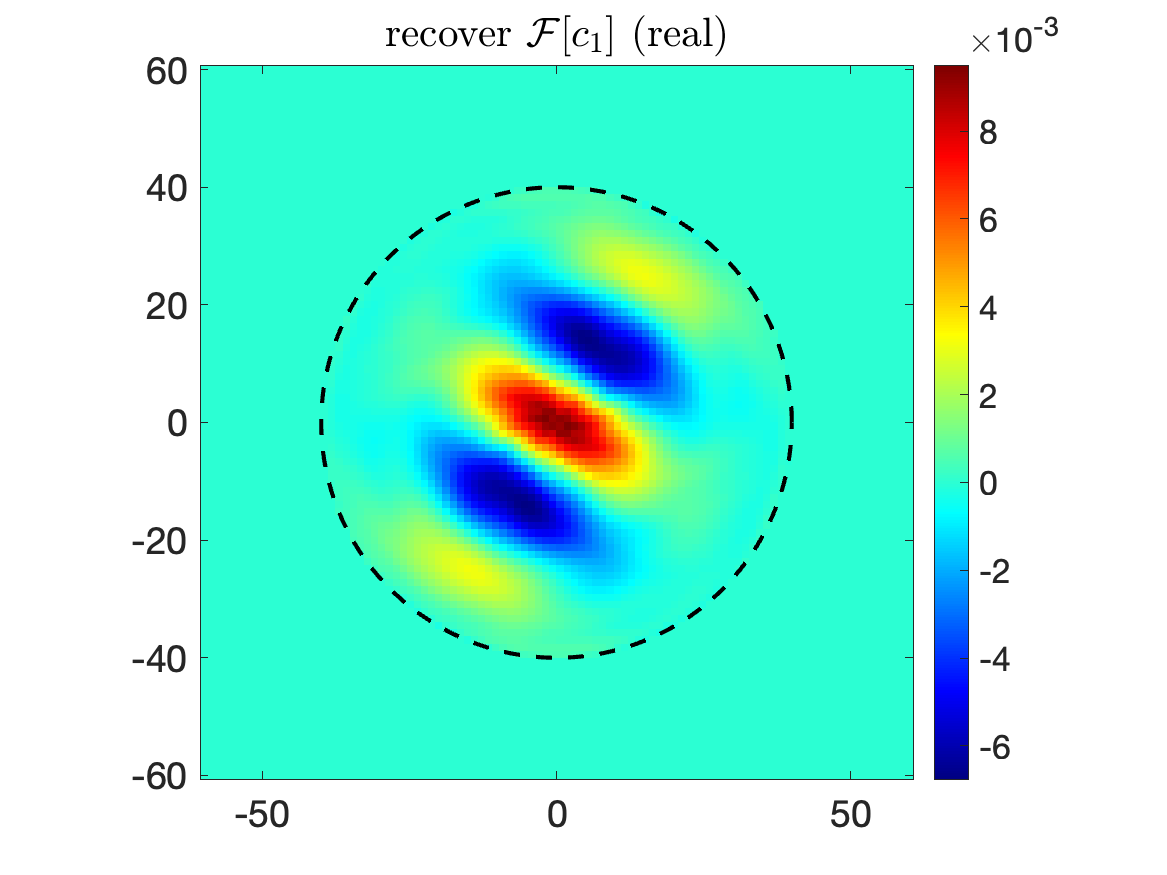} &
\includegraphics[width=0.3\textwidth]{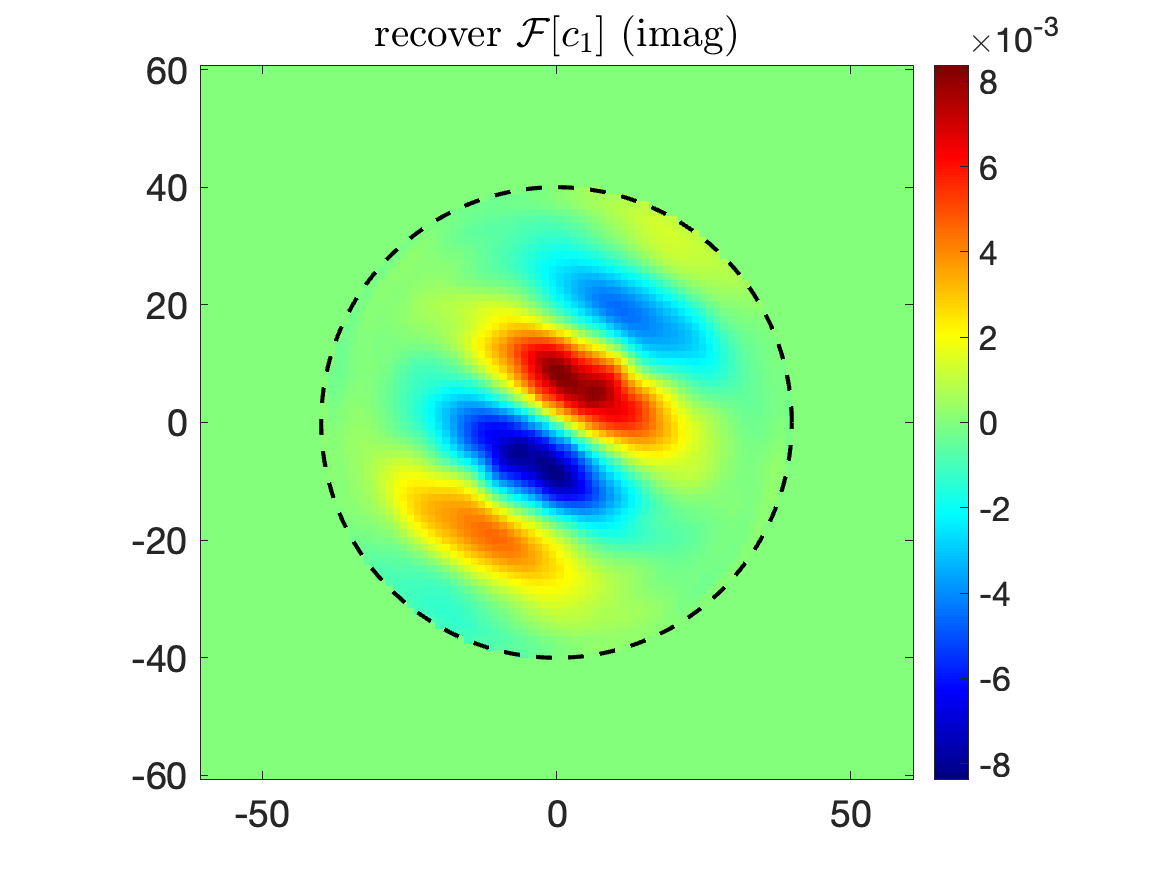} \\
\textbf{(1)} $\widetilde{c}_{1}(x)$ &
\textbf{(1-a)} $\widetilde{d}_{1}(\xi)$ (real) &
\textbf{(1-b)} $\widetilde{d}_{1}(\xi)$ (imaginary) \\[3ex]
\includegraphics[width=0.33\textwidth]{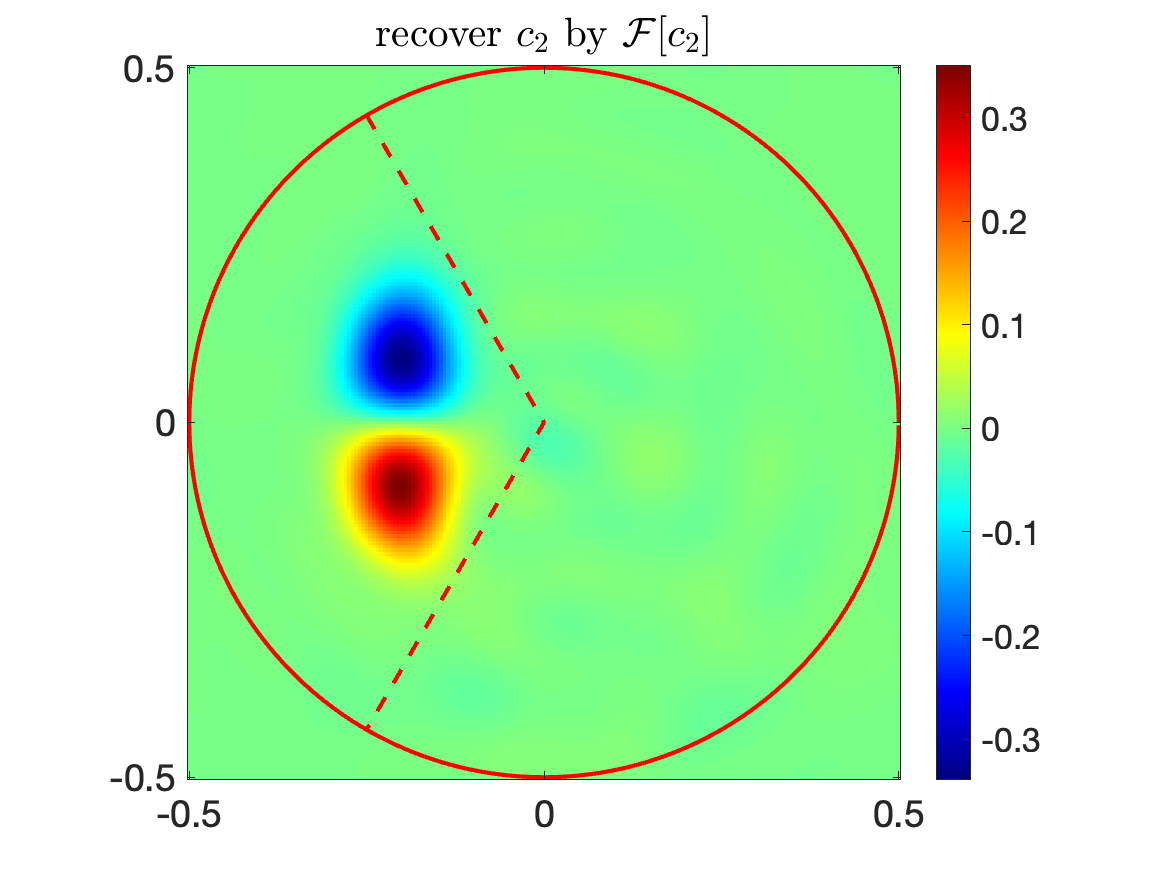} &
\includegraphics[width=0.3\textwidth]{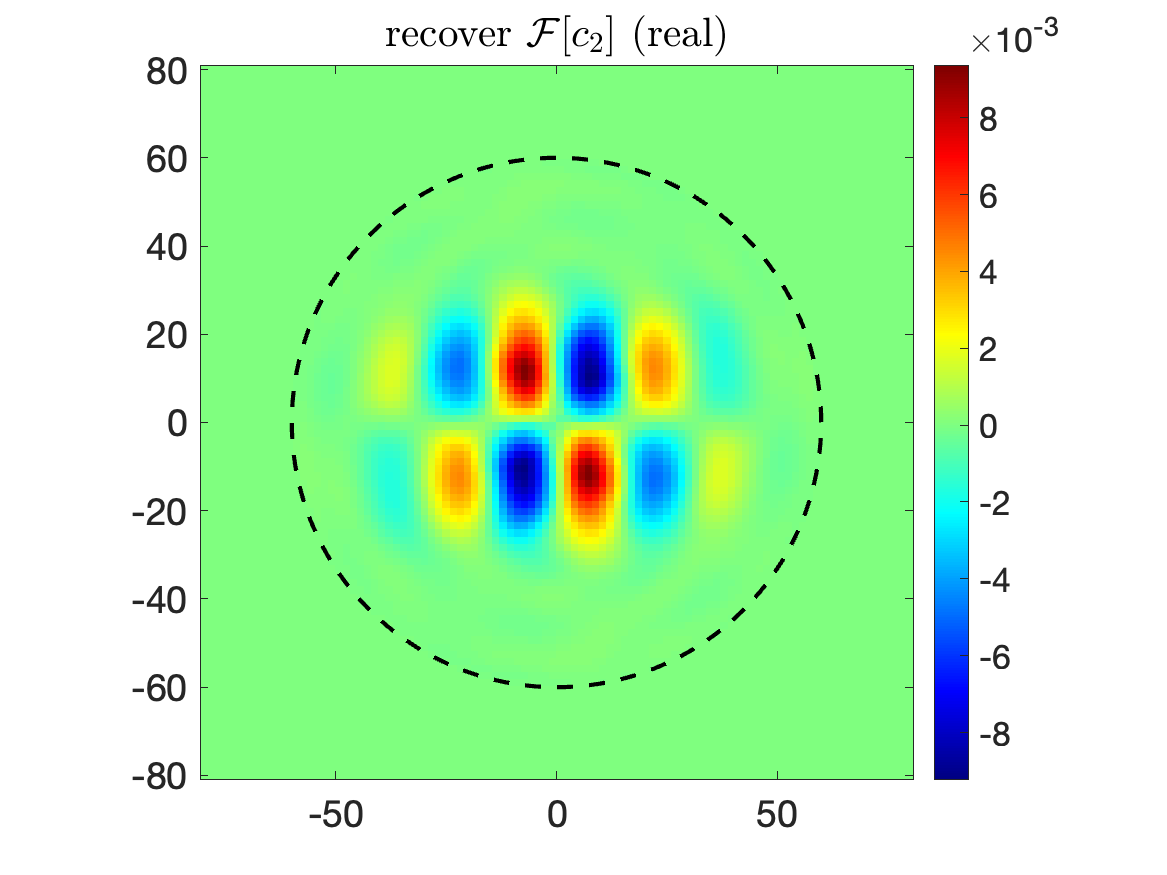} &
\includegraphics[width=0.3\textwidth]{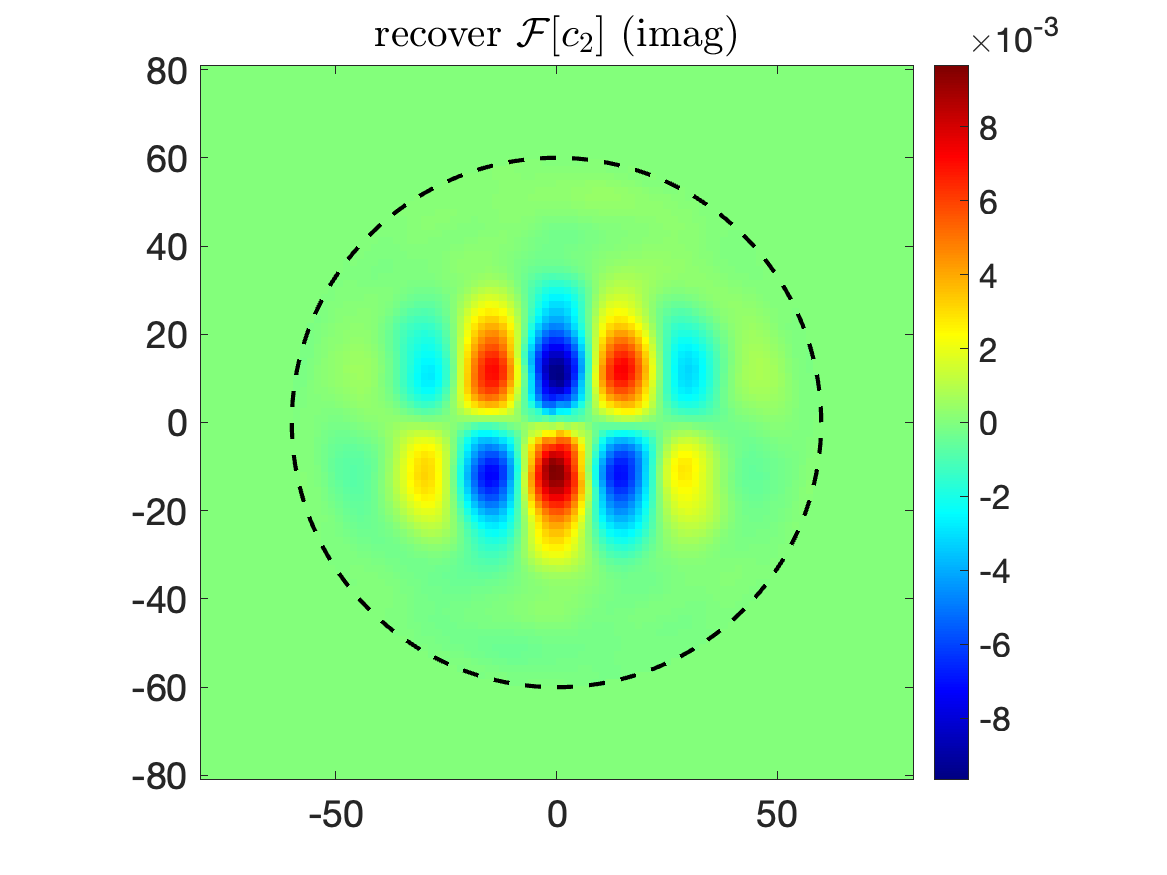} \\
\textbf{(2)} $\widetilde{c}_{2}(x)$ &
\textbf{(2-a)} $\widetilde{d}_{2}(\xi)$ (real) &
\textbf{(2-b)} $\widetilde{d}_{2}(\xi)$ (imaginary) \\
\end{tabular}
\caption{%
\textsf{Left}: The recovered coefficients $\widetilde{c}_{\ell}(x) = \mathcal{F}^{-1}[\widetilde{d}_{\ell}](x)$ by using the inverse transform of the modified data $\widetilde{d}_{1}(\xi)$, $\ell = 1,2$.
\textsf{Middle \& Right}: The modified data $\widetilde{d}_{\ell}(\xi)$ updated by using previously modified data $\big\{ \widetilde{d}_{\ell+a}(\xi) : a = 1,2,\dots,m-\ell \big\}$ and the data $d_{\ell}(\xi)$, $\ell = 1,2$.
\textbf{(1)} The recovered coefficient $\widetilde{c}_{1}(x)$, \textbf{(1-a)} the real part and \textbf{(1-b)} the imaginary part of the modified data $\widetilde{d}_{1}(\xi)$.
\textbf{(2)} The recovered coefficient $\widetilde{c}_{2}(x)$, \textbf{(2-a)} the real part and \textbf{(2-b)} the imaginary part of the modified data $\widetilde{d}_{2}(\xi)$.
}
\label{figs:poly_m3_Fc_l}
\end{figure}

In this example, for $m = 3$, we obtain the modified data $\widetilde{d}_{\ell}(\xi)$ ($\ell = 1,2,3$) by
\begin{align*}
\widetilde{d}_{3}(\xi) = d_{3}(\xi),
\end{align*}
and
\begin{align*}
\widetilde{d}_{2}(\xi) &= d_{2}(\xi) - \tfrac{3}{2} \widetilde{d}_{3}(\xi+\zeta_{2,1}) - \tfrac{3}{2} \widetilde{d}_{3}(\xi+\zeta_{2,2}), \\
\widetilde{d}_{1}(\xi) &= d_{1}(\xi) - \widetilde{d}_{3}(\xi+2\zeta_{1,1}) - \widetilde{d}_{2}(\xi+\zeta_{1,1}),
\end{align*}
respectively.
We also show the modified data $\widetilde{d}_{1}(\xi)$ and $\widetilde{d}_{2}(\xi)$ in \textbf{Figure \ref{figs:poly_m3_Fc_l} (1-a)--(1-b)} and \textbf{(2-a)--(2-b)}, and show the recovered coefficients $\widetilde{c}_{1}(x)$ and $\widetilde{c}_{2}(x)$ in \textbf{Figure \ref{figs:poly_m3_Fc_l} (1)} and \textbf{(2)}.
Again, we verify the equalities
\begin{align*}
\mathcal{F}[c_{2}] = \widetilde{d}_{2},
\qquad
\mathcal{F}[c_{1}] = \widetilde{d}_{1},
\end{align*}
in agreement with \textbf{Lemma \ref{lmm:formulas}}: \textbf{Case $c_{\ell}$ with $0 < \ell < m$}.

%


%
%
%

\bibliographystyle{siamplain}
\bibliography{PolyInv}

\end{document}